\documentclass[draftcls,onecolumn]{IEEEtran}
\newcommand{\argmin}{\operatornamewithlimits{arg~min}}

\newtheorem{fact}{Fact}

\usepackage{ifpdf}

\usepackage{cite}

\ifCLASSINFOpdf
 \usepackage[pdftex]{graphicx}

\else

   \usepackage[dvips]{graphicx}
  
\fi

\makeatletter
\newcommand\footnoteref[1]{\protected@xdef\@thefnmark{\ref{#1}}\@footnotemark}
\makeatother

\usepackage[cmex10]{amsmath}
\usepackage{amsfonts}

\usepackage{subfigure}
\usepackage{multirow}
\usepackage{multicol}
\usepackage{algorithm}
\usepackage{algorithmic}
\usepackage{wrapfig}
\usepackage{bbm}

\usepackage{url}

\begin{document}
%
\title{Fast Singular Value Shrinkage with\\ Chebyshev Polynomial Approximation\\ Based on Signal Sparsity}
%
%
%

\author{Masaki~Onuki,~\IEEEmembership{Student~Member,~IEEE,}
        Shunsuke~Ono,~\IEEEmembership{Member,~IEEE,}
        Keiichiro~Shirai,~\IEEEmembership{Member,~IEEE,}
        and~Yuichi~Tanaka,~\IEEEmembership{Member,~IEEE.}
\thanks{M. Onuki and Y. Tanaka are with the Grad. School of BASE, Tokyo Univ. of Agri. and Tech., Koganei, Tokyo, 184-8588 Japan (e-mail: masaki.o@msp-lab.org; ytnk@cc.tuat.ac.jp).

S. Ono is with the Lab. for Future Interdisciplinary Res. of Sci. and Tech. (FIRST), Tokyo Inst. of Tech., Midori, Kanagawa, 226-8503 Japan (e-mail: ono@isl.titech.ac.jp).

K. Shirai is with the Dept. of Elec. and Compt. Eng., Shinshu Univ., Wakasato, Nagano, 380-8553 Japan (keiichi@shinshu-u.ac.jp).}
\thanks{Manuscript received ; revised.}}

%
%

\markboth{}%
{Onuki \MakeLowercase{\textit{et al.}}}

\maketitle

%
%
\begin{abstract}
We propose an approximation method for thresholding of singular values using Chebyshev polynomial approximation (CPA).
Many signal processing problems require iterative application of singular value decomposition (SVD) for minimizing the rank of a given data matrix with other cost functions and/or constraints, which is called matrix rank minimization.
In matrix rank minimization, singular values of a matrix are shrunk by hard-thresholding, soft-thresholding, or weighted soft-thresholding.
However, the computational cost of SVD is generally too expensive to handle high dimensional signals such as images; hence, in this case, matrix rank minimization requires enormous computation time.
In this paper, we leverage CPA to (approximately) manipulate singular values without computing singular values and vectors.
The thresholding of singular values is expressed by a multiplication of certain matrices, which is derived from a characteristic of CPA.
The multiplication is also efficiently computed using the sparsity of signals.
As a result, the computational cost is significantly reduced.
Experimental results suggest the effectiveness of our method through several image processing applications based on matrix rank minimization with nuclear norm relaxation in terms of computation time and approximation precision.
\end{abstract}

\begin{IEEEkeywords}
Chebyshev polynomial approximation, nuclear norm relaxation, singular value thresholding
\end{IEEEkeywords}

\IEEEpeerreviewmaketitle

%
%
\section{Introduction}

\begin{figure}[htb]
\small
\begin{tabular}{c}
\begin{minipage}[h]{0.95 \linewidth}
  \centering
  \centerline{\includegraphics[width=7.5cm]{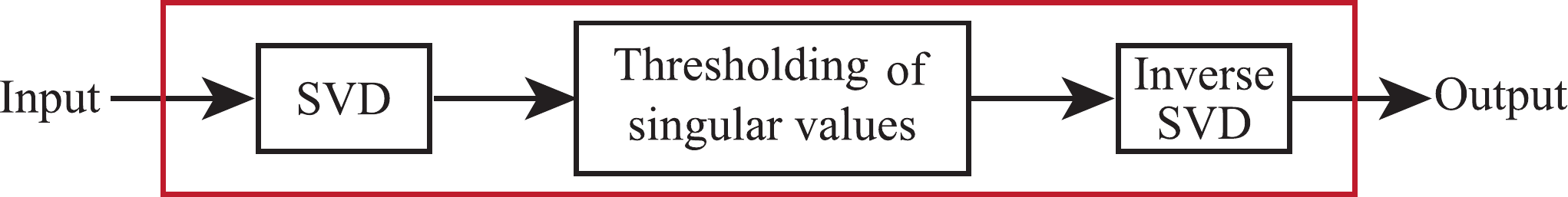}}
  \centerline{(a) Explicit singular value shrinkage}\medskip
\end{minipage}
\\  \\
\begin{minipage}[h]{0.95 \linewidth}
  \centering
  \centerline{\includegraphics[width=6.7cm]{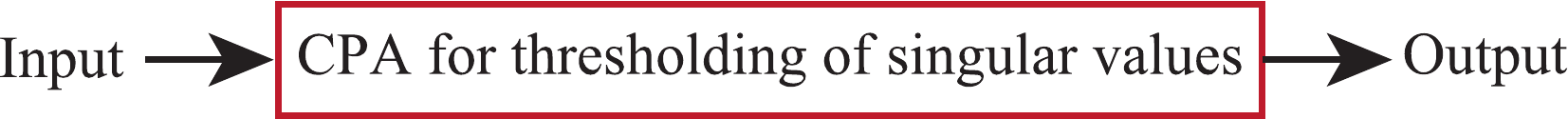}}
  \centerline{(b) CPA-based singular value shrinkage}\medskip
\end{minipage}
\\ \\
\begin{minipage}[h]{0.95 \linewidth}
  \centering
  \centerline{\includegraphics[width=7.2cm]{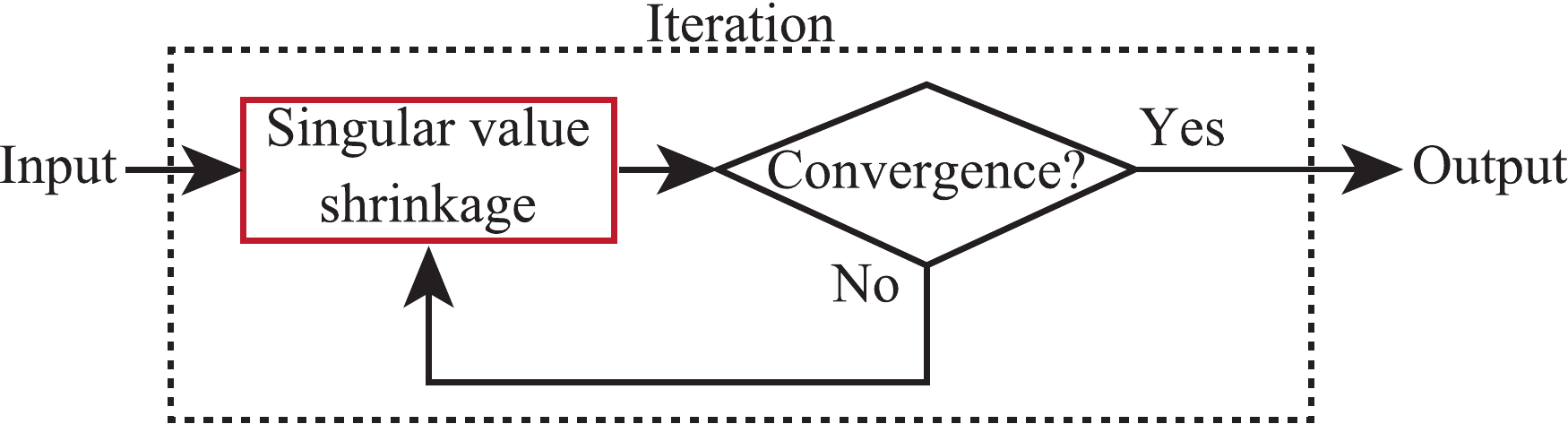}}
  \centerline{(c) Matrix rank minimization}\medskip
\end{minipage}\\
\end{tabular}
\caption{Descriptions of singular value shrinkage and matrix rank minimization with singular value shrinkage.
Our paper is focused on red parts.}
\label{fig:methods_shrinkage}
\end{figure}

\IEEEPARstart{T}{he} low-rank structure inherent in various signals has been widely exploited in many signal processing applications, such as matrix and tensor completion  \cite{matrix_tensor_completion1,matrix_tensor_completion2,matrix_tensor_completion3,matrix_tensor_completion4}, image decomposition \cite{ADMM,image_decomposition}, photometric stereo \cite{photometoric_stereo1,photometoric_stereo2}, image alignment \cite{image_alignment1,image_alignment2}, colorization \cite{colorization1}, inpainting \cite{repaire_texture1,repaire_texture2}, background modeling \cite{background_mod1,background_mod2,background_mod3,background_mod4,background_mod5}, color artifact removal \cite{color_artifact}, cognitive radio \cite{cognitive_radio}, and voice separation \cite{voice_separation}.
In such applications, the low-rank structure is incorporated into a minimization problem involving the rank function or its continuous relaxation.
The problem is solved using some iterative algorithms, with which the thresholding of singular values is usually required at each iteration.
We refer to this methodology as \textit{matrix rank minimization}.

There are two representative approaches of matrix rank minimization.
One is the exact method.
It is an ideal formulation, but the resulting problem is very difficult to solve due to the non-convexity and combinatorial nature of the rank function.
The other is the nuclear norm relaxation \cite{nuclear_norm}.
Since the nuclear norm, the sum of the singular values of a matrix, is the tightest convex relaxation of the rank function, we can efficiently solve the resulting problem via convex optimization techniques.
Weighted nuclear norm relaxation \cite{weightedSoft2,weightedSoft} has recently been proposed as a non-convex but continuous approximation of the rank function.

Essentially, both of the above methods require the thresholding of singular values, which we call \textit{singular value shrinkage}, at each iteration of certain optimization methods (Fig.~\ref{fig:methods_shrinkage}).
That is, most methods for matrix rank minimization must carry out singular value decomposition (SVD) many times.
This is a serious problem in terms of computational cost when we handle large matrices, even with high-spec computers.

Several methods have been proposed to tackle this issue \cite{SVD_echon2, SVD_echon3,Fast_singular_thresh}.
The basic concept of \cite{SVD_echon2, SVD_echon3} is to approximately compute partial singular values and/or vectors.
These methods can drastically reduce the computation time of singular value shrinkage but would not be suitable for the matrix rank minimization.
Since the number of singular values above a threshold is not identified without the full decomposition, many singular values above a threshold are reduced to zero in each iteration.
As a result, large approximation errors are produced, which results in an unstable convergence in the matrix rank minimization.
With the other method \cite{Fast_singular_thresh}, singular value shrinkage is carried out by computing neither singular values nor vectors, but the reduction in the computation time is still limited.
This is because the method requires a complete orthogonal decomposition \cite{SVD2} and the calculation of the inverse of a large matrix.
It also leads to large approximation errors, i.e., an unstable convergence in the matrix rank minimization.
We consider a method similar to that by Cai and Osher \cite{Fast_singular_thresh}: We only need a ``processed'' matrix with thresholded singular values.

In this paper, we propose a fast singular value shrinkage method for reducing the computational cost in the matrix rank minimization of high dimensional matrices.
Note that the proposed method computes neither singular values nor vectors during the process of singular value shrinkage, similar to the method by Cai and Osher \cite{Fast_singular_thresh}.
Furthermore, our method maintains computational precision to lead matrix rank minimization algorithms to stable convergence.
The two key tools of our method are described as follows.
\begin{itemize}
\item \textbf{Chebyshev polynomial approximation (CPA) \cite{cheby1,cheby2,cheby3}}:
This tool is often used for designing filters in signal processing \cite{kaiser1,okd_method_CPA} and is a key tool for reducing computational cost.
The applications of CPA have been studied by Saad et~al.\cite{polynomial_filtering3,polynomial_filtering1,polynomial_filtering2,CJPA3}.
With the applications by Saad et~al., CPA is used to calculate a \textit{vector} after being transformed by a matrix with singular value shrinkage.
That is, it requires the iterative multiplications of a matrix and vector to derive the Chebyshev polynomials.
The concept of the applications has recently been used for improving the performance of image filtering methods such as bilateral filter, non-local means, and BM3D \cite{eig_cheby1,eig_cheby2,eig_cheby3}.
In contrast, we propose a method to obtain a \textit{matrix} whose singular values are processed by using CPA.
Since CPA results in truncation errors, such as \textit{ripples} in the lower-order approximations, we also investigate the designs of thresholding functions and appropriate approximation order for reducing approximation errors.
\item \textbf{Sparsity of signals}:
By using CPA, our method can represent singular value shrinkage as a multiplication of matrices.
Since the multiplication can be computed efficiently when the matrices are sparse, our method exploits the inherent sparsity of signals in their frequency domain for further acceleration.
\end{itemize}

Since matrix rank minimization plays a central role in various signal processing tasks, our method offers many promising applications.
For this study, we validated the proposed method by using two image processing applications: \textit{image inpainting} \cite{repaire_texture1} and \textit{background modeling}\cite{background_mod1,background_mod2,background_mod3,background_mod4,background_mod5}.
In these applications, target problems are formulated as convex optimization problems involving the nuclear norm so that they can be efficiently solved using the \textit{alternating direction method of multiplier} (ADMM) \cite{ADMM2} with our method.

Although the ADMM is widely known as a robust method for computation errors in each iteration, optimization methods (including the ADMM) with the other fast singular value shrinkage methods \cite{SVD_echon2,SVD_echon3,Fast_singular_thresh} do not converge well due to their large approximation errors.
In contrast, our CPA-based singular value shrinkage method leads optimization methods to stable convergence.
We validated this advantage experimentally by comparing our method with the other fast singular value shrinkage methods in several image processing tasks and a synthetic data.

The preliminary version of this study, without using signal sparsity, analysis of our method, and new applications, has previously been published \cite{im_col_cheby}.

This paper is organized as follows.
Section \ref{sec:notations_and_mathematical_tools} defines notations and preliminaries.
We discuss our CPA-based singular value shrinkage method, which is the main contribution in this paper, in Section \ref{sec:singular_filter}.
We discuss an approximation order of CPA for reducing the size of approximation errors in Section \ref{sec:chebyshev_filter_design} and verification of our method through applications in Section \ref{sec:applications}.
Finally, we conclude the paper in Section \ref{sec:conclusion}.

%
%
\section{Notations and Preliminaries}
\label{sec:notations_and_mathematical_tools}

%
%

\subsection{Notations}

Bold-face capital and small letters indicate a matrix and a vector, respectively.
Superscript $\cdot^\top$ is the transpose of a matrix and a vector, and superscript $\cdot^{-1}$ is the inverse of a non-singular matrix.
The matrices $\mathbf{Id}$ and $\mathbf{O}$ are the identity matrix and null matrix, respectively.
The vector $\mathbbm{1}_n:=[\underbrace{1,\ldots,1}_{n}]^\top$.
The $\ell_p$ norm for $p\ge1$ is defined as $\| \mathbf{x} \|_p \!:=\! (\sum^n_{i=1} | x_i |^p)^{\frac{1}{p}}~(\forall \mathbf{x}\!\in\! \mathbb{R}^n)$.
We also use CPA as follows.

%
%
\subsection{Chebyshev Polynomial Approximation}

Let $h(x)$ and $\widehat{h}(x)$ be a real-valued function defined on the interval $x\!\in\! [-1,1]$ and its approximated function by using CPA, respectively.
Chebyshev polynomial approximation \cite{cheby1,cheby2,cheby3} gives an approximate solution of $h(x)$ by using the \textit{truncated Chebyshev series}:
\begin{equation}
\widehat{h}(x):=\frac{1}{2}c_0 + \sum^{\alpha-1}_{k=1} c_k \, \psi_k(x),
\label{chebychev_series}
\end{equation}
where $c_k$ and $\alpha$ denote a Chebyshev coefficient (described later) and an approximation order, respectively.
Additionally, $\psi_k(\cdot)$ denotes the $k$-th order \textit{Chebyshev polynomials} of the first kind, defined as
\begin{equation}
\psi_k(x):=\cos \bigl( k \arccos (x) \bigr).
\label{chebychev_polynomial}
\end{equation}
It can also be computed using the stable recurrence relation:
\begin{equation}
\begin{split}
&\psi_k(x)=2x \, \psi_{k-1}(x)-\psi_{k-2}(x),\\
&\psi_0(x)=1,\quad \psi_1(x)=x.
\end{split}
\label{chebyche2}
\end{equation}
The initial condition is defined as $\psi_0(x)$ and $\psi_1(x)$.
Since the polynomials consist of cosine functions, the value of $\psi_k(x)$ is bounded between $-1$ and $1$ for $x \!\in\! [-1,1]$.
By using $\psi_k(x)$ and the orthogonality of the cosine function, $c_k$ is calculated as
\begin{equation}
c_k := \frac{2}{\alpha}\sum^\alpha_{l=1} \cos( k \theta(l) ) \, h( \cos \theta(l)  ),
\label{chebychev_coeff_discrete}
\end{equation}
where $\displaystyle\theta(l) \!:=\!\frac{\pi(l-\frac{1}{2})}{\alpha}$.

%
%

\section{Singular Value Shrinkage using Chebyshev Polynomial Approximation by Exploiting Sparsity}
\label{sec:singular_filter}

We discuss singular value shrinkage using CPA.
First, the CPA of a matrix form, which can approximately shrink the eigenvalues of a matrix (\textit{eigenvalue shrinkage}), is indicated then extended to the singular one.

%
%
\subsection{Chebyshev Polynomial Approximation for Matrix}

Let $\mathbf{A} \!\in\! \mathbb{R}^{n\times n}$ be a full rank matrix and $\mathbf{A} \!=\! \mathbf{P} \mathbf{\Lambda}_\mathrm{A} \mathbf{P}^{-1}$ be its eigendecomposition (EVD), where $\mathbf{P} \!\in\! \mathbb{R}^{n\times n}$ is the matrix composed of eigenvectors and $\mathbf{\Lambda}_\mathrm{A}\!=\!\mathrm{diag}(\lambda^\mathrm{A}_1, \ldots , \lambda^\mathrm{A}_i, \ldots, \lambda^\mathrm{A}_n)$ is the diagonal matrix with the corresponding eigenvalues.
We assume that the eigenvalues are bounded between $0$ and $\lambda^\mathrm{A}_\mathrm{max}$, where $\lambda^\mathrm{A}_\mathrm{max} \!>\! 1$.
Hence, the eigenvalues of $\mathbf{A}$ are shrunk as
\begin{equation}
\mathcal{H}(\mathbf{A}):=
\mathbf{P}\,
\mathrm{diag}\bigl( h(\lambda^\mathrm{A}_1), \ldots , h(\lambda^\mathrm{A}_n)\bigr)
\mathbf{P}^{-1},
\label{eigen_shrink_normal}
\end{equation}
where $\mathcal{H}(\cdot)$ is the eigenvalue shrinkage function, and $h(x)$ is the filter kernel defined in $x\!\in\! [0, \lambda^\mathrm{A}_\mathrm{max}]$.
In this subsection, we consider the approximated solution of \eqref{eigen_shrink_normal} using the CPA.

The CPA of the matrix form \cite{cheby2,polynomial_filtering3,polynomial_filtering2,CJPA3} gives an approximated solution of the eigenvalue shrinkage function $\mathcal{H}(\cdot)$ by using truncated Chebyshev series as
\begin{equation}
\widehat{\mathcal{H}}(\mathbf{A}):=\frac{1}{2}\widehat{c}_0 \, \mathbf{Id} + \sum^{\alpha-1}_{k=1} \widehat{c}_k \mathit{\Psi}_k(\widehat{\mathbf{A}}),
\label{chebychev_series_matrix}
\end{equation}
where $\widehat{c}_k$ and $ \mathit{\Psi}_k(\widehat{\mathbf{A}})$ are Chebyshev coefficients and Chebyshev polynomials, respectively, which are defined later.
Additionally, $\widehat{\mathbf{A}}$ is the eigenvalue-shifted matrix given by
\begin{equation}
\widehat{\mathbf{A}} := \frac{2}{\lambda^\mathrm{A}_\mathrm{max}}\mathbf{A}-\mathbf{Id},
\label{transform_lamB}
\end{equation}
whose eigenvalues are obviously within $[-1,1]$.
Thanks to \eqref{transform_lamB}, the $k$-th order Chebyshev polynomial of $\widehat{\mathbf{A}}$ is computed as 
\begin{align}
\mathit{\Psi}_k(\widehat{\mathbf{A}}) &=\mathit{\Psi}_k\left(\frac{2}{\lambda^\mathrm{A}_\mathrm{max}}\mathbf{A}-\mathbf{Id}\right)\nonumber \\
&=\mathbf{P}\mathit{\Psi}_k\left(\frac{2}{\lambda^\mathrm{A}_\mathrm{max}}\mathbf{\Lambda}_\mathrm{A}-\mathbf{Id}\right)\mathbf{P}^{-1}\nonumber \\
&=\mathbf{P}\,\mathrm{diag}(\cos k\theta_1,\ldots ,\cos k\theta_n)\mathbf{P}^{-1},
\label{transform_lamA3}
\end{align}
where $\theta_i := \arccos \left(\frac{2}{\lambda^{\mathrm{A}}_\mathrm{max}}\lambda^{\mathrm{A}}_i-1\right)$.
Similarly to \eqref{chebyche2}, the Chebyshev polynomials are obtained using the recurrence relation:
\begin{equation}
\begin{split}
&\mathit{\Psi}_k(\widehat{\mathbf{A}})=2\widehat{\mathbf{A}} \, \mathit{\Psi}_{k-1}(\widehat{\mathbf{A}})-\mathit{\Psi}_{k-2}(\widehat{\mathbf{A}}),\\
&\mathit{\Psi}_0(\widehat{\mathbf{A}})=\mathbf{Id},\quad \mathit{\Psi}_1(\widehat{\mathbf{A}})=\widehat{\mathbf{A}}.
\end{split}
\label{cheby_matrix_eigen_saiki3}
\end{equation}
Recall that $\mathit{\Psi}_k(\widehat{\mathbf{A}})$ is defined only in the interval $[-1,1]$.
Therefore, the range of the filter kernel is modified by deriving $\widehat{c}_k$ as
\begin{equation}
\widehat{c}_k = \frac{2}{\alpha}\sum^\alpha_{l=1} \cos( k\theta(l)  ) \, h\!\hspace{0.4mm}\Bigl( \frac{\lambda^\mathrm{A}_\mathrm{max}}{2}( \cos \theta(l)  +1 ) \Bigr).
\label{chebychev_coeff_shifted}
\end{equation}
The term $h\left(\lambda^\mathrm{A}_\mathrm{max}/2 \, (\cos \theta(l)  + 1 ) \right)$ returns the shifted range back to the original range $[0,\lambda^\mathrm{A}_\mathrm{max}]$.
From \eqref{chebychev_coeff_shifted}, $\widehat{\mathcal{H}}(\mathbf{A})$ can also be represented using $\widehat{h}(\lambda^\mathrm{A}_i)$ as
\begin{equation}
\widehat{\mathcal{H}}(\mathbf{A})=
\mathbf{P}\,
\mathrm{diag}\bigl( \widehat{h}(\lambda^\mathrm{A}_1), \ldots, \widehat{h}(\lambda^\mathrm{A}_n) \bigr)
\mathbf{P}^{-1}.
\label{chebychev_series_matrix_shift2}
\end{equation}
The function $\widehat{\mathcal{H}}(\cdot)$, which is referred to as the CPA-based eigenvalue shrinkage function, results in approximate eigenvalue shrinkage.
The CPA-based eigenvalue shrinkage actually computes neither eigenvalues nor vectors thanks to the recurrence relation \eqref{cheby_matrix_eigen_saiki3}.

%
%
\subsection{CPA-based Singular Value Shrinkage}
Let $\mathbf{B}\!\in\! \mathbb{R}^{m\times n}~(m\!>\!n)$ be a rectangular matrix and $\mathbf{B}\!=\!\mathbf{U}\mathbf{\Sigma}\mathbf{V}^\top$ be its singular value decomposition, where $\mathbf{U} \!\in\! \mathbb{R}^{m\times m}$ and $\mathbf{V}\!\in\! \mathbb{R}^{n\times n}$ are orthogonal matrices.
The $\mathbf{\Sigma}\!\in\! \mathbb{R}^{m\times n}$ is the singular value matrix represented as
\begin{equation}
\mathbf{\Sigma}=
\begin{bmatrix}
 	\sigma_1 &&\mathbf{O}\\
	 &\ddots&\\
	 &&\sigma_n\\
	 \mathbf{O}&&\\
	\end{bmatrix}.
\label{sin_val_matrix_A}
\end{equation}
Without loss of generality, we can assume $\sigma_1\geq \ldots \geq \sigma_n$.
The singular values of $\mathbf{B}$ are shrunk with the singular value shrinkage function $\mathcal{G}(\cdot)$ as
\begin{equation}
\mathcal{G}(\mathbf{B}):=\mathbf{U}
\begin{bmatrix}
 	g(\sigma_1) & &\mathbf{O}\\
	 &\ddots&\\
	 &&g(\sigma_n)\\
	 \mathbf{O}&&\\
	\end{bmatrix}
\mathbf{V}^\top,
\label{sing_shrink_exact}
\end{equation}
where $g(\cdot)$ is an arbitrary function.

The eigenvalue shrinkage in \eqref{eigen_shrink_normal} can be extended to $\mathcal{G}(\mathbf{B})$ in \eqref{sing_shrink_exact} as \cite{polynomial_filtering1}
\begin{equation}
\mathcal{G}(\mathbf{B})=\mathbf{B} \, \mathcal{H}(\mathbf{B}^{\!\top}\!\mathbf{B}),
\label{sing_shrink_from_eig_shrink}
\end{equation}
where $h(x) \!:=\! g(\sqrt{x})/\sqrt{x}$ in $\mathcal{H}(\cdot)$ in \eqref{eigen_shrink_normal}.
Equation \eqref{sing_shrink_from_eig_shrink} is derived as follows.
First, \eqref{sing_shrink_exact} can be expanded as
\begin{equation}
\begin{aligned}
\mathcal{G}(\mathbf{B}) & = \mathbf{U}\mathbf{\Sigma}\,\mathrm{diag}\left(\frac{g(\sigma_1)}{\sigma_1}, \ldots , \frac{g(\sigma_n)}{\sigma_n}\right)\mathbf{V}^\top
\\
& = \mathbf{U}\mathbf{\Sigma}\mathbf{V}^\top\mathbf{V}\,\mathrm{diag}\left(\frac{g(\sigma_1)}{\sigma_1}, \ldots , \frac{g(\sigma_n)}{\sigma_n}\right)\mathbf{V}^\top
\\
& = \mathbf{B}\mathbf{V}\,\mathrm{diag}\left(\frac{g(\sigma_1)}{\sigma_1}, \ldots , \frac{g(\sigma_n)}{\sigma_n}\right)\mathbf{V}^\top.
\end{aligned}
\label{sing_shrink_exact_from_eig}
\end{equation}
When the eigenvalue matrix of $\mathbf{B}^{\!\top}\!\mathbf{B}$ is defined as $\mathbf{\Lambda}_{\mathrm{B}^{\!\top}\!\mathrm{B}}\!=\! \mathrm{diag}(\lambda^{\mathrm{B}^{\!\top}\!\mathrm{B}}_1,\ldots, \lambda^{\mathrm{B}^{\!\top}\!\mathrm{B}}_n)$, it is obviously represented using the singular values of $\mathbf{B}$ as $\lambda^{\mathrm{B}^{\!\top}\!\mathrm{B}}_i \!=\! \sigma^2_i$.
Consequently, \eqref{sing_shrink_exact_from_eig} is equally calculated using $\mathcal{H}(\cdot)$ in \eqref{eigen_shrink_normal} as
\begin{equation}
\begin{aligned}
\mathcal{G}(\mathbf{B}) & = \mathbf{B}\mathbf{V}\,\mathrm{diag}\left(\frac{g(\sigma_1)}{\sigma_1}, \ldots , \frac{g(\sigma_n)}{\sigma_n}\right)\mathbf{V}^\top
\\
& = \mathbf{B}\mathbf{V}\,\mathrm{diag}\bigl( h(\sigma^2_1), \ldots , h(\sigma^2_n) \bigr)\mathbf{V}^\top
\\
& = \mathbf{B} \, \mathcal{H}(\mathbf{B}^{\!\top}\!\mathbf{B}).
\end{aligned}
\label{sing_shrink_exact_from_eig2}
\end{equation}
Note that \cite{polynomial_filtering1} aims to calculate a vector represented as
\begin{equation}
\widehat{\mathbf{x}} = \mathbf{B}\mathcal{H}(\mathbf{B}^{\!\top}\!\mathbf{B})\mathbf{x},
\label{prop_another_form}
\end{equation}
where $\mathbf{x} \in \mathbb{R}^{n}$ and $\widehat{\mathbf{x}} \in \mathbb{R}^{m}$ are the input and output vectors, respectively.
The CPA is applied to $\mathcal{H}(\mathbf{B}^{\!\top}\!\mathbf{B})\mathbf{x}$ to quickly derive $\widehat{\mathbf{x}}$ in \cite{polynomial_filtering1}.
In contrast, our method is focused on deriving the matrix in \eqref{sing_shrink_exact_from_eig2} itself.
When the matrix, whose singular values are shrunk using CPA, is represented as $\mathbf{B}\mathcal{H}(\mathbf{B}^{\!\top}\!\mathbf{B})$, the explicit SVD of $\mathbf{B}$ can be avoided.
However, deriving the matrix, not the vectors, usually requires enormous computation time because of multiplication of dense matrices.
To accelerate the calculation, the sparseness of a matrix is exploited with our method, as indicated below.

Assume that $\mathbf{B}$ is a matrix composed of an inherently sparse signal.
Let $\mathbf{T}\!\in\! \mathbb{R}^{n\times n}$ be an arbitrary orthogonal matrix that efficiently sparsifies $\mathbf{B}$, e.g., $\mathbf{T}$ is considered as the discrete Fourier transform \cite{dft}, discrete cosine transform (DCT) \cite{dct}, and discrete wavelet transform (DWT) \cite{DWT}.
Note that the matrix $\mathbf{T}$ is the forward transform, i.e., when let $\mathbf{y}\in \mathbb{R}^n$ be a column vector, the forward transform is represented as $\mathbf{Ty}$.
From the above, \eqref{sing_shrink_from_eig_shrink} is further rewritten with $\mathbf{T}$ as
\begin{equation}
\mathcal{G}(\mathbf{B}) =\mathbf{B}\mathbf{T}^\top\mathcal{H}(\mathbf{T}\mathbf{B}^{\!\top}\!\mathbf{B}\mathbf{T}^\top)\mathbf{T} = \mathbf{B}\mathbf{T}^\top\mathcal{H}( \mathbf{\Phi} )\mathbf{T},
\label{sing_shrink_exact_sparse}
\end{equation}
where $\mathbf{\Phi} := \mathbf{T}\mathbf{B}^{\!\top}\!\mathbf{B}\mathbf{T}^\top$ for simplicity.
With the CPA-based eigenvalue shrinkage function $\widehat{\mathcal{H}}(\cdot)$ in \eqref{chebychev_series_matrix} and \eqref{chebychev_series_matrix_shift2}, $\mathcal{H}(\mathbf{\Phi})$ in \eqref{sing_shrink_exact_sparse} is efficiently approximated as
\begin{equation}
\widehat{\mathcal{H}}(\mathbf{\Phi})= \mathbf{T}\mathbf{V}\,\mathrm{diag}\bigl( \widehat{h}(\sigma^2_1), \ldots, \widehat{h}(\sigma^2_n) \bigr) \mathbf{V}^\top\mathbf{T}^\top.
\label{ATA_eig_fil}
\end{equation}
The form of \eqref{ATA_eig_fil} enables us to use the sparsity of a signal in its frequency domain.

For further enhancing the sparsity of $\mathbf{\Phi}$, its components are thresholded as
\begin{equation}
\underline{\Phi_{ij}}=
\begin{cases}
\Phi_{ij} &\text{if}~|\Phi_{ij}|\ge\varepsilon, \\
0 & \text{otherwise},
\end{cases}
\label{threshold_components}
\end{equation}
where $\Phi_{ij}$ and $\underline{\Phi_{ij}}$ are the $i$-th row and $j$-th column of $\mathbf{\Phi}$ and its truncated coefficient, and $\varepsilon\!\in\! \mathbb{R}$ is an arbitrary small value.
We show that this truncation has little effect on the performance of our method and provides recommended settings of $\varepsilon$ in Section~\ref{sec:applications}.
As a result, the CPA-based eigenvalue shrinkage of $\widehat{\mathcal{H}}(\mathbf{\Phi})$ is approximately given by
\begin{equation}
\widehat{\mathcal{H}}(\mathbf{\Phi}) \approx \widehat{\mathcal{H}}\left(\underline{\mathbf{\Phi}}\right).
\label{approx_eigen_shrink_trans_coeff}
\end{equation}
\begin{algorithm}
\caption{CPA-based singular value shrinkage}         
\label{alg1}
\begin{algorithmic}[1]
\REQUIRE $\mathbf{B}$
\ENSURE $\widehat{\mathcal{G}}(\mathbf{B})$
\STATE $\mathbf{\Phi} \leftarrow \mathbf{T}\mathbf{B}^{\!\top}\!\mathbf{B}\mathbf{T}^\top$.
\STATE Derive $\underline{\mathbf{\Phi}}$ from $\mathbf{\Phi}$ by using an arbitrary $\varepsilon$ in \eqref{threshold_components}.
\STATE Compute the maximum eigenvalue $\Lambda_\mathrm{max}$ of $\underline{\mathbf{\Phi}}$.
\STATE $\widehat{c}_k \leftarrow \frac{2}{\alpha}\sum^\alpha_{l=1} \cos( k\theta(l)  ) \, h\!\hspace{0.4mm}\Bigl( \frac{\Lambda_\mathrm{max}}{2}( \cos \theta(l)  +1 ) \Bigr)$,\\where $h(x) \!:=\! g(\sqrt{x})/\sqrt{x}$.
\STATE $\widehat{\mathbf{B}} \leftarrow \frac{2}{\Lambda_\mathrm{max}}\underline{\mathbf{\Phi}} - \mathbf{Id}$.
\STATE $\mathit{\Psi}_0(\widehat{\mathbf{B}})\leftarrow \mathbf{Id}$, $\mathit{\Psi}_1(\widehat{\mathbf{B}})\leftarrow \widehat{\mathbf{B}}$.
\STATE  $\hspace{2pt}\widehat{\mathcal{H}}\left(\underline{\mathbf{\Phi}}\right)\leftarrow \frac{1}{2}\widehat{c}_0\mathit{\Psi}_0(\widehat{\mathbf{B}})+\widehat{c}_1\mathit{\Psi}_1(\widehat{\mathbf{B}})$.
\FOR{$k=2$ to $\alpha-1$}
\STATE $\mathit{\Psi}_k (\widehat{\mathbf{B}}) \leftarrow2\widehat{\mathbf{B}} \, \mathit{\Psi}_{k-1} (\widehat{\mathbf{B}}) - \mathit{\Psi}_{k-2} (\widehat{\mathbf{B}})$.
\STATE $\hspace{2.5pt}\widehat{\mathcal{H}}\left(\underline{\mathbf{\Phi}}\right)\leftarrow \widehat{\mathcal{H}}\left(\underline{\mathbf{\Phi}}\right)+\widehat{c}_k\mathit{\Psi}_k(\widehat{\mathbf{B}})$.
\ENDFOR
\STATE $\widehat{\mathcal{G}}(\mathbf{B})\leftarrow \mathbf{B}\mathbf{T}^\top\widehat{\mathcal{H}}\left(\underline{\mathbf{\Phi}}\right)\mathbf{T}$.
\end{algorithmic}
\end{algorithm}
In summary, the singular value shrinkage of $\mathbf{B}$ is approximately represented as
\begin{align}
\mathcal{G}(\mathbf{B}) & = \mathbf{U}
\begin{bmatrix}
 	g(\sigma_1) & &\mathbf{O}\\
	 &\ddots&\\
	 &&g(\sigma_n)\\
	\mathbf{O}&&\\
	\end{bmatrix}
\mathbf{V}^\top
\nonumber
\\
& = \mathbf{B} \, \mathcal{H}(\mathbf{B}^{\!\top}\!\mathbf{B})
\nonumber
\\
& = \mathbf{B}\mathbf{T}^\top\mathcal{H}(\mathbf{\Phi})\mathbf{T}
\nonumber
\\
& \approx \mathbf{B}\mathbf{T}^\top\widehat{\mathcal{H}}(\mathbf{\Phi})\mathbf{T}\nonumber
\\
& \approx \mathbf{B}\mathbf{T}^\top\widehat{\mathcal{H}}\left(\underline{\mathbf{\Phi}}\right)\mathbf{T}= \widehat{\mathcal{G}}(\mathbf{B}), 
\label{SVD_fil}
\end{align}
where the function $\widehat{\mathcal{G}}(\cdot)$ is the CPA-based singular value shrinkage function.
It can be calculated with the recurrence relation as
\begin{equation}
\begin{split}
&\mathit{\Psi}_k (\widehat{\mathbf{B}}) =2\widehat{\mathbf{B}} \, \mathit{\Psi}_{k-1} (\widehat{\mathbf{B}}) - \mathit{\Psi}_{k-2} (\widehat{\mathbf{B}}),
\\
&\mathit{\Psi}_0 (\widehat{\mathbf{B}}) = \mathbf{Id},\quad \mathit{\Psi}_1 (\widehat{\mathbf{B}}) = \widehat{\mathbf{B}},
\end{split}
\label{recurrence_our_method}
\end{equation}
where
\begin{equation}
\widehat{\mathbf{B}} := \frac{2}{\Lambda_\mathrm{max}}\underline{\mathbf{\Phi}} - \mathbf{Id},
\end{equation}
in which $\Lambda_\mathrm{max}\!=\!\lambda^{\underline{\mathbf{\Phi}}}_\mathrm{max}$ is the maximum eigenvalue of $\underline{\mathbf{\Phi}}$.
The pseudocode of the CPA-based singular value shrinkage is indicated in Algorithm~\ref{alg1}.

%
%
\subsection{Computational Complexity of CPA-based Singular Value Shrinkage}

We now discuss the computational complexity of our method.
Assume that matrices $\widehat{\mathbf{B}}\in \mathbb{R}^{n\times n}$ and $\mathit{\Psi}_k(\widehat{\mathbf{B}})\in \mathbb{R}^{n\times n}$ have $M$ and $M_k$ nonzero elements, respectively.
The maximum number of multiplications of nonzero elements required to calculate $\widehat{\mathbf{B}}\mathit{\Psi}_k(\widehat{\mathbf{B}})$ is represented as $MM_k$ in the case of a sparse matrix.
The computational complexity of line 9 in Algorithm~\ref{alg1} can be represented as $\mathcal{O}(\sum^{\alpha-2}_{k=1}MM_k)$ due to the multiplication $\widehat{\mathbf{B}}\mathit{\Psi}_{k-1}(\widehat{\mathbf{B}})$.
At line 10 in Algorithm~\ref{alg1}, the computation takes $\mathcal{O}(\sum^{\alpha-1}_{k=2}M_k)$ from the multiplication $c_k \mathit{\Psi}_k(\widehat{\mathbf{B}})$.
That is, the total computational complexity is represented as $\mathcal{O}((M+1) \max_k \{M_k\})$, where $\max_k \{ M_k\}$ represents the maximum value among $M_k$.
From the above, when $\max_k \{ M_k\}$ becomes small, the computational cost is also reduced.
\begin{figure}[htb]
\small
\tabcolsep = 0.4mm
\centering
\begin{tabular}{cc}
\includegraphics[trim=0 0 0 0, width=0.48\linewidth, clip]{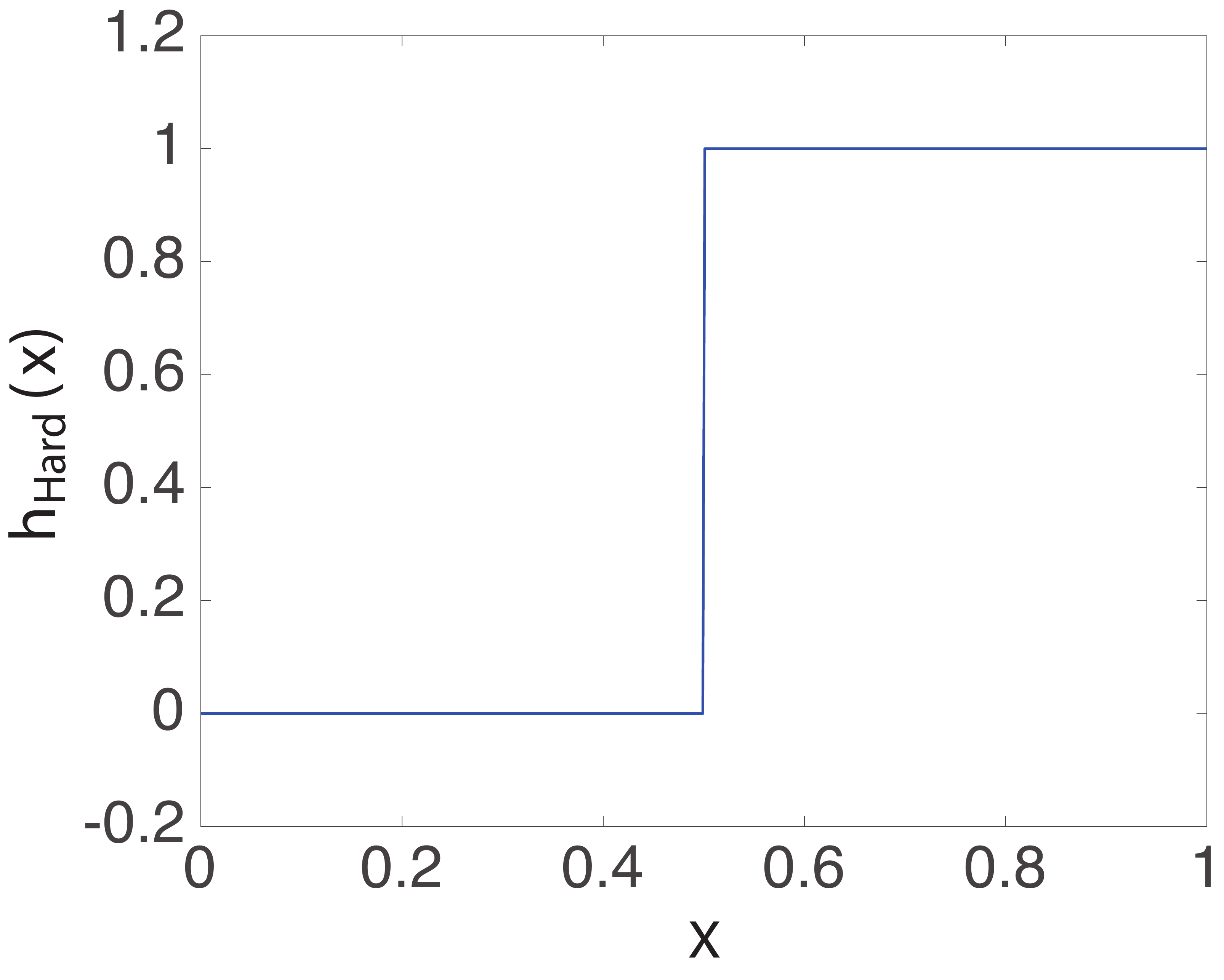}&
\includegraphics[trim=0 0 0 0, width=0.48\linewidth, clip]{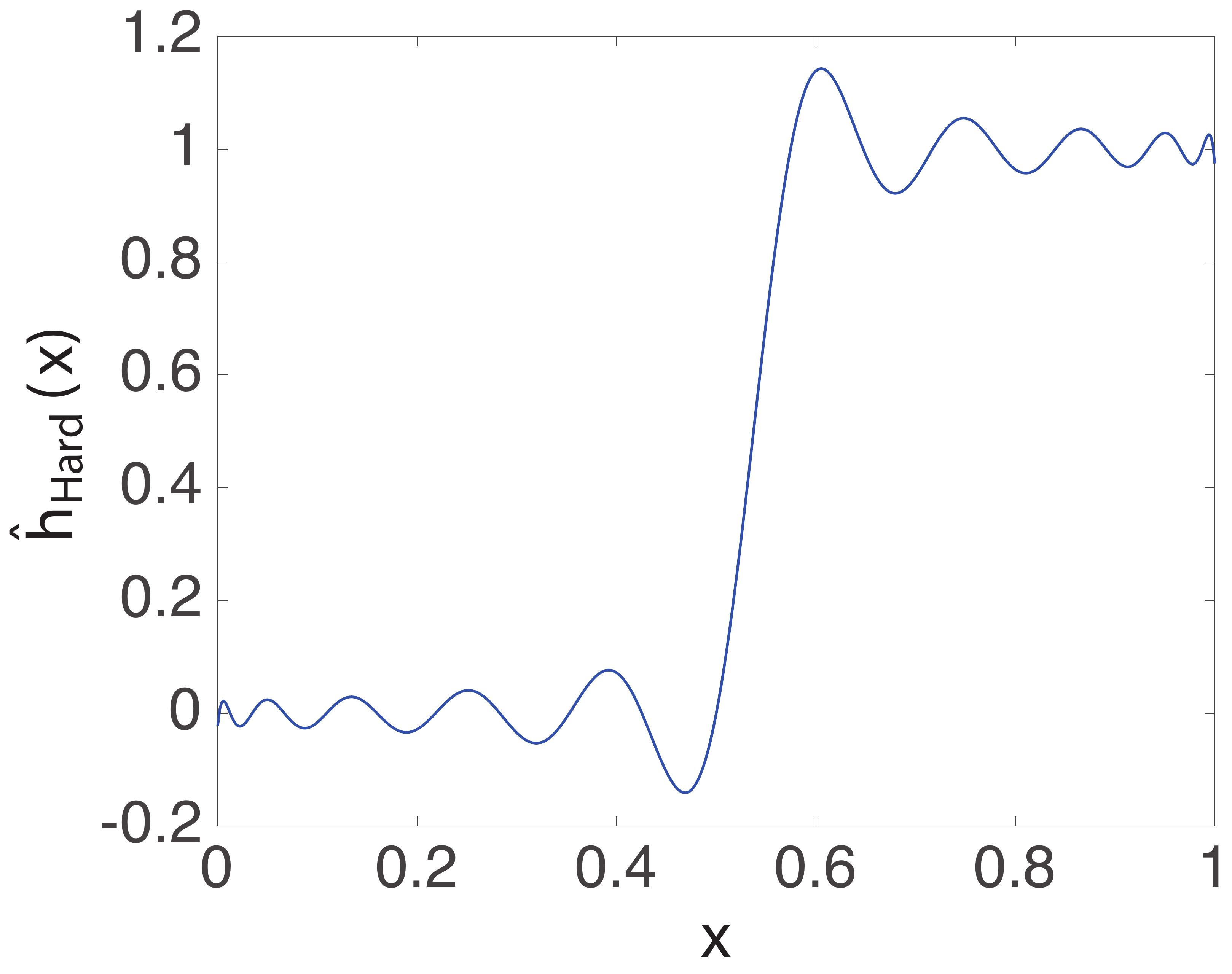}
\\
(a) & (b)
\end{tabular}
\caption{Example shrinkage responses ($s_\mathrm{c} \!=\! 0.5$).
(a) Hard-shrinkage response.
(b) Approximate shrinkage response ($\alpha = 30$).}
\label{fig:ideal_filter_res}
\end{figure}For low computational complexity, the matrix should be constructed so as not to increase the number of its nonzero elements as much as possible in the multiplication of matrices.

%
%
\section{Shrinkage Functions and Approximation Order}
\label{sec:chebyshev_filter_design}

In this section, we discuss suitable approximation orders for shrinkage functions approximated by CPA, which has small truncation errors.
Additionally, we argue that CPA is a reasonable choice for our method among a variety of polynomial approximation methods.

As an introduction, we consider the shrinkage function shown in Fig.~\ref{fig:ideal_filter_res}(a).
Let $h_\mathrm{hard}(x)$ be the hard shrinkage response defined as
\begin{equation}
h_\mathrm{hard}(x) :=
\begin{cases}
1 & \text{if}~x > \tau_\mathrm{hard},\\
0 & \text{otherwise},
\end{cases}
\label{ideal_filter_response}
\end{equation}
where $\tau_\mathrm{hard}$ is an arbitrary real value and $x\!\in\! [0,1]$.
Chebyshev polynomial approximation gives an approximate response of $h_\mathrm{hard}(x)$ in \eqref{chebychev_series}.
As in Fig.~\ref{fig:ideal_filter_res}(b), the approximated response $\widehat{h}_\mathrm{hard}(x)$ has \textit{ripples}, which is widely known in digital filter design \cite{minimax2,minimax3,kaiser1,kaiser2,kaiser3}.
Therefore, studying the design of appropriate shrinkage responses and approximation orders is an important topic, even for our method.

%
%
\subsection{Approximation Order}

Possible shrinkage responses handled with our method can be expressed as the following generic form:
\begin{figure}[htb]
\small
\tabcolsep = 0.4mm
\centering
\begin{tabular}{cc}
\includegraphics[trim=0 0 165 0, width=0.48\linewidth, clip]{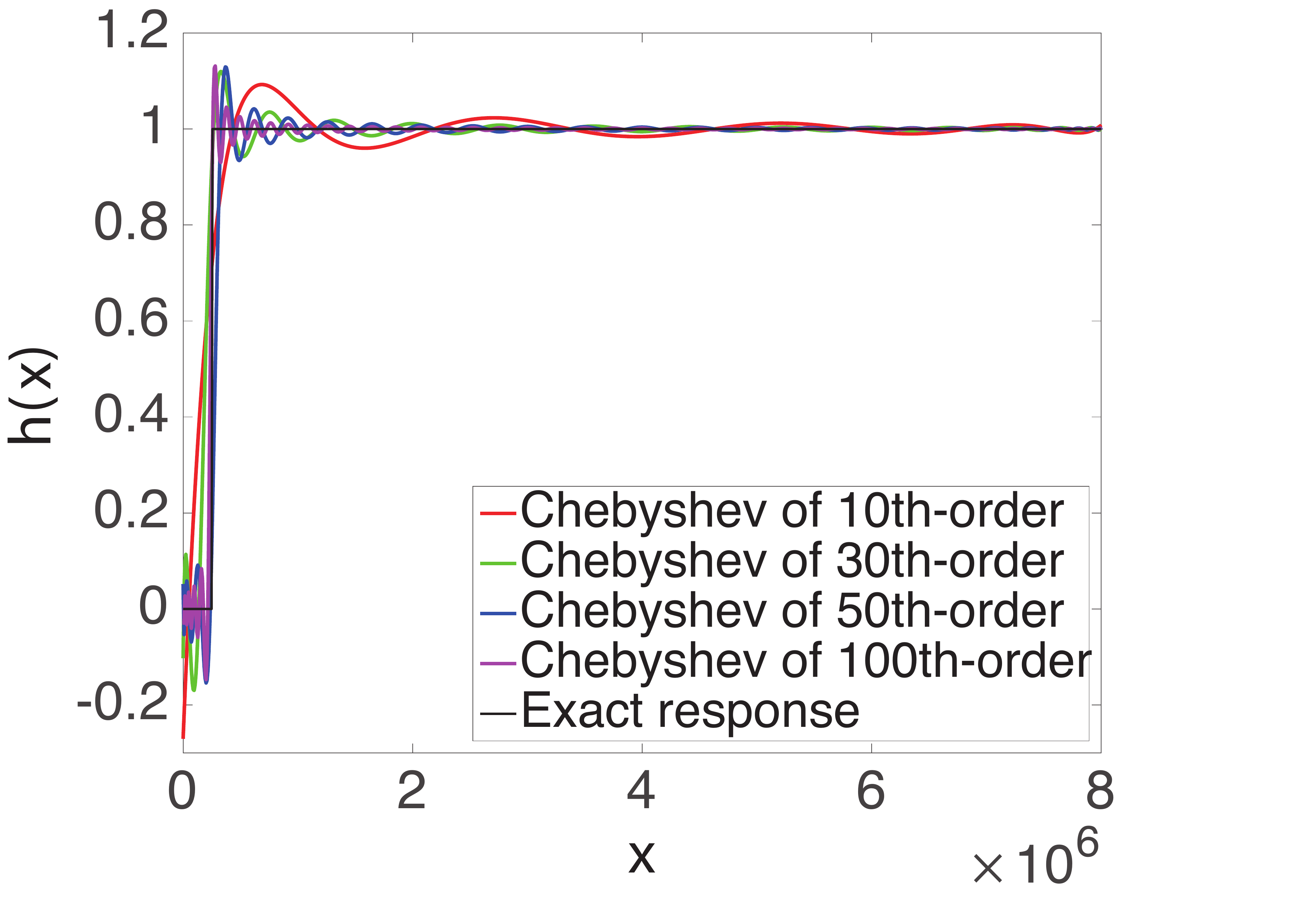}&
\includegraphics[trim=0 0 165 0, width=0.48\linewidth, clip]{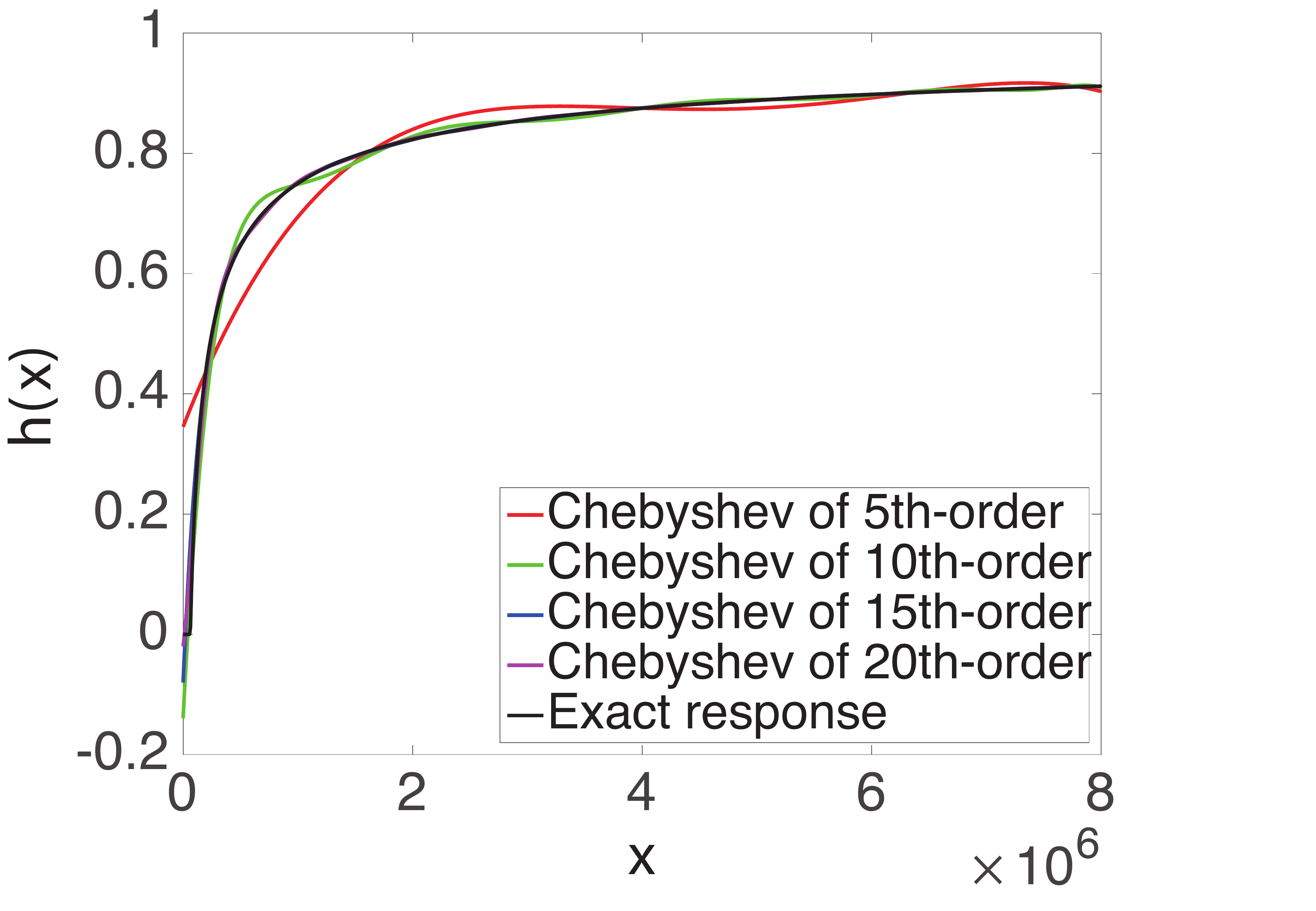}
\\
(a) & (b)
\\
\includegraphics[trim=0 0 165 0, width=0.48\linewidth, clip]{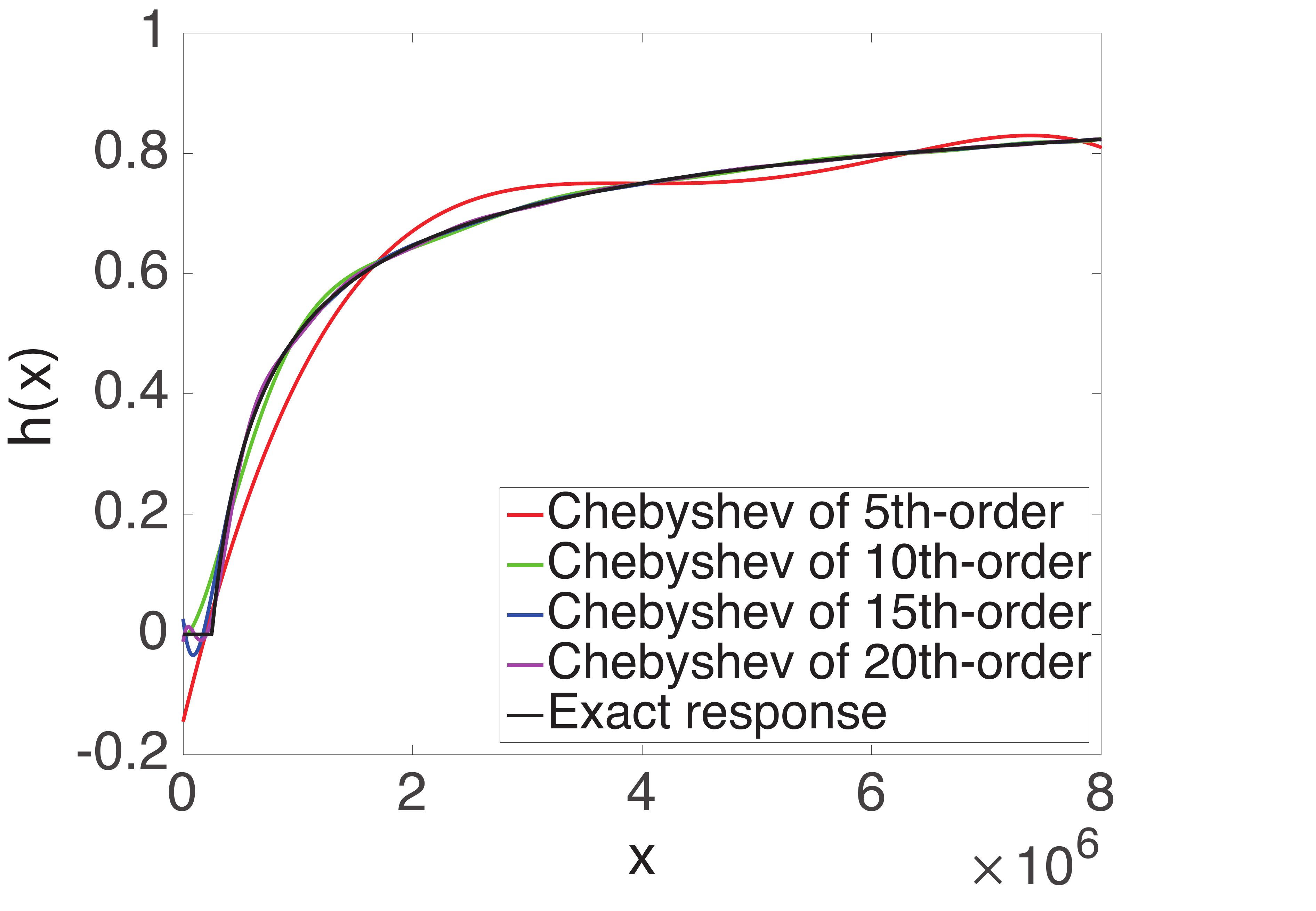}&
\includegraphics[trim=0 0 165 0, width=0.48\linewidth, clip]{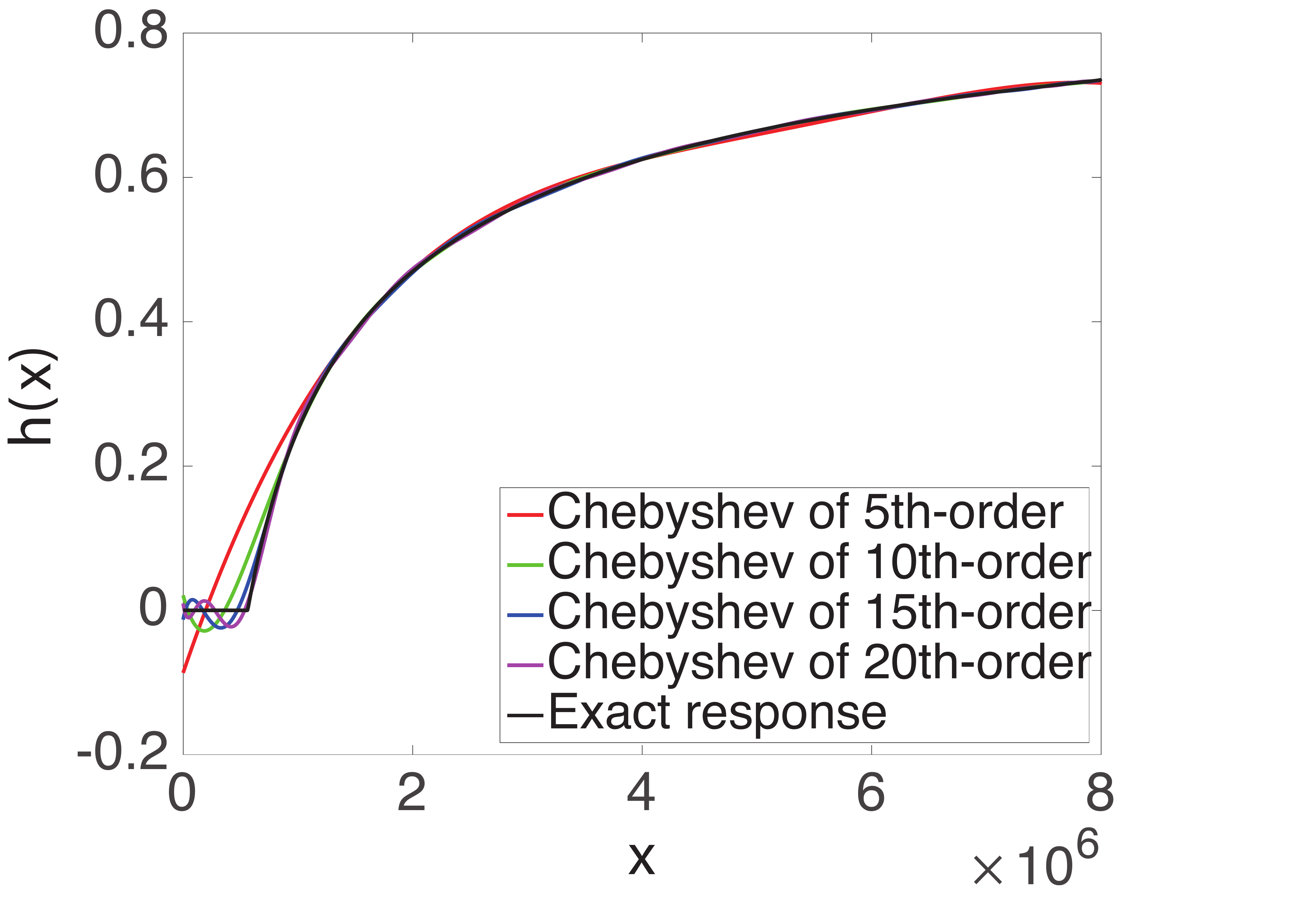}
\\
(c) & (d)
\end{tabular}
\caption{Comparison of responses between certain cases of \eqref{filter_response_ours} and their approximations by CPA.
(a) Hard-shrinkage.
(b) Weighted soft-shrinkage with $w(x)\!=\!0.5$.
(c) Soft-shrinkage.
(d) Weighted soft-shrinkage with $w(x)\!=\!1.5$.}
\label{fig:test_chebyshev_response}
\end{figure}
\begin{equation}
h\left(x;\frac{w(x)}{\rho},\tau\right) :=
\begin{cases}
\displaystyle\frac{\sqrt{x}-\frac{w(x)}{\rho}}{\sqrt{x}} & \text{if}~\sqrt{x} > \tau, \\
0 & \text{otherwise},
\end{cases}
\label{filter_response_ours}
\end{equation}
where $w(x)$ is a weight function, and $\tau$ and $\rho$ are arbitrary thresholding values.
The choices of $w(x)$ and $\tau$ determine the characteristics of \eqref{filter_response_ours} as follows:
\begin{itemize}
\item
$\displaystyle h\left(x;0,\tau_\mathrm{hard}\right)$ 
\quad : Hard-shrinkage.
\vspace{1mm}
\item
$\displaystyle h\left(x;\frac{w(x)}{\rho},\frac{w(x)}{\rho}\right)$ 
\quad : Weighted soft-shrinkage.
\vspace{1mm}
\item
$\displaystyle h\left(x;\frac{1}{\rho},\frac{1}{\rho}\right)$ 
\quad : Soft-shrinkage\footnote{Soft-shrinkage is widely known as $g(x) \!:=\! \max (x\!-\!1/\rho,0)$ in which $\max (x_1,x_2)$ is an operator choosing the greater one out of $x_1$ and $x_2$.
However, we call $h(x; 1/\rho, 1/\rho)$ soft-shrinkage because it is finally transformed into $g(x) \!:=\! \max (x\!-\!1/\rho,0)$ in \eqref{SVD_fil}.}.
\end{itemize}
Note that we defined $g(\sqrt{x})=\sqrt{x}h(x)$ in \eqref{sing_shrink_from_eig_shrink}, where $h(x)$ is an arbitrary shrinkage function for the eigenvalues of $\mathbf{B}^{\!\top}\!\mathbf{B}$.
\begin{figure*}[htb]
\small
\tabcolsep = 1.5mm
\centering
\begin{tabular}{ccc|c}
\hline
\multicolumn{3}{c|}{Responses} & Differences\\
\hline \hline
\includegraphics[trim=0 0 165 0, width=0.23\linewidth, clip]{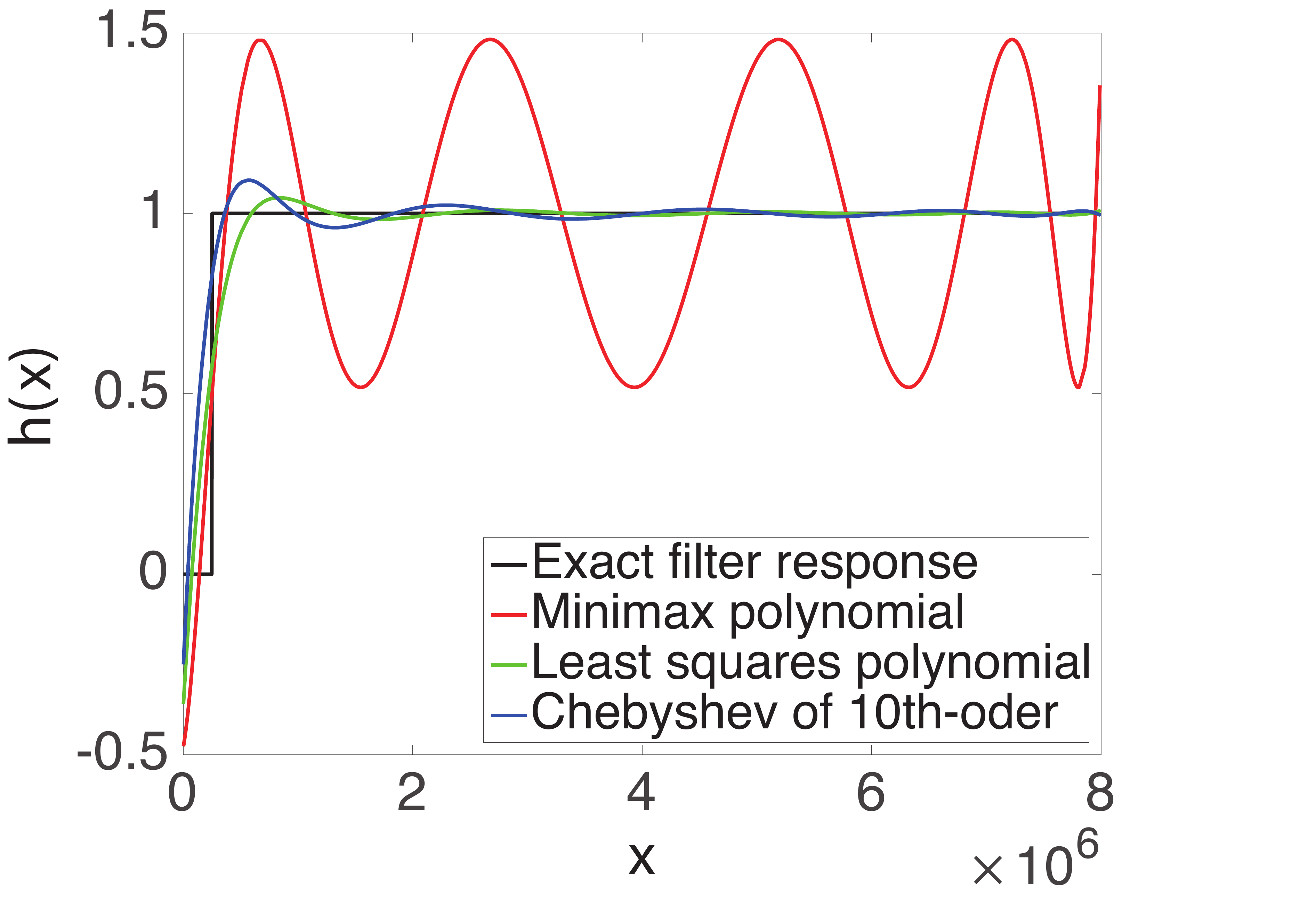}&
\includegraphics[trim=0 0 165 0, width=0.23\linewidth, clip]{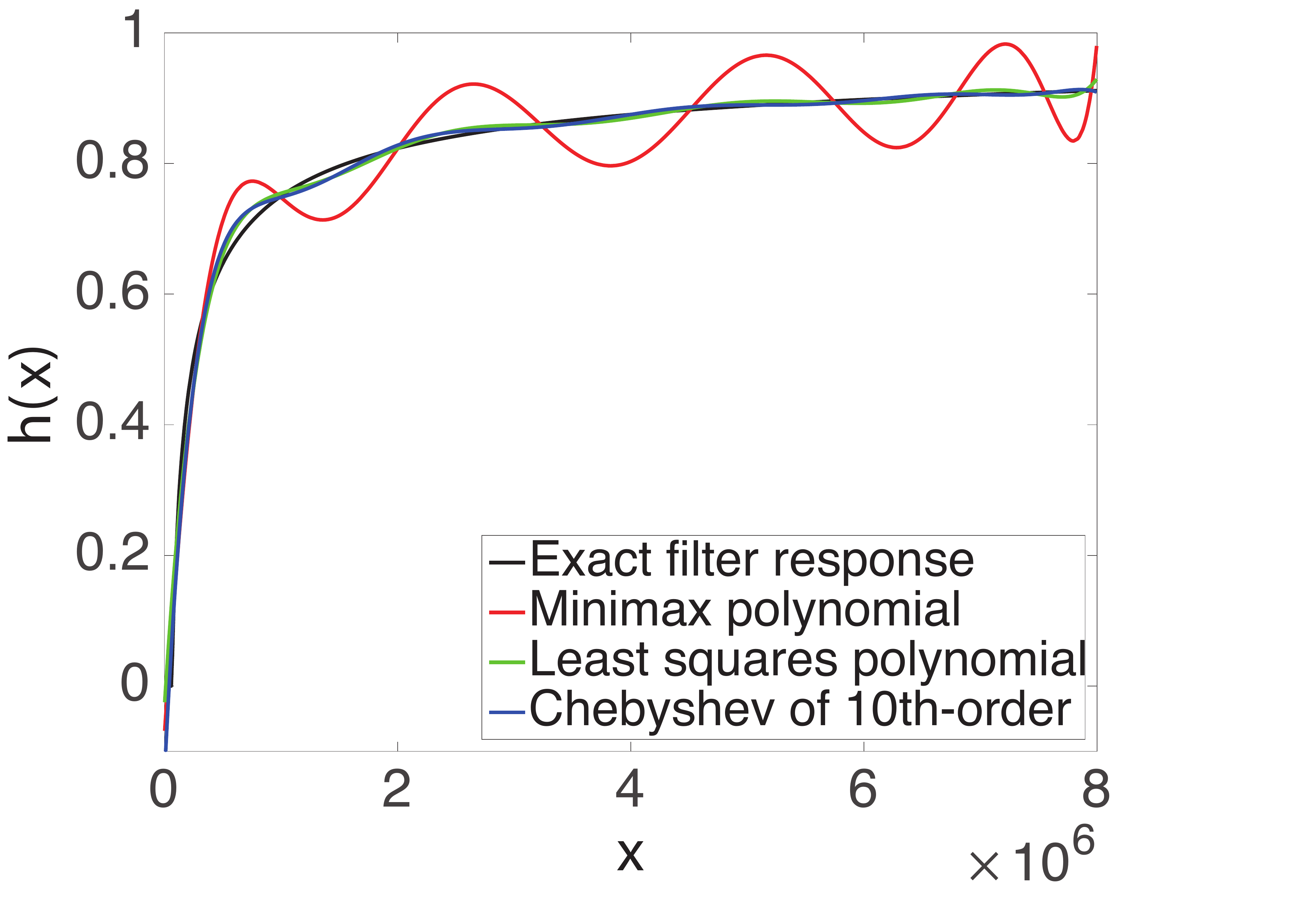}&
\includegraphics[trim=0 0 165 0, width=0.23\linewidth, clip]{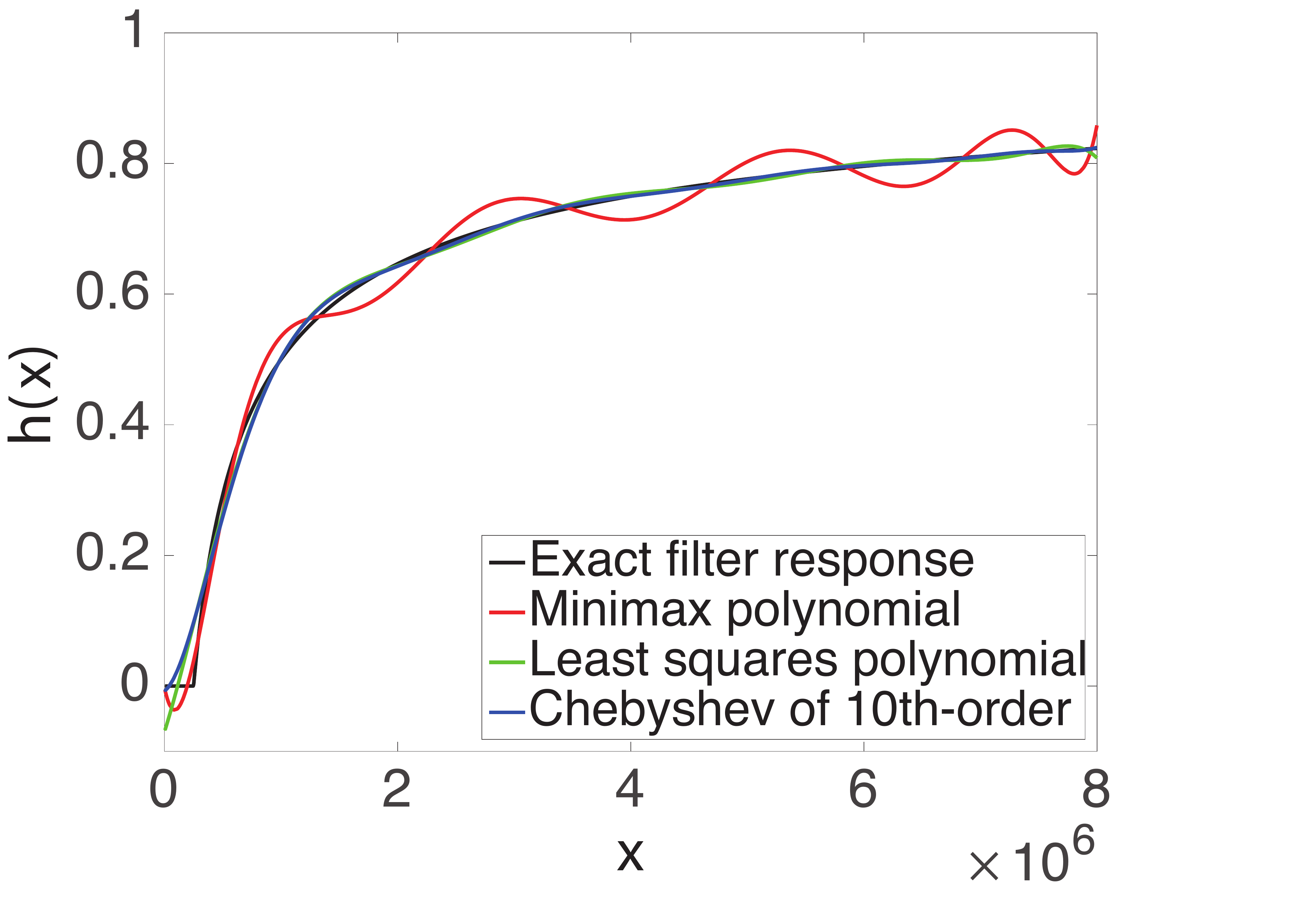}&
\includegraphics[trim=0 0 165 0, width=0.23\linewidth, clip]{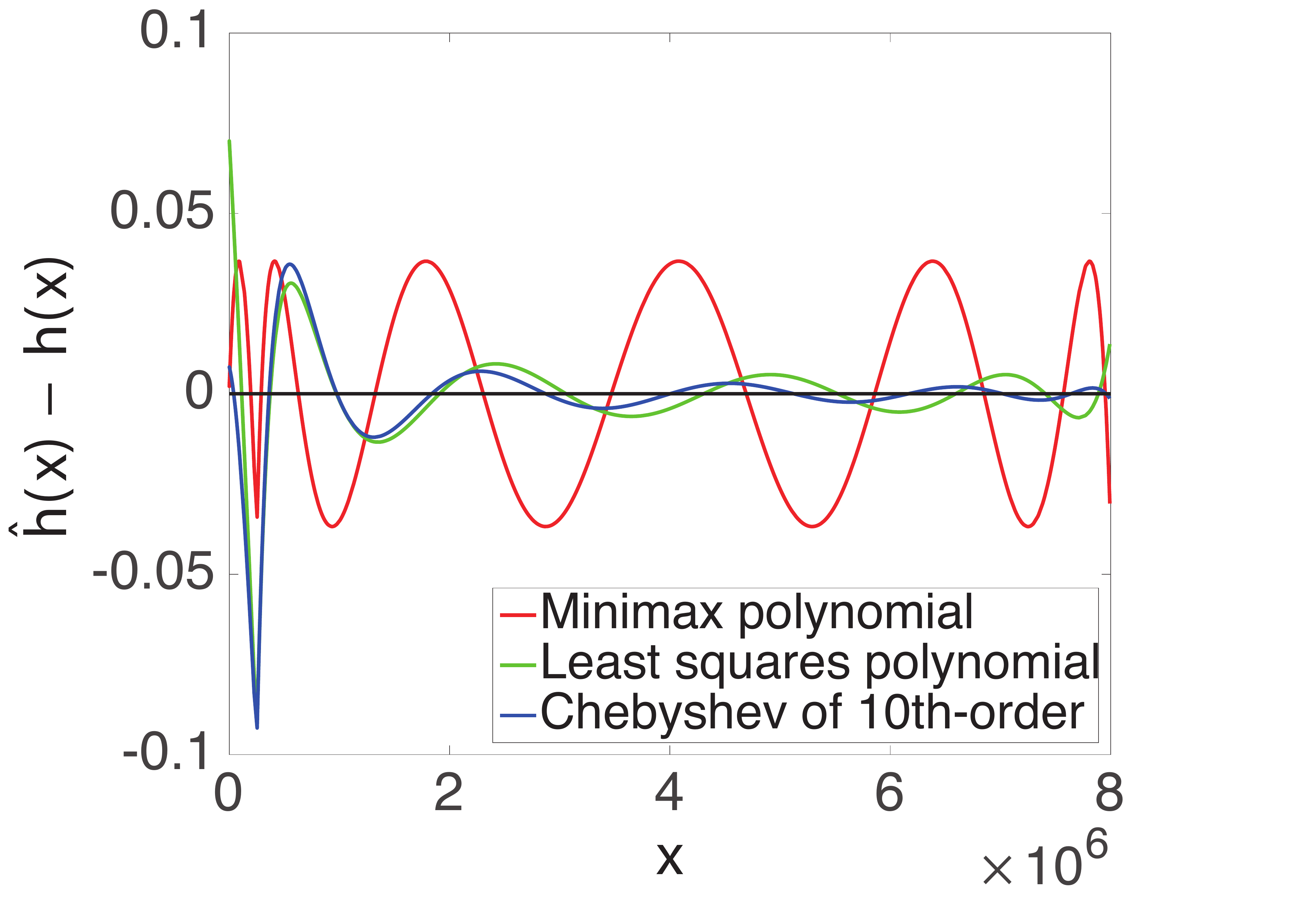}
\\
(a) Hard-shrinkage&
(b) Weighted soft-shrinkage&
(c) Soft-shrinkage&
Difference in (c)
\\
\includegraphics[trim=0 0 165 0, width=0.23\linewidth, clip]{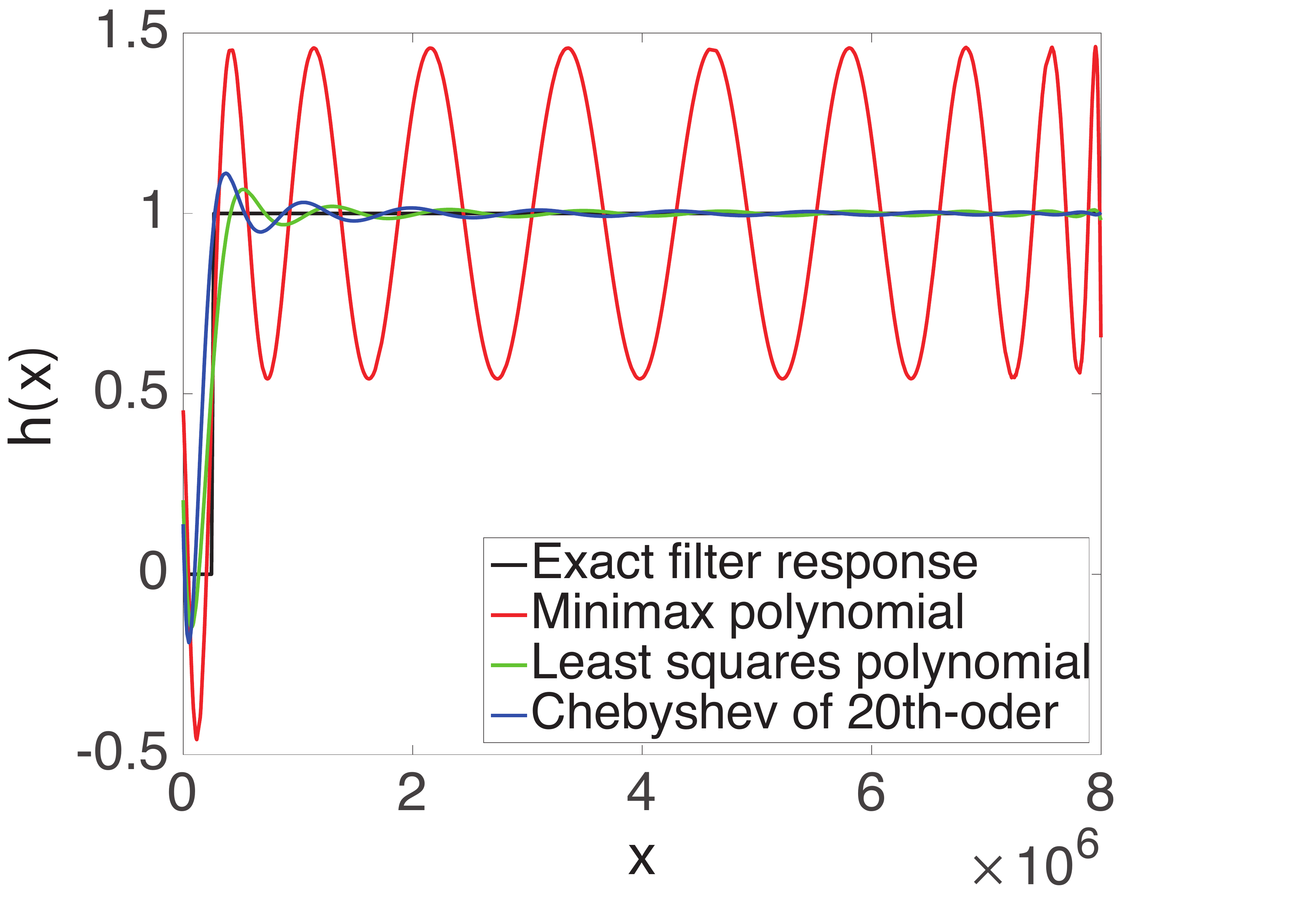}&
\includegraphics[trim=0 0 165 0, width=0.23\linewidth, clip]{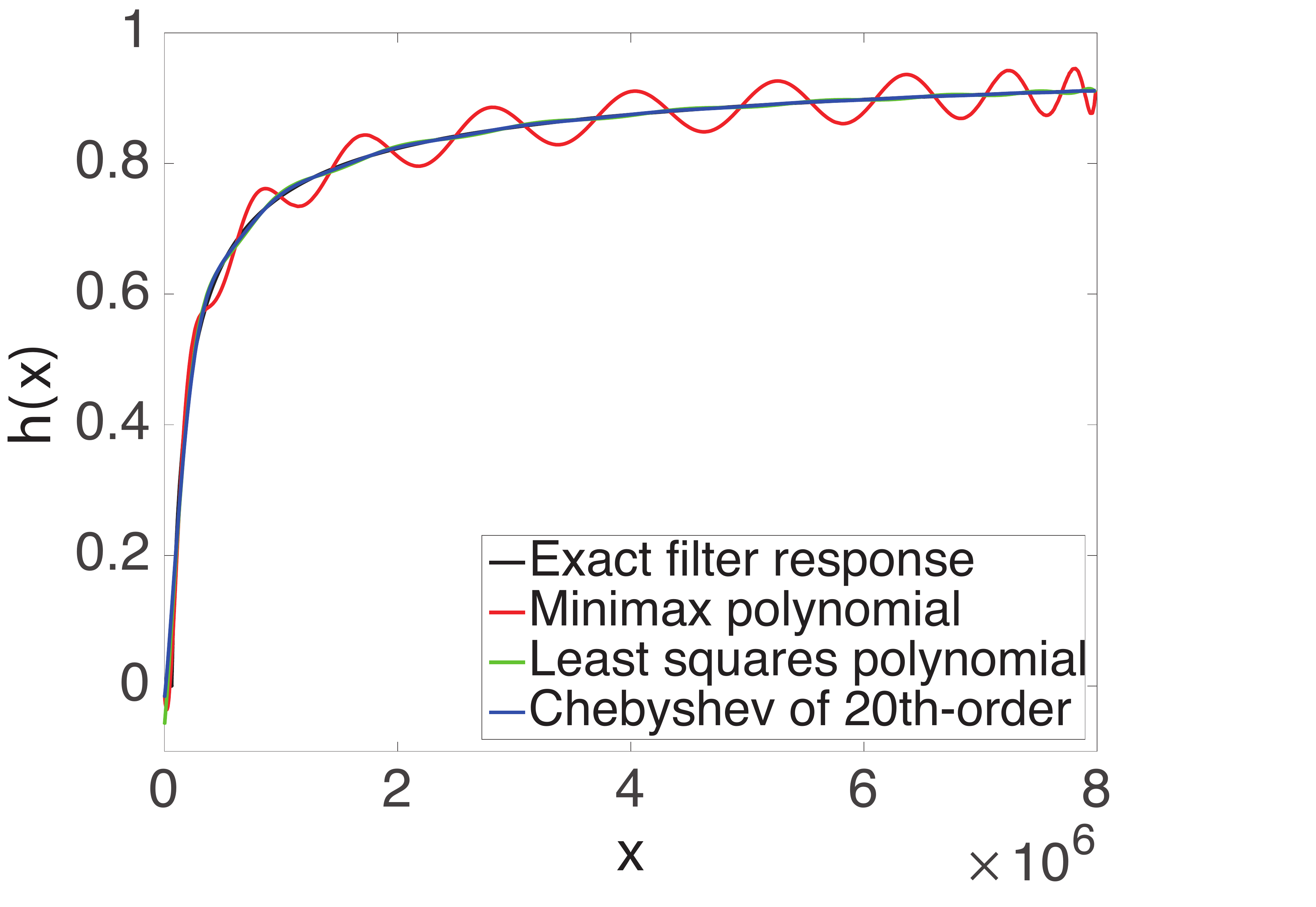}&
\includegraphics[trim=0 0 165 0, width=0.23\linewidth, clip]{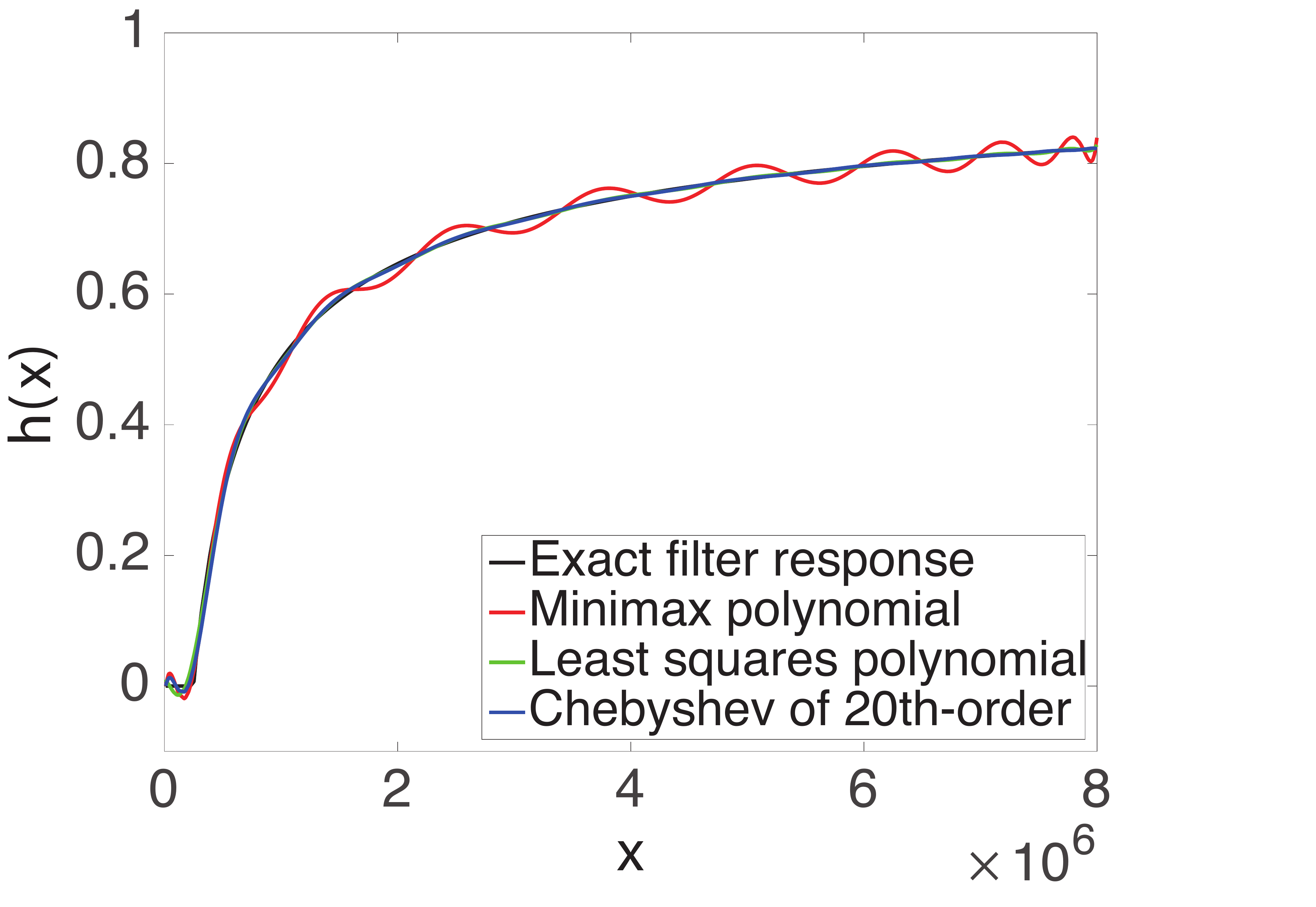}&
\includegraphics[trim=0 0 165 0, width=0.23\linewidth, clip]{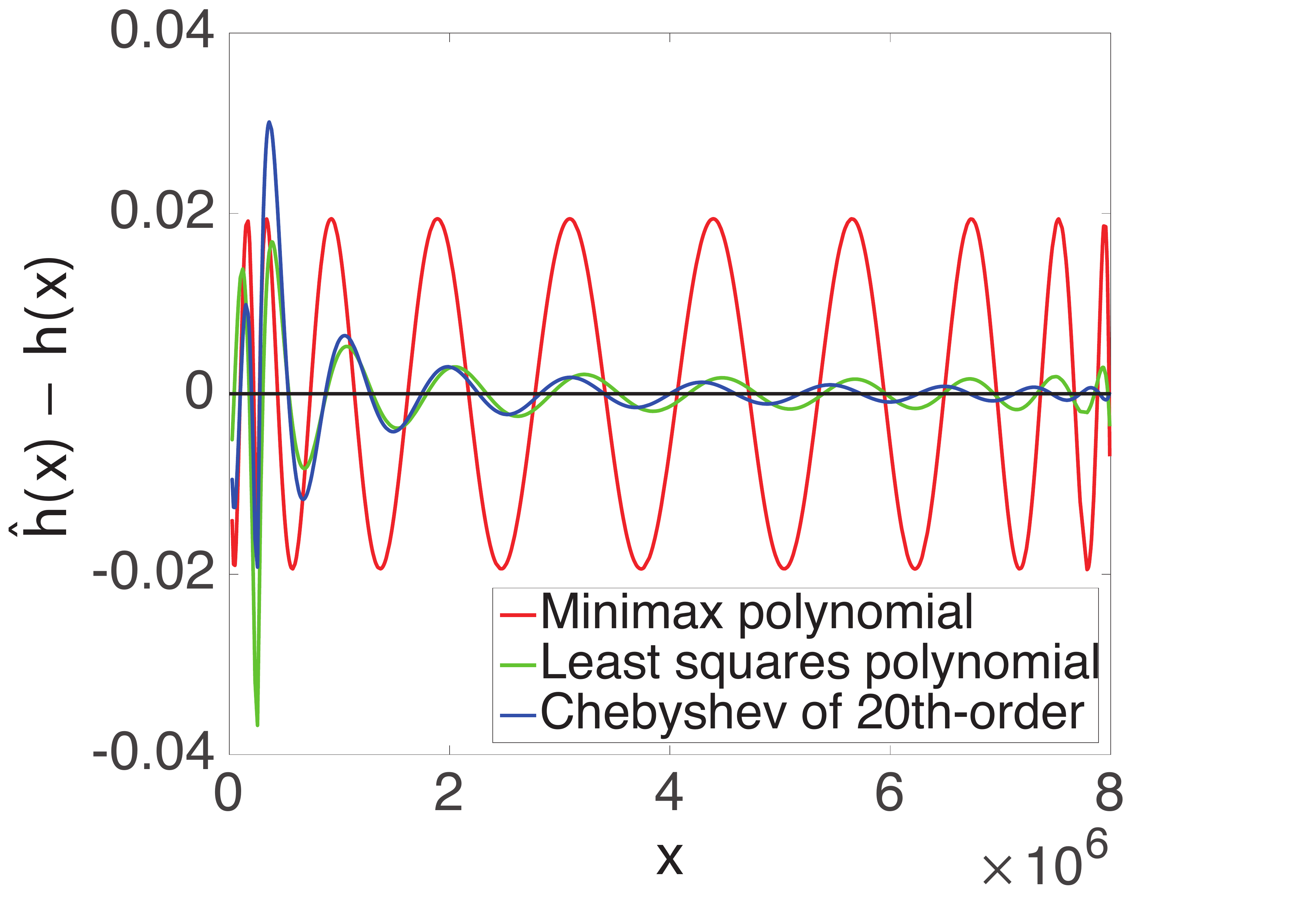}
\\
(d) Hard-shrinkage&
(e) Weighted soft-shrinkage&
(f) Soft-shrinkage&
Difference in (f)
\\ \hline
\end{tabular}
\caption{CPA compared with minimax polynomial and least squares polynomial.
(a)-(c) are shown in case of $10$th-order approximation and (d)-(f) are indicated in case of $20$th-order approximation.
In (b) and (e), $w(x) \!=\! 0.5$ was used.}
\label{fig:test_chebyshev_minimax}
\end{figure*}As a result, $g(x)$ becomes the hard-shrinkage, weighted soft-shrinkage, or soft-shrinkage functions when $h(x)$ is set as above.

These choices among the shrinkages are shown in Fig.~\ref{fig:test_chebyshev_response}.
It is clear that hard-shrinkage has a sharp transition band (see Fig.~\ref{fig:test_chebyshev_response}(a)); therefore, CPA, which is computed as a linear combination of cosine functions, may not approximate it well.
In contrast, one can expect that a response that has a smooth transition band is suitable for CPA.
To verify this numerically, the approximated responses were compared among hard-shrinkage, weighted soft-shrinkage, and soft-shrinkage.
In this experiment, $\rho \!=\! 0.002$ was used.
For weighted soft-shrinkage, $w(x) \!=\! 0.5$ and $1.5$ were used.
Additionally, the thresholding value for hard-shrinkage was set to $\tau_\mathrm{hard} \!=\! 500$.

Figure~\ref{fig:test_chebyshev_response} also shows the approximated shrinkage responses of \eqref{filter_response_ours} for various shrinkage conditions.
Hard-shrinkage yields larger errors than soft ones.
Empirically, hard-shrinkage requires more than the 50th-oder approximation.
In contrast, soft-shrinkages only require 10--20th-order approximations.
To be more specific, $\alpha \!=\! 20$ is recommended for a small weight shown in Fig.~\ref{fig:test_chebyshev_response}(b), whereas $\alpha \!=\! 10$ for the soft-shrinkage response shown in Fig.~\ref{fig:test_chebyshev_response}(d).

%
%
\subsection{Suitability of CPA}

There are many polynomial approximations.
Even among them, minimax polynomial approximation \cite{minimax1,minimax2,minimax3,minimax4,minimax5} and least squares approximation \cite{LMSE} are well known as the best approximation in the sense of the minimization of the infinity norm and the least squares error w.r.t the difference between an exact and approximated responses, respectively.
To derive polynomial coefficients, their optimization requires a minimization of $\ell_p$ norm represented as
\begin{equation}
\min_{h'(x)\in \mathbb{R}} \| h(x) - h'(x) \|_p,
\label{l_p_norm_for_approximation}
\end{equation}
where $h'(x)$ is an approximated shrinkage response with the above two polynomial approximations.
Clearly, \eqref{l_p_norm_for_approximation} requires the exact response $h(x)$ for $x\!\in\! \mathbb{R}$.
When $h(x)$ is precisely represented using many sampling points, $h'(x)$ exhibits good performance.
However, computational complexity becomes high when many sampling points are used, especially in the case of least squares approximation.
Let $\mathbf{c}_\alpha\!:=\![c_0,c_1,\ldots,c_{\alpha-1}]^\top$ be the column vector of coefficients for the polynomial approximation.
That is, the approximated shrinkage response can be calculated as $h'(x) \!=\! c_0 + c_1 x + c_2 x^2 +\ldots + c_{\alpha-1} x^{\alpha-1}$.
Additionally, let $\mathbf{x}\!:=\![x_1,\ldots,x_n]^\top$ and $\mathbf{h}\!:=\![h(x_1),\ldots,h(x_n)]^\top$ be the column vectors composed of real values, respectively.
The Vandermonde matrix $\mathbf{\Upsilon}\!\in\! \mathbb{R}^{n\times \alpha}$ is defined as
\begin{equation}
\mathbf{\Upsilon}:=
\begin{bmatrix}
1&x_1^{2}&\cdots&x_1^{\alpha-1}\\
1&x_2^{2}&\cdots&x_2^{\alpha-1}\\
\vdots&\vdots&\ddots&\vdots\\
1&x_n^{2}&\cdots&x_n^{\alpha-1}
\end{bmatrix}.
\end{equation}
From the above definitions, coefficients of the least squares approximation are calculated as $\mathbf{c}_\alpha \!=\! \mathbf{\Upsilon}^{+}\mathbf{h}$, where $\cdot^{+}$ is the pseudo inverse of a matrix.
The calculation requires high computational cost when $n$ and/or $\alpha$ are large.
In contrast, CPA only requires the inner product of $[h(\cos\theta_1),\ldots,h(\cos\theta_\alpha)]^\top$ and $[\cos k\theta_1,\ldots,\cos k\theta_\alpha]^\top$ to derive coefficients of polynomials from \eqref{chebychev_coeff_discrete}, where $\theta_i \!\in\! [0,\pi]$.
Additionally, CPA performs better approximation than other optimization methods.
To verify the exellent approximation, CPA was compared with minimax approximation and least squares approximation, as shown in Fig.~\ref{fig:test_chebyshev_minimax}.
As can be seen, CPA and least squares approximations have a similar oscillation pattern.
Chebyshev polynomial approximation sufficiently attenuates ripples, compared with the other methods in the stopband, as shown in the differences of Fig.~\ref{fig:test_chebyshev_minimax}.

%
%
\section{Applications}
\label{sec:applications}

We compared our CPA-based singular value shrinkage method with the exact and approximate singular value shrinkage methods.
Specifically, we applied our method to two applications using nuclear norm relaxation, i.e., inpainting of texture images and background subtraction of videos.
Additionally, we compared our CPA-based method with the existing methods, i.e., the exact partial singular value decomposition (PSVD) based method and fast singular value shrinkage methods \cite{SVD_echon2,SVD_echon3,Fast_singular_thresh}, in Section~\ref{sec:comparison_existing_meth}.
The computation time and approximation precision were indicated for the comparisons.

%
%
\subsection{Experimental Conditions}

The applications were implemented with MATLAB R2015b and run on a 3.2-GHz Intel Xeon E5-2667 processor with 512-GB RAM.
We compared our method with the SVD-based naive method (denoted as SVD-based method) in \eqref{sing_shrink_exact} and EVD-based methods in \eqref{sing_shrink_exact_from_eig} with respect to approximation precision and computation time.
Both SVD and EVD-based methods are exact singular value soft-shrinkage methods.
The EVD-based method\footnote{The EVD-based method could lead to loss of computational precision compared with the SVD-based one. Though the errors may affect the performance of applications, we did not encounter such a problem in the experiments described in this paper.} is usually faster than the SVD-based method and is widely used in many applications.
Therefore, the computation time of only the EVD-based method is indicated for the results of the exact methods.
With the SVD-based method, the SVD of an arbitrary matrix $\mathbf{X}\!\in\! \mathbb{R}^{m\times n}$ is first performed, then the obtained singular values are shrunk as $\max (\sigma_i (\mathbf{X})\!-\!1/\rho,0)$, where $\sigma_i (\mathbf{X})$ indicates the $i$-th largest singular value of $\mathbf{X}$.
The EVD-based method uses the relation between singular value shrinkage and eigenvalue shrinkage: the EVD of $\mathbf{X}^{\!\top}\!\mathbf{X}$ is first computed, then the obtained eigenvalues are shrunk 
as $\max \bigl(\sqrt{\lambda^{\mathrm{X}^{\!\top}\!\mathrm{X}}_i}-1/\rho,0\bigr) / \sqrt{\lambda^{\mathrm{X}^{\!\top}\!\mathrm{X}}_i}$, where $\lambda_i^{\mathbf{X}^{\!\top}\!\mathbf{X}}$ denotes the $i$-th largest eigenvalue of $\mathbf{X}^{\!\top}\!\mathbf{X}$, to derive singular value shrinkage.
Also, the shrinkage function $h(\cdot)$ in \eqref{sing_shrink_from_eig_shrink} is defined as the soft-shrinkage case given by $h\left(x;\frac{1}{\rho},\frac{1}{\rho}\right)$ from \eqref{filter_response_ours} for our method.
The DWT \cite{DWT} was used in \eqref{SVD_fil} to sparsify the signals.
We used Haar wavelet transform as the DWT.
In the DWT, one level transform was performed and all high frequency components were set to $0$.
The selection of a transform method naturally affects the computation time of our method.
Therefore, we indicate the effect of the selection in Section~\ref{sec:influence}.
To indicate the approximation precision, root mean squared error (RMSE) was used, which was computed using the results of our method and those of the SVD/EVD-based methods.
Furthermore, the computation times of all the methods are shown, and the average computation times of the CPA-based/exact singular value shrinkage in each iteration are also indicated.
In all applications, we used the 5th, 10th, 15th, and 20th-order approximations.
We also used the following optimization tools to solve the above applications.

%
%
\subsection{Optimization Tools}
\label{sec:optimization_tool}

\subsubsection{Proximity Operator}

Let $ \boldsymbol{\Gamma}_0(\mathbb{R}^N)$ be the set of all proper lower semicontinuous convex functions\footnote{A function $f : \mathbb{R}^N \!\rightarrow \mathbb{R} \cup \{ \infty \}$ is called \textit{proper lower semicontinuous convex} if $\mbox{dom}(f)\!:=\! \{ \mathbf{x}\!\in\!\mathbb{R}^N |~f(\mathbf{x}) \!<\! \infty \}\!\ne\! \emptyset$, $\mbox{lev}_{\le a}(f)\!:=\!\{ \mathbf{x} \!\in\! \mathbb{R}^N |~f(\mathbf{x}) \!\le\! a \}$ is closed in $\forall a\!\in\! \mathbb{R}$, and $f(\eta \mathbf{x}+(1-\eta)\mathbf{y}) \!\le\! \eta f(\mathbf{x})+(1-\eta)f(\mathbf{y})$ in $\forall\mathbf{x},\mathbf{y}\!\in\! \mathbb{R}^N$ and $\forall \eta\!\in\! (0,1)$, respectively.} over $\mathbb{R}^N$.
The \textit{proximity operator} \cite{prox1} of a function $f \!\in\! \boldsymbol{\Gamma}_0(\mathbb{R}^N)$ of index $\gamma \!>\! 0$ is defined as
\begin{equation}
\mbox{prox}_{\gamma f}:\mathbb{R}^N\!\rightarrow\mathbb{R}^N:\mathbf{x}
\mapsto
\argmin_{\mathbf{y} \in \mathbb{R}^N} \, f(\mathbf{y})+\frac{1}{2 \gamma} \| \mathbf{x} - \mathbf{y} \|^2.
\end{equation}
The proximity operator plays a central role in the optimization of applications, as discussed in this section.
When function $f$ is defined as the nuclear norm, i.e., $\mbox{prox}_{\gamma \| \cdot \|_*}$, the proximity operator can be calculated by singular value shrinkage with the thresholding parameter $\gamma$ \cite{matrix_tensor_completion2}.
Therefore, our CPA-based method is applied to the operator in the case of the nuclear norm.

\subsubsection{Alternating Direction Method of Multipliers}

The ADMM \cite{ADMM2} is an algorithm for solving a convex optimization problem represented as
\begin{equation}
\min_{\mathbf{x} \in \mathbb{R}^{n_1}, \mathbf{z} \in \mathbb{R}^{n_2}} f(\mathbf{x}) + g(\mathbf{z})
\quad \text{s.t.} \quad \mathbf{z} = \mathbf{K}\mathbf{x},
\label{ADMM_prob}
\end{equation}
where $f \!\in\! \boldsymbol{\Gamma}_0 (\mathbb{R}^{n_1})$, $g\!\in\! \boldsymbol{\Gamma}_0(\mathbb{R}^{n_2})$ and $\mathbf{K}\!\in\! \mathbb{R}^{n_2\times n_1}$.
For arbitrary $\mathbf{z}_{0}$, $\mathbf{p}_{0}\!\in\! \mathbb{R}^{n_2}$, and $\rho\!>\!0$, the ADMM algorithm is given by
\begin{figure}[htb]
\small
\begin{minipage}{1 \linewidth}
  \centering
  \centerline{\includegraphics[width=2.5cm]{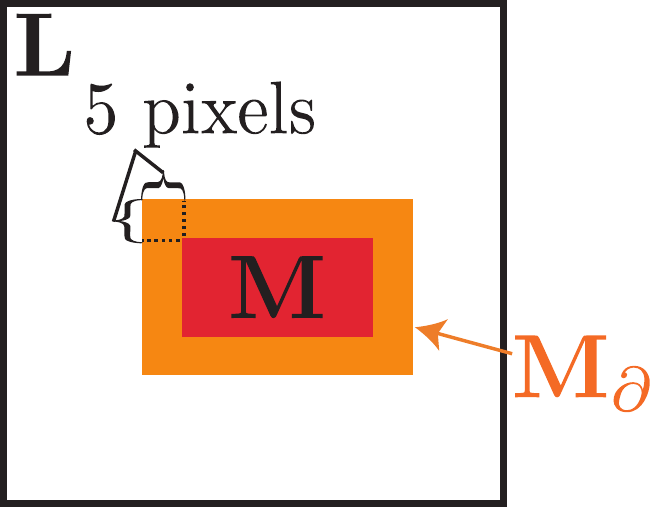}}
\end{minipage}
\caption{Missing region ($\mathbf{M}$) and its surrounding region ($\mathbf{M}_\partial$).}
\label{fig:t-t_near}
\end{figure}
\begin{equation}
\left\lfloor{\begin{split}
\mathbf{x}_{t+1} &:= \argmin_{\mathbf{x}} \, f(\mathbf{x}) + \frac{\rho}{2}\| \mathbf{z}_{t} - \mathbf{K}\mathbf{x} - \mathbf{u}_{t} \|^2_2\\
\mathbf{z}_{t+1}&:=\mbox{prox}_{1/\rho g}(\mathbf{K}\mathbf{x}_{t+1}+\mathbf{u}_{t})\\
\mathbf{u}_{t+1}&:=\mathbf{u}_{t}+\mathbf{K}\mathbf{x}_{t+1}-\mathbf{z}_{t+1}.\end{split}}\right.
\label{ADMM_algorithm}
\end{equation}
We recall a convergence analysis of the ADMM by Eskstein and Bertsekas \cite{ADMM2}.
\begin{fact}[Convergence of the ADMM \cite{ADMM2}]
Consider Prob.~\eqref{ADMM_prob}.
Assume that $\mathbf{K}^{\!\top}\!\mathbf{K}$ is invertible and that a saddle point of its unaugmented Lagrangian $L_0(\mathbf{x},\mathbf{z},\mathbf{u}')\!:=\!f(\mathbf{x})+g(\mathbf{z})-\langle\mathbf{u}',\mathbf{K}\mathbf{x}\!-\!\mathbf{z} \rangle$ exists, where $\mathbf{u}' \!:=\! \rho \mathbf{u}$.
Then the sequence $(\mathbf{x}_t)_{(t \ge 1)}$ generated using \eqref{ADMM_algorithm} converges to a solution of Prob.~\eqref{ADMM_prob}.
\end{fact}

We used the ADMM algorithm to practically solve the following applications.
In all applications, the stopping criterion\footnote{For example, in \eqref{admm_algorithm_inp}, which is indicated in Appendix~\ref{sec:admm_inpaint}, the criterion is evaluated using $\| \mathbf{l}_{t+1} - \mathbf{l}_{t} \|_2 / \| \mathbf{l}_{t+1}  \|_2$.} in the ADMM algorithm was set to $1.0\!\times\!10^{-4}$.

%
%
\subsection{Texture Image Inpainting \cite{repaire_texture1,repaire_texture2}}
\label{subsec:texture_image_inpainting}

\begin{figure*}[htb]
\small
\tabcolsep = 0.3mm
\centering
\begin{tabular}{ccccc}
\includegraphics[width=3.2cm,bb=0 0 2560 1920]{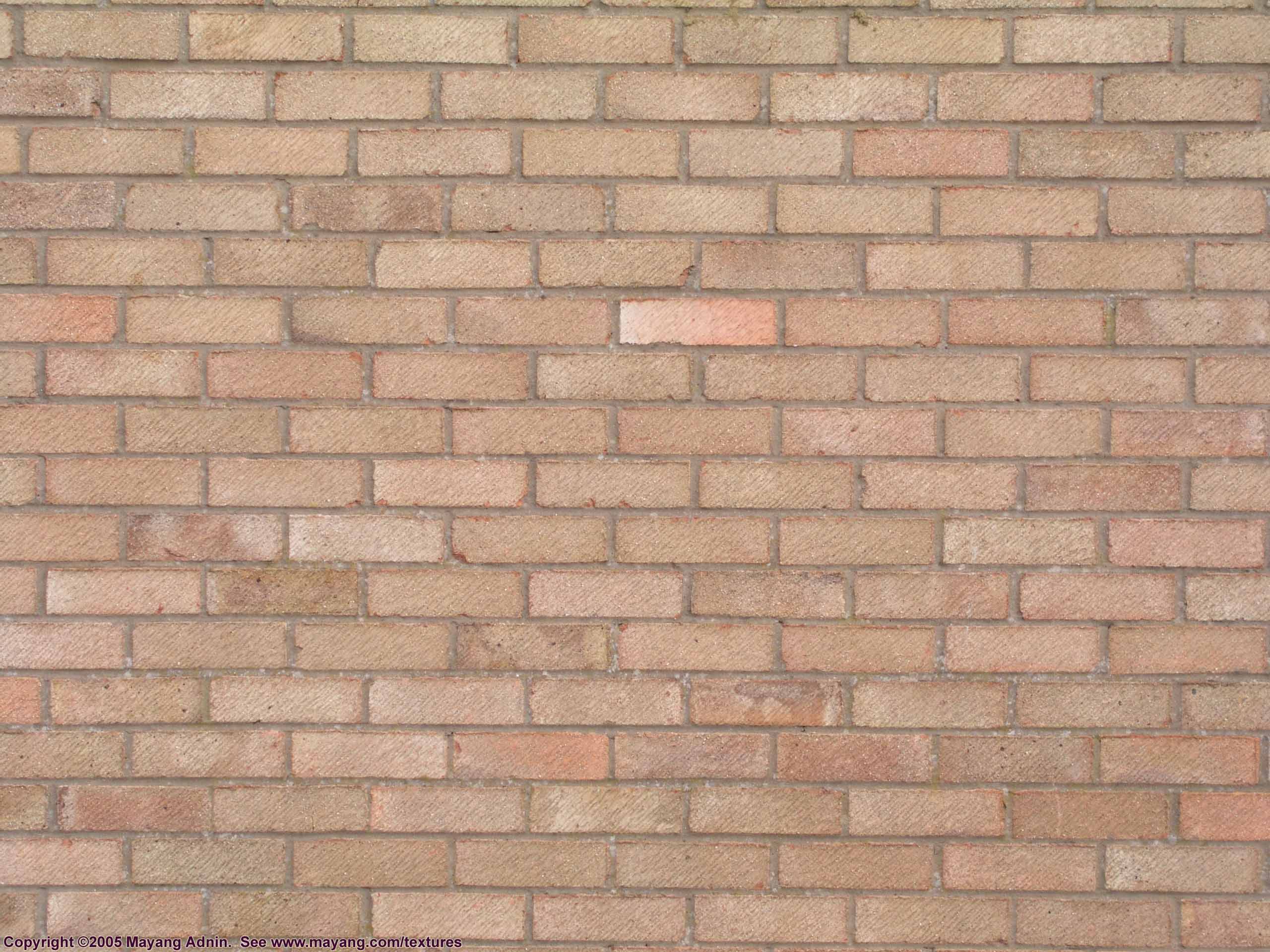}&
\includegraphics[width=3.2cm,bb=0 0 2560 1920]{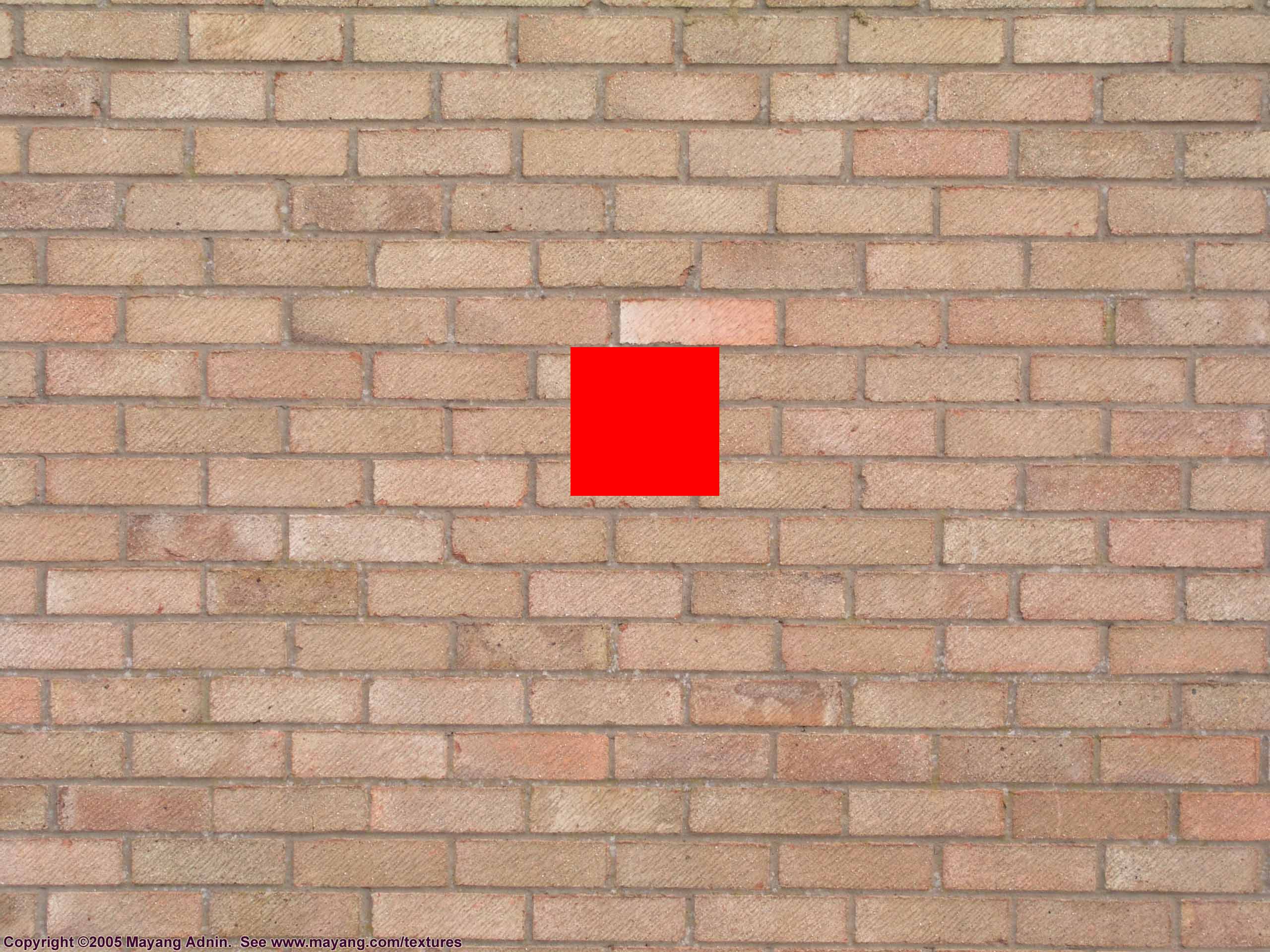}&
\includegraphics[width=3.2cm,bb=0 0 2560 1920]{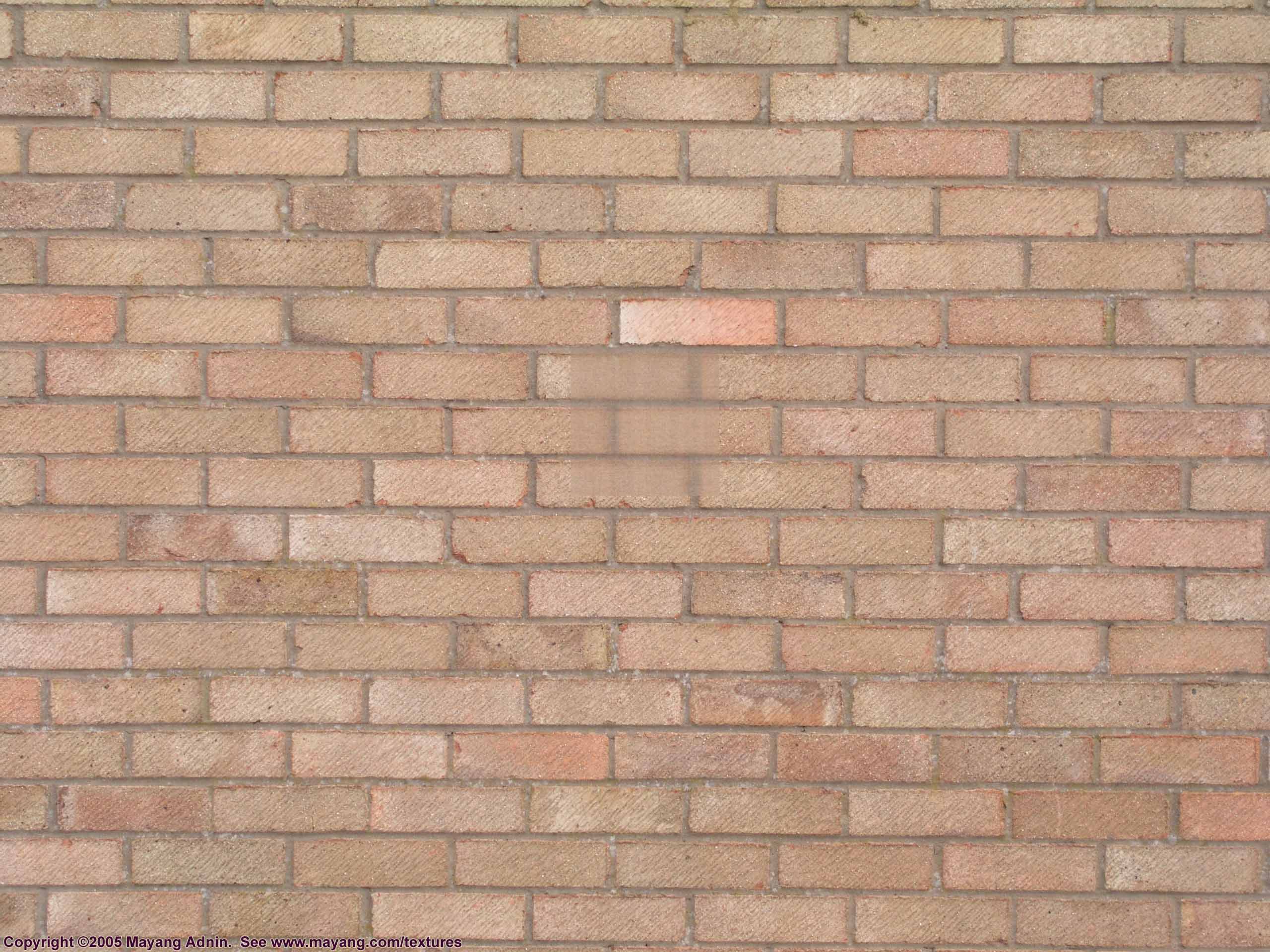}&
\includegraphics[width=3.2cm,bb=0 0 2560 1920]{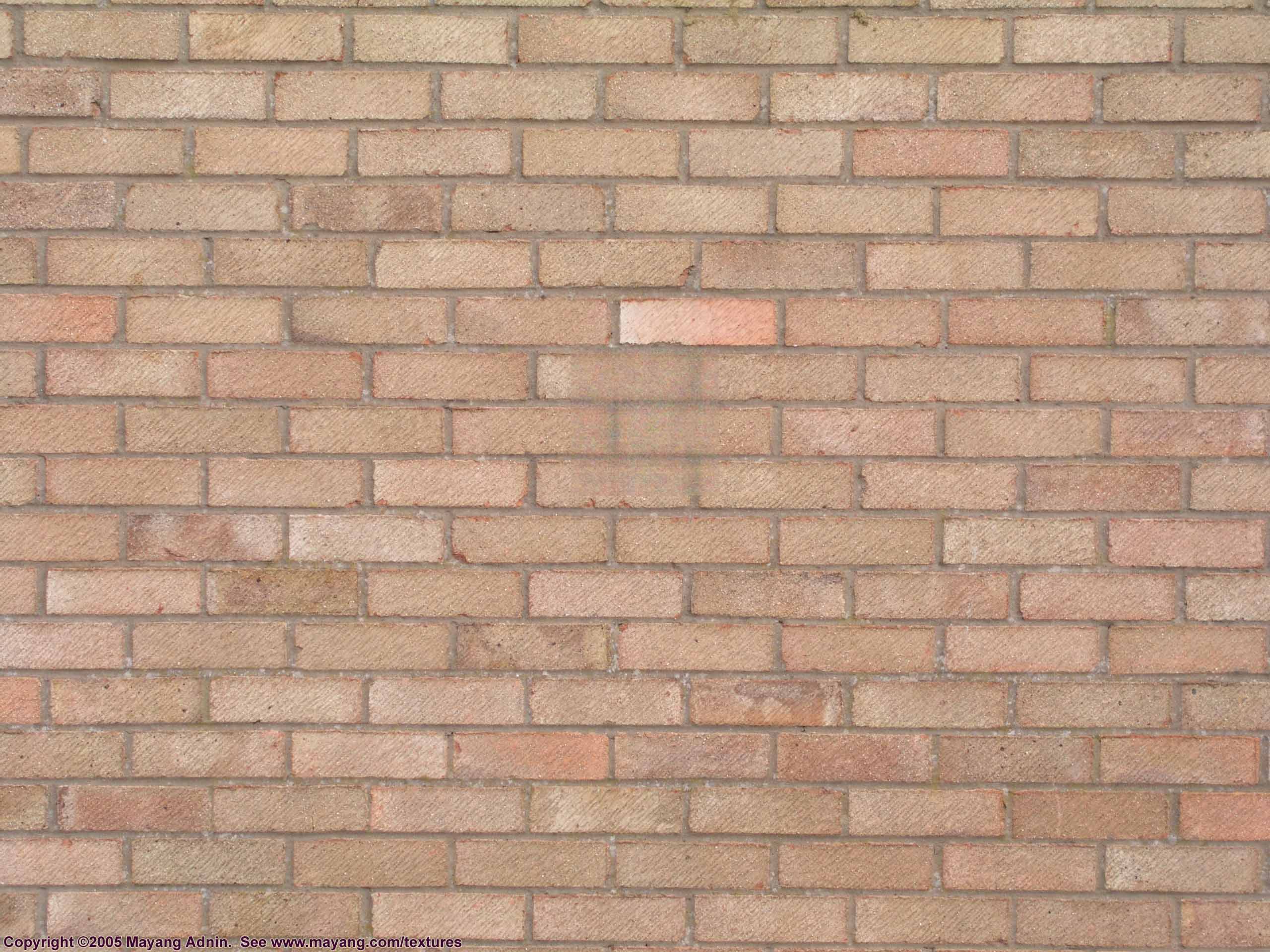}&
\includegraphics[width=3.2cm,bb=0 0 2560 1920]{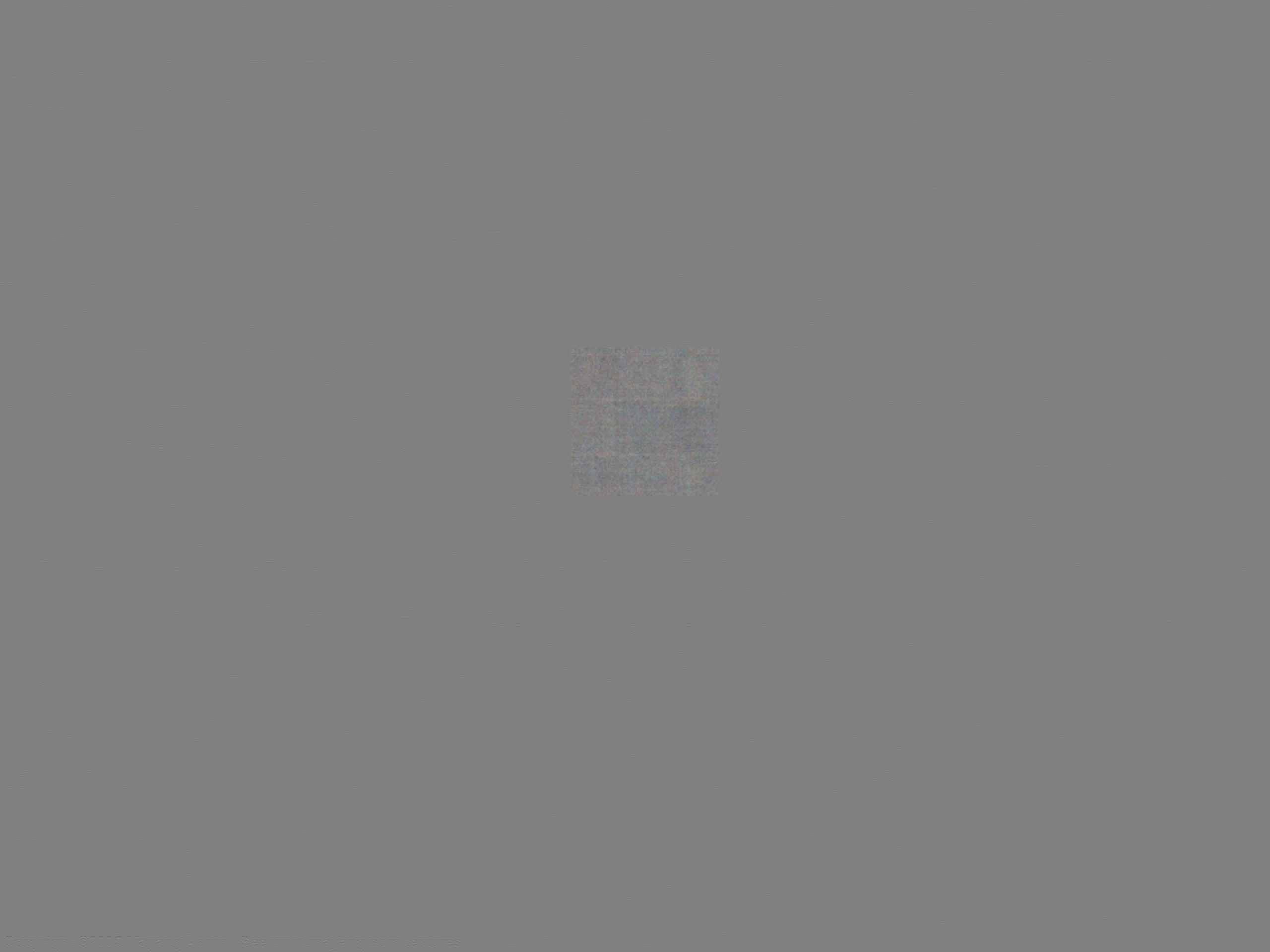}
\\
(a) & (b) & (c) & (d) & (e)
\\
\includegraphics[width=3.2cm,bb=0 0 2560 1920]{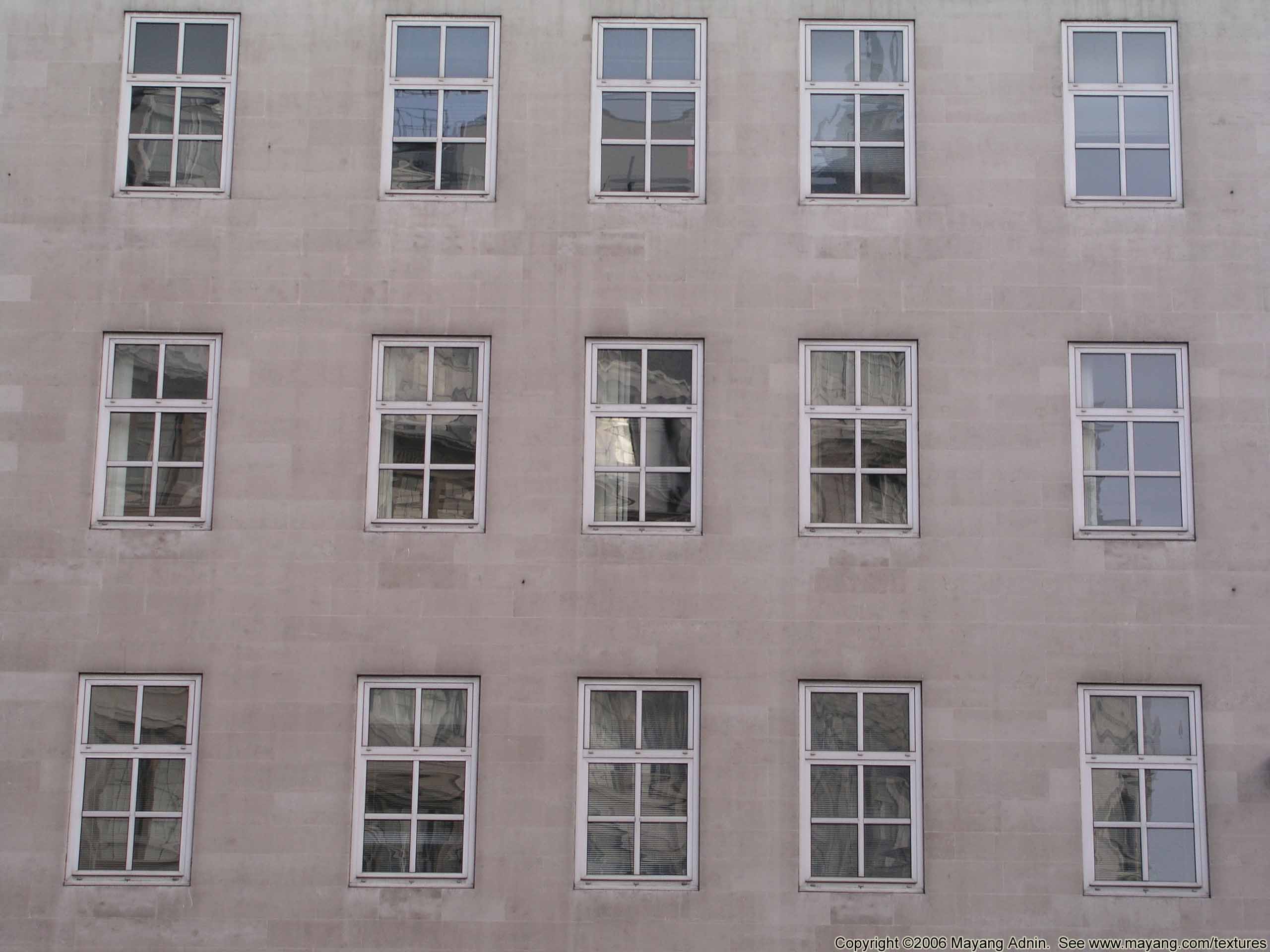}&
\includegraphics[width=3.2cm,bb=0 0 2560 1920]{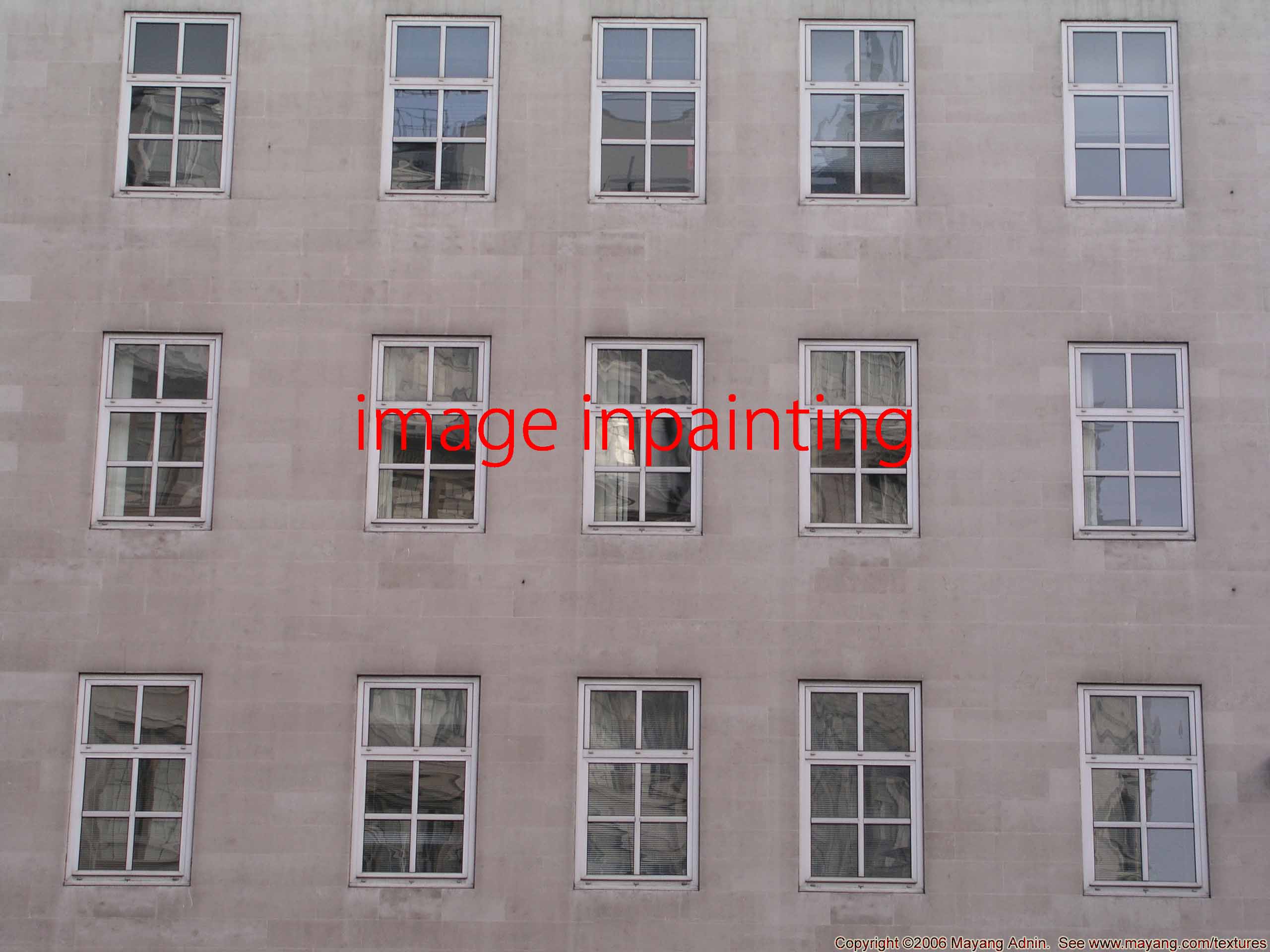}&
\includegraphics[width=3.2cm,bb=0 0 2560 1920]{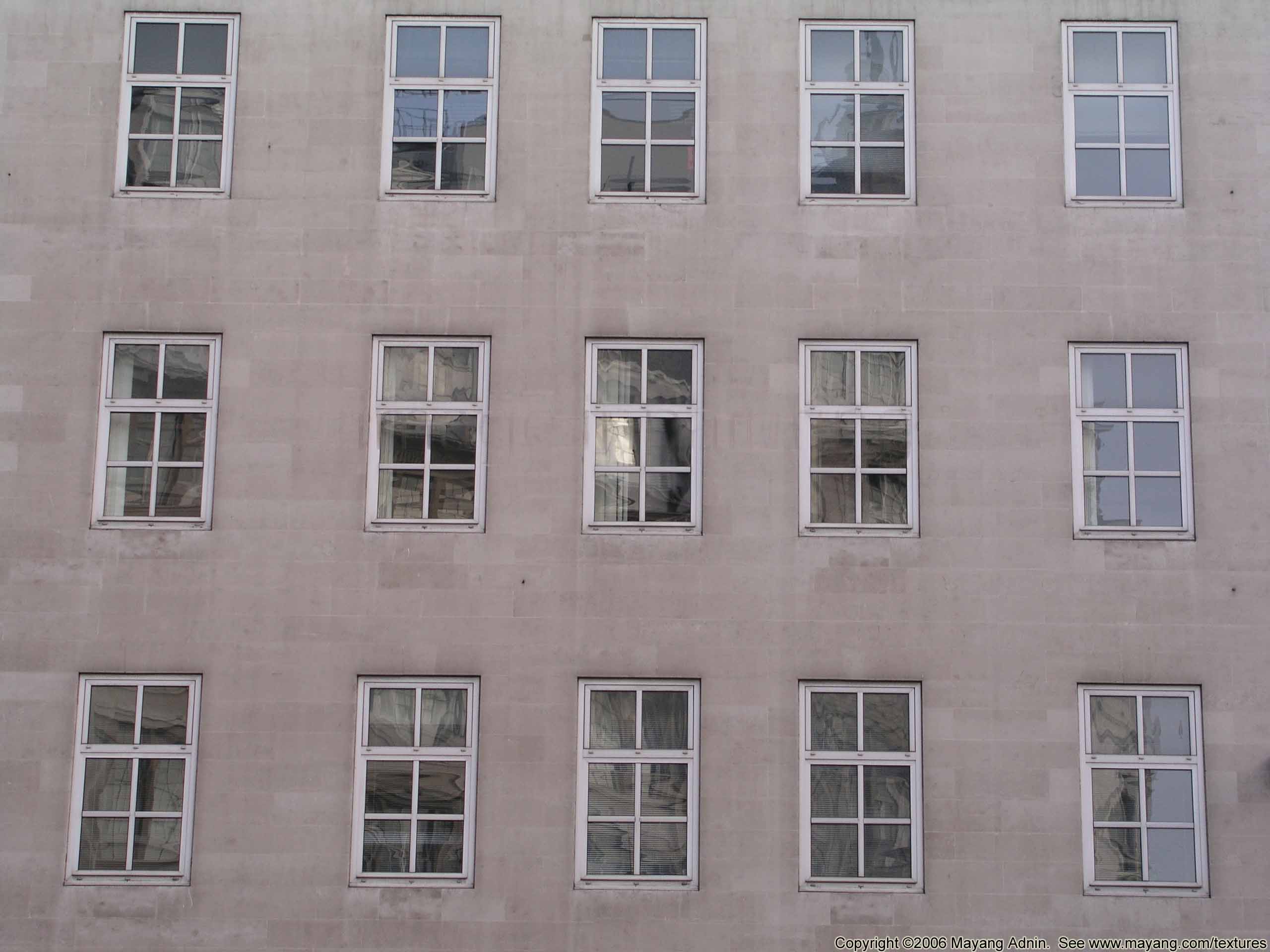}&
\includegraphics[width=3.2cm,bb=0 0 2560 1920]{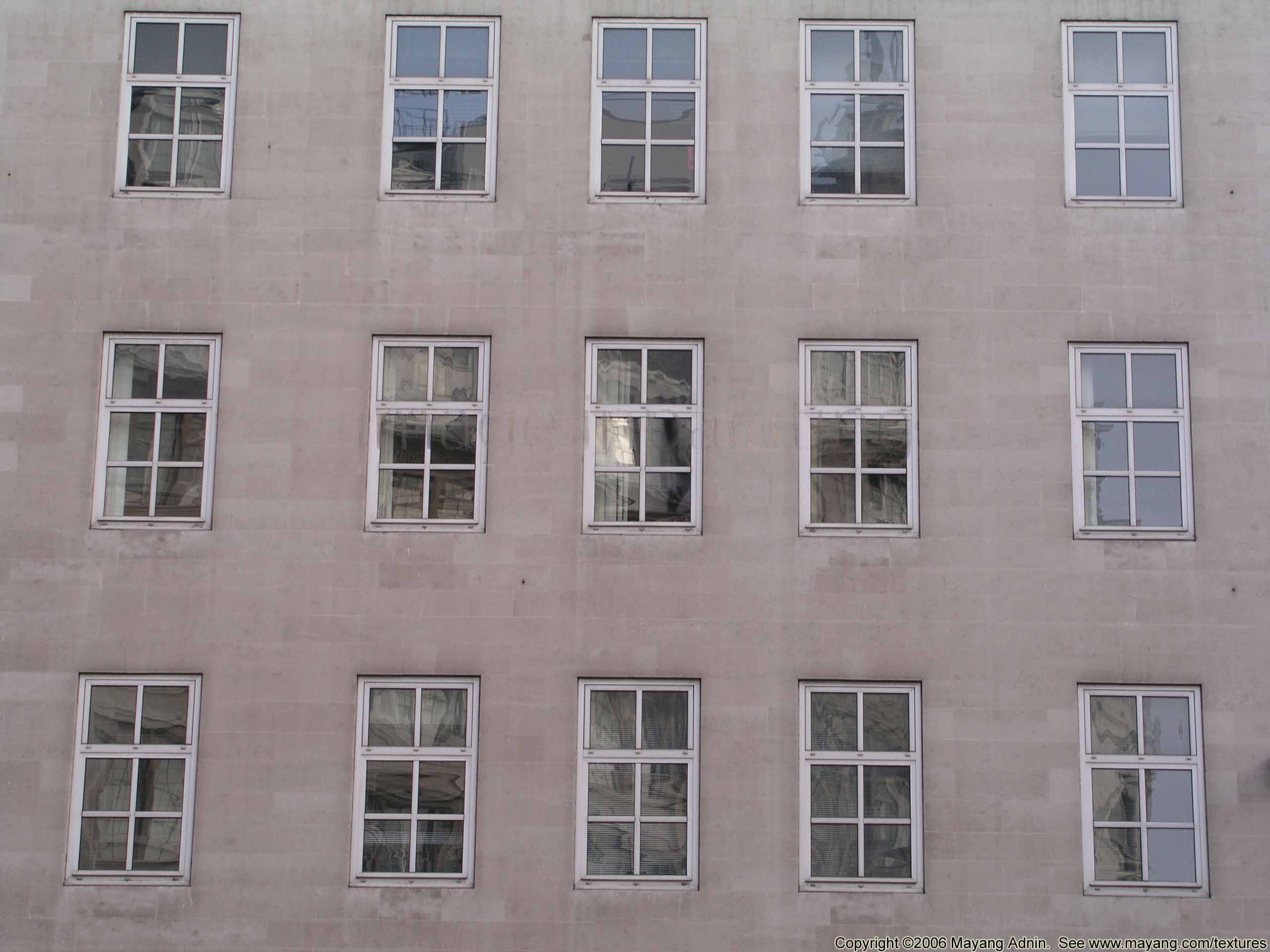}&
\includegraphics[width=3.2cm,bb=0 0 2560 1920]{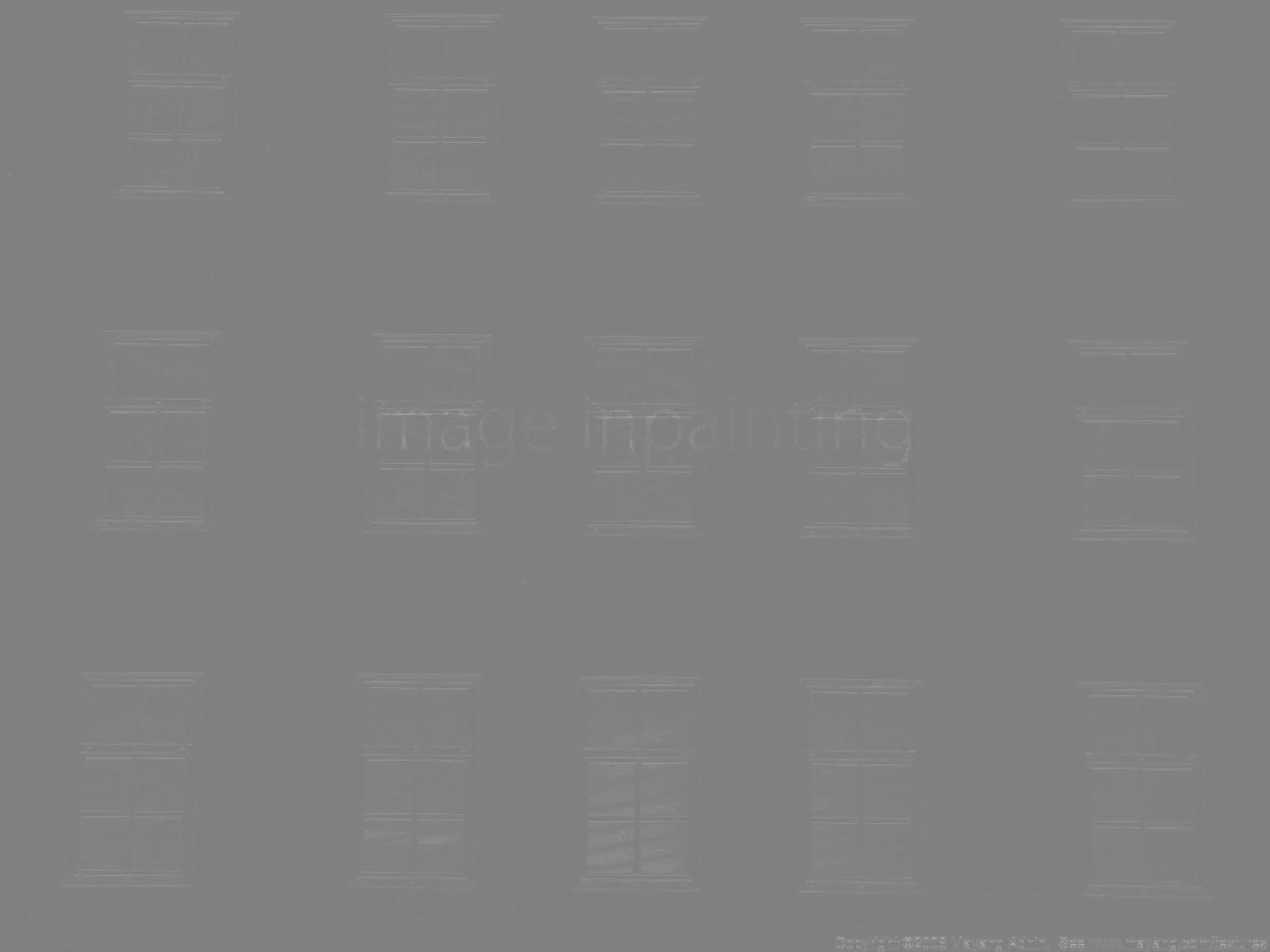}
\\
(f) & (g) & (h) & (i) & (j)
\end{tabular}
\caption{Image inpainting results with $20$th-order approximation.
(a) and (f) Original images \textit{Bricks} and \textit{Office windows}.
(b) and (g) Observed image with missing region.
(c) and (h) SVD-based method.
(d) and (i) Proposed method.
(e) Difference between (c) and (d).
(j) Difference between (h) and (i).
In these difference images, $0.5$ is added to all pixel values.}
\label{fig:image_inpaint}
\end{figure*}

The objective with this application is to recover a missing region (as shown in the later Fig.~\ref{fig:image_inpaint}(b)).

Let $\mathbf{L}$ and $\mathbf{I} \!\in\! \mathbb{R}^{m\times n}$ be a texture image and a given image with missing regions, respectively.
Then, let $\Omega$ and $\overline{\Omega}$ be observed and missing regions and $P_{\Omega}(\cdot)$ and $P_{\overline{\Omega}}(\cdot)$ be linear operators extracting pixels in their regions.
From the notations, the missing region is represented as $\mathbf{M}\!=\!P_{\overline{\Omega}}(\mathbf{L})$.
The pixels surrounding $\mathbf{M}$ with the size of five pixels, as shown in Fig.~\ref{fig:t-t_near}, are defined as $\mathbf{M}_\partial$.
Let $\mathbf{T}_1 \!\in\! \mathbb{R}^{m\times m}$ and $\mathbf{T}_2 \!\in\! \mathbb{R}^{n\times n}$ be the DCT matrices in the horizontal and vertical matrix directions, i.e., these matrices transform an image to its frequency domain.
Since a regular texture image is basically sparse in its frequency domain, it can be represented as $\mathbf{L}\!=\!\mathbf{T}_1 \mathbf{S} \mathbf{T}_2^\top$, where $\mathbf{S}$ is the coefficients on the frequency domain of $\mathbf{L}$.
Additionally, the set of a normalized dynamic range constraint is defined as $\mathcal{D} \!:=\! \{ \mathbf{x}\!:=\![x_i]_{i=1}^{mn}  |~x_i \!\in\! [0,1] \}$.
When $\mathbf{L}$ and $\mathbf{S}$ are assumed to be low rank and sparse, the reconstruction problem can approximately be solved using the nuclear norm\footnote{The nuclear norm of $\mathbf{X}\!\in\! \mathbb{R}^{m\times n}$ is defined as $\| \mathbf{X} \|_* := \sum^K_{i=1} \sigma_i(\mathbf{X})$, where $i \!\in\! \{ 1,2,\ldots,K \} (K := \min (m, n))$.} and the $\ell_1$ norm as
\begin{table*}[htb]
  \begin{center}
    \caption{Total and Average Computation Time and RMSEs (Image Inpainting Shown in Fig.~\ref{fig:image_inpaint})}
    \label{computation_time_psnr_inpaint}
    \begin{tabular}{c||c|c|c|c||c}\hline
    \multicolumn{6}{c}{\textit{Bricks}}\\ \hline
    Approximation order& 5 & 10&15&20&EVD-based method\\ \hline \hline
     Total computation time (s) &401.85&373.91&347.53&348.38&464.81 \\ \hline
     Average computation time of singular value soft-shrinkage (s) &1.35&1.64&1.81&2.00&2.48 \\ \hline\hline
     RMSE between our method and SVD-based method ($\times$10$^{-3}$)&9.92&9.83&9.77&9.73&-\\\hline
     \end{tabular}
     \\[6pt]
     \begin{tabular}{c||c|c|c|c||c}\hline
     \multicolumn{6}{c}{\textit{Office windows}}\\ \hline
      Approximation order& 5 & 10&15&20&EVD-based method\\ \hline \hline
     Total computation time (s) &418.50&384.95&366.74&369.93&464.70 \\ \hline
     Average computation time of singular value soft-shrinkage (s) &1.39&1.59&1.85&1.87&2.48 \\ \hline\hline
     RMSE between our method and SVD-based method ($\times$10$^{-3}$)&8.73&8.75&8.75&8.72&-\\\hline
    \end{tabular}
  \end{center}
\end{table*}
\begin{equation}
\begin{aligned}
	& \min_{\mathbf{L}, \mathbf{S},\mathrm{vec}(\mathbf{L})\in \mathcal{D}} \ \| \mathbf{L}\|_{*}+\eta  \| \mathbf{S} \|_1
	\\
	&
	\begin{aligned}
	\text{s.t.} \quad & P_{\Omega}(\mathbf{I}) = P_{\Omega}(\mathbf{L}), \quad \mathbf{L}=\mathbf{T}_1 \mathbf{S} \mathbf{T}_2^\top,\\
	& \mathrm{ave}\bigl( \mathrm{vec}(\mathbf{M}) \bigr) = \mathrm{ave}\bigl( \mathrm{vec}(\mathbf{M}_\partial) \bigr),
	\end{aligned}
\end{aligned}
\label{texture_repaire4}
\end{equation}
where a positive real value $\eta$ is a regularization parameter, $\mathrm{ave}(\cdot)$ calculates the arithmetic average, and $\mathrm{vec}(\cdot)$ is the operator vectorizing a matrix.
The average pixel value on the recovered region is assumed to be identical to that of its surrounding pixel values.
As can be seen, \eqref{texture_repaire4} is composed of nuclear norm relaxation so that our method can be used for its efficient calculation.
Hereafter, we discuss the validation of our method by applying the ADMM to \eqref{texture_repaire4} to obtain the optimal solution.
In Appendix~\ref{sec:admm_inpaint}, \eqref{texture_repaire4} is converted to the form to which the ADMM is applicable.

Eight-bit color images \textit{Bricks} and \textit{Office widows}\footnote{The images are available at \url{http://www.mayang.com/textures/}.}, as shown in Figs.~\ref{fig:image_inpaint}(a) and (f), were used for the application, where the size of each color component of the images was 1920 $\times$ 2560.
The pixel values of each color component of the images were in the range from 0 to 1.
In the application, each color component was inpainted separately.
The observed image $\mathbf{I}\!\in\! \mathbb{R}^{2560\times 1920}$ with the missing regions were defined, as shown in Fig.~\ref{fig:image_inpaint}(b) and (g)\footnote{In the experiment, the images are transposed to ``portrait''.}.
The number of missing pixels in Fig.~\ref{fig:image_inpaint}(b) was $\mathbf{M}\!\in\! \mathbb{R}^{300 \times 300}$ and that in Fig.~\ref{fig:image_inpaint}(g) was $\mathbf{M}\!\in\! \mathbb{R}^{180\times 1250}$.
In the ADMM applicable form (see Appendix~\ref{sec:admm_inpaint}), the column vectors $\mathbf{l}_0$ and $\mathbf{u}_0$ were initialized by all-one vectors.
Additionally, for $\mathrm{prox}_{1/\rho\| \cdot \|_*}$ and $\mathrm{prox}_{\eta/\rho\| \cdot \|_1}$, the thresholding parameters $(1/\rho,\eta/\rho)$ were set to $(6, 0.1)$ in \textit{Bricks} and $(5, 0.1)$ in \textit{Office windows}, where the parameters were determined for the fast and stable convergence of the optimization.
For fast computation, parallel processing\footnote{The MATLAB function \textit{parfor}, which is contained in the parallel computing toolbox, was used for the parallel computing only in the image inpainting method.} was performed in the application of the color components.

Figure~\ref{fig:image_inpaint} shows the results of image inpainting with the 20th-order approximation.
The resulting image recovered using our method was practically equivalent to that with the SVD-based method by comparing Figs.~\ref{fig:image_inpaint}(c), (d), (h), and (i).
In Figs.~\ref{fig:image_inpaint}(e) and (j), it is clear that the exact and approximated solutions had little differences, which visually indicates the high approximation precision of our method.

Table~\ref{computation_time_psnr_inpaint} lists the computation time and RMSE comparisons.
Our method was faster than the EVD-based method while maintaining reconstruction performance.
Regarding the total computation times of \textit{Bricks} and \textit{Office windows}, our method with the 5th-order and 10th-order approximations was slower than with the 15th-order approximations.
This is because our method with the low-order approximations did not converge well due to the low approximate precision.

%
%
\subsection{Background Modeling of Video\cite{background_mod1,background_mod2,background_mod3,background_mod4,background_mod5}}

\begin{figure*}[htb]
\small
\begin{tabular}{c|c|c}\hline
Original and difference&CPA-based method&EVD-based method\\ \hline \hline
\begin{minipage}[h]{0.33\linewidth}
  \centering
\centerline{\includegraphics[width=3.5cm,bb=0 0 360 240]{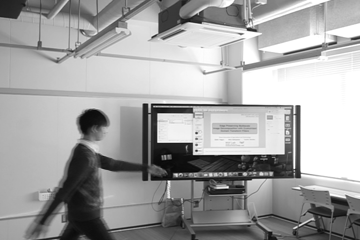}}
    \centerline{Original frame}\medskip
\end{minipage}
&
\begin{minipage}[h]{0.3\linewidth}
\centering
\includegraphics[width=3.5cm,bb=0 0 360 240]{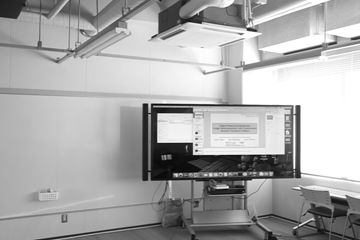}
\centerline{(a) Low rank}\medskip
 \end{minipage}
 &
\begin{minipage}[h]{0.3\linewidth}
  \centering
\centerline{\includegraphics[width=3.5cm,bb=0 0 360 240]{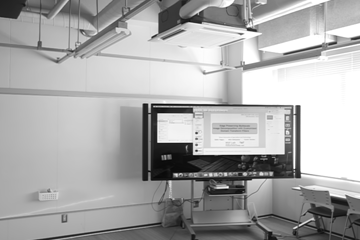}}
    \centerline{(c) Low rank}\medskip
\end{minipage}
\\
\begin{minipage}[h]{0.3\linewidth}
\centering
\includegraphics[width=3.5cm,bb=0 0 360 240]{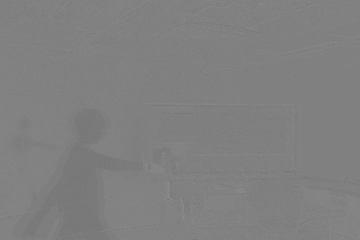}
\centerline{Difference between (a) and (c)}\medskip
 \end{minipage}
&
\begin{minipage}[h]{0.3\linewidth}
\centering
\includegraphics[width=3.5cm,bb=0 0 360 240]{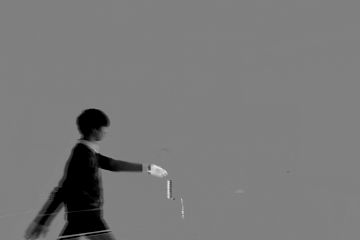}
\centerline{(b) Sparse}\medskip
 \end{minipage}
&
\begin{minipage}[h]{0.3\linewidth}
  \centering
\centerline{\includegraphics[width=3.5cm,bb=0 0 360 240]{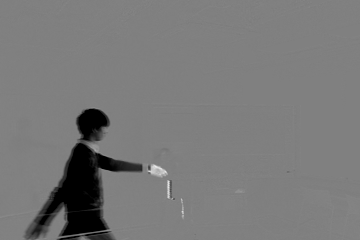}}
    \centerline{(d) Sparse}\medskip
\end{minipage}\\ \hline
\end{tabular}
\caption{Background modeling results with 20th-order approximation.
2$\times$ amplified difference between (a) and (c) is also shown, where 0.5 is added to all pixel values.}
\label{fig:back_modeling}
\end{figure*}

The objective with this application is to divide a video sequence into background and object sequences (as shown in Fig.~\ref{fig:back_modeling}(a) and (b)).

Let $\mathbf{I}^{(i)}\!\in\! \mathbb{R}^{m\times n}$ be the $i$-th frame of a video sequence.
The sequence is rearranged into a matrix $\mathbf{I}\!\in\! \mathbb{R}^{mn\times K}$ as
\begin{equation}
\mathbf{I} :=
\begin{bmatrix}
\mathrm{vec}(\mathbf{I}^{(1)})~\mathrm{vec}(\mathbf{I}^{(2)})~\dots~ \mathrm{vec}(\mathbf{I}^{(K)})
\end{bmatrix}.
\label{observed_sequences}
\end{equation}
Then, let $\mathbf{L}$ and $\mathbf{S} \!\in\! \mathbb{R}^{mn\times K}$ be the background sequence and sequence of moving objects of a video.
In $\mathbf{L}$, pixel values corresponding to $\mathbf{S}$ are zero and vice versa.
The background and moving objects can be assumed to be low rank and sparse; hence, the background modeling is solved as the following convex optimization problem:
\begin{equation}
\min_{\mathbf{L}, \mathbf{S}} \ 
\| \mathbf{L} \|_* + \eta \| \mathbf{S} \|_1\quad 
\text{s.t.} \quad 
\mathbf{I} =  \mathbf{L} + \mathbf{S}.
\label{convex_opt_back_mod}
\end{equation}
Problem~\eqref{convex_opt_back_mod} is also solved using the ADMM.

In the background modeling, an eight-bit grayscale video \textit{Laboratory}\footnote{This video was recorded with our video camera. It was downsampled and transformed into grayscale for the experiment.} was used.
The $\mathbf{I}^{(i)}\!\in\! \mathbb{R}^{360\times 240}$ was the $i$-th frame of the video in $i\!\in\! \{ 1,2,\ldots,5000 \}$; hence, the matrix of the sequence was $\mathbf{I} \!\in\! \mathbb{R}^{86400\times 5000}$.
The pixel values of the video were in the range from 0 to 1.
In the ADMM applicable form (see Appendix~\ref{sec:admm_background}), the column vectors $\mathbf{l}_0$, $\mathbf{s}_0$, and $\mathbf{u}_0$ were initialized by all-one vectors.
\begin{table*}[htb]
  \begin{center}
    \caption{Total and Average Computation Time and RMSEs (Background Modeling Shown in Fig.~\ref{fig:back_modeling})}
    \label{computation_time_psnr_background}
    \begin{tabular}{c||c|c|c|c||c} \hline
    Approximation order& 5 & 10&15&20&EVD-based method\\ \hline \hline
     Total computation time (s) &1200.21&1127.60&1150.27&2192.24&6457.53 \\ \hline
     Average computation time of singular value soft-shrinkage (s) &20.95&21.21&21.79&22.20&30.77 \\ \hline\hline
     RMSE between our method and EVD-based method ($\times$10$^{-3}$)&56.27&23.89&10.86&3.71&-\\\hline
    \end{tabular}
  \end{center}
\end{table*}Additionally, the thresholding parameters $(1/\rho,\eta/\rho)$ were set to $(480,0.12)$, where the parameters were determined for the fast and stable convergence of the optimization.
For comparison, low rank and sparse components in the 160th frame are shown in Fig.~\ref{fig:back_modeling} with the 20th-order approximation.

Our method effectively decomposed the video sequences to low rank and sparse sequences, as shown in Figs.~\ref{fig:back_modeling}(a) and (b).
They are almost equivalent to those with the EVD-based method; hence, the difference between the low rank images are not displayed even though the difference is amplified (bottom left of Fig.~\ref{fig:back_modeling}).

Table~\ref{computation_time_psnr_background} summarizes the computation times and RMSEs between the background modeling of our method and that of the EVD-based method\footnote{The SVD-based method ran out of memory in our machine so that the results of the EVD-based method were used for the RMSEs.}.
Our method was sufficiently faster than the EVD-based method in all approximation orders.
\begin{table*}[htb]
  \begin{center}
    \caption{Total and Average Computation Time and RMSEs of Proposed and Existing Methods.}
    \label{computation_time_dct_wav_blk}
    \begin{tabular}{c||c|c|c} \hline
    \multicolumn{4}{c}{Proposed method (see Section~\ref{sec:influence} for the explanation)}\\ \hline
    Used transformation methods of CPA-based method&DWT& Block DCT&DCT\\ \hline \hline
     Total computation time (s) &347.53&344.03&380.23\\ \hline
     Average computation time of singular value soft-shrinkage (s) &1.81&1.06&1.33 \\ \hline\hline
     RMSE between our method and SVD-based method ($\times$10$^{-3}$)&9.77&3.82&3.81\\\hline
     \end{tabular}
     \\[6pt]
     \begin{tabular}{c||c|c|c|c} \hline
     \multicolumn{5}{c}{Existing method (see Section~\ref{sec:comparison_existing_meth} for the explanation)}\\ \hline
    Used algorithms&Exact PSVD&FRSVS \cite{SVD_echon2}&NSVS\cite{SVD_echon3}&FSVS\cite{Fast_singular_thresh}\\ \hline \hline
     Total computation time (s) &1935.73&334.59&336.22&1926.12 \\ \hline
     Average computation time of singular value soft-shrinkage (s) &6.70&1.63&1.67&21.26 \\ \hline\hline
     RMSE between existing methods and SVD-based method ($\times$10$^{-3}$)&12.17&13.94&117.61&26.40\\\hline
    \end{tabular}\\
    \end{center}
\end{table*}
\begin{table*}[htb]
  \begin{center}
    \caption{Comparisons of DCT with Adaptive and Fixed $k$. Recall that $k_1=(0.5\!\times\!10^{-4})n^2$, $k_2=(1.5\!\times\!10^{-4})n^2$, $k_3=(0.5\!\times\!10^{-3})n^2$, and $k_4=n^2$.}
    \label{computation_time_adaptive_DCT}
    \begin{tabular}{c||c||c|c|c|c} \hline
                 Methods&DCT with adaptive $k$&\multicolumn{3}{c}{DCT with fixed $k$}\\ \hline
    Used $k$ for $\mathcal{K}(|\underline{\mathbf{\Phi}}|,k)$& Defined in \eqref{recommendation_epsilon}&$k_1$&$k_2$&$k_3$&$k_4$\\ \hline \hline
     Total computation time (s) &358.16&Not converged&416.52&380.23&765.84\\ \hline
     Average computation time of singular value soft-shrinkage (s) &1.26&1.25&1.26&1.33&5.92\\ \hline\hline
     RMSE between our method and SVD-based method ($\times$10$^{-3}$)&3.81&-&3.82&3.81&3.81\\\hline
    \end{tabular}
  \end{center}
\end{table*}This is because the reduction in computational complexity by thresholding the transformed coefficients described in \eqref{approx_eigen_shrink_trans_coeff} is effective for our method to have low computational complexity while retaining high approximation precision.
However, our method with the 5th-order approximations took more time than that with 10th-order approximation because it did not converge well due to its low approximate precision.

%
%
\subsection{Effects of Selections: Transform Matrix and Thresholding Value}
\label{sec:influence}

The selections of a transform matrix $\mathbf{T}$ in \eqref{SVD_fil} and a thresholding value $\varepsilon$ in \eqref{threshold_components} affect computation time and size of approximation error.
In this subsection, we indicate these effects experimentally by using the image inpainting method for \textit{Bricks}.
In all experiments, the 15th-order approximation was used for our method.

%
%
\subsubsection{Effects of Selected Transform Matrix}
We compared the DWT with the DCT and the block diagonal forms of the DCT (block DCT) \cite{kaiser1} whose block size was $8\times 8$ for indicating the differences among chosen transform matrices $\mathbf{T}$ in \eqref{SVD_fil}.
The threshold values in \eqref{threshold_components} were fixed to $\varepsilon \!=\! 250$ for the DCTs.
The other experimental conditions were the same as those discussed in Section~\ref{subsec:texture_image_inpainting}.

The results of the proposed method in Table~\ref{computation_time_dct_wav_blk} show the performance comparisons.
Our method with the DWT is as fast as that with the block DCT in the total computation time, though the singular value shrinkage with our method with the DWT takes more time than the others.
%
%
This is because the maximum iteration of our method with the DWT is only 83, whereas those of our method with the block DCT and DCT are 101 and 97, respectively.
Therefore, the DWT leads our method to be stable convergence.
However, our method with the DWT indicates a higher RMSE than the others since all the high frequency components were removed in it.
The fact is not fatal problem because the images derived by using our method with the DWT were very similar to that of the exact one as shown in the previous section.
Additionally, the block DCT can substantially sparsify an image compared to the DWT and DCT.
Therefore, our method with the block DCT is faster than the others.
In spite of fast computation, the results with the block DCT shows an RMSE as low as that of the DCT.

%
%
\subsubsection{Effects of Thresholding Value}

When $\varepsilon$ is an excessively large value, our method becomes fast but matrix rank minimization cannot converge well due to the errors w.r.t. the reduction in the number of components.
For achieving the fast computation and stable convergence, we present a recommended guideline on the thresholding values.

Let $E_{t}$ be an error at the $t$-th iteration of the ADMM, i.e., $E_{t} \!=\! \| \mathbf{l}_{t} - \mathbf{l}_{t-1} \| / \| \mathbf{l}_{t} \|$ from Appendix~\ref{sec:admm_inpaint} and $\mathcal{K}(\mathbf{X},k)$ be a function that returns the $k$-th largest element in a matrix $\mathbf{X}$, where  $\mathcal{K}(\mathbf{X},k)$ is used as the threshold.
Basically, when $E_t$ is a large value, a small $k$ does not have any problem to decrease the error.
In contrast, when $E_t$ is a small value, i.e., the optimization almost converges, $k$ should be large for a stable convergence.
For this purpose, we recommend the thresholding percentage for $\mathcal{K}(|\underline{\mathbf{\Phi}}|,k)$ as
\begin{equation}
\varepsilon(E_t)=
\begin{cases}
\mathcal{K}\bigl(|\underline{\mathbf{\Phi}}|,(0.5\!\times\!10^{-4})n^2\bigr) & \text{if}~E_{t} \!>\! E_\ell,\\
\mathcal{K}\bigl(|\underline{\mathbf{\Phi}}|,(1.5\!\times\!10^{-4})n^2\bigr) & \text{if}~E_\ell \!\ge\! E_{t} \!>\! E_\mathrm{m},\\
\mathcal{K}\bigl(|\underline{\mathbf{\Phi}}|,(0.5\!\times\!10^{-3})n^2\bigr) & \text{otherwise},
\end{cases}
\label{recommendation_epsilon}
\end{equation}
\begin{figure}[tb]
\small
\tabcolsep = 0.2mm
\centering
\begin{tabular}{cc}
\includegraphics[trim=0 0 50 0, width=0.3\linewidth, clip]{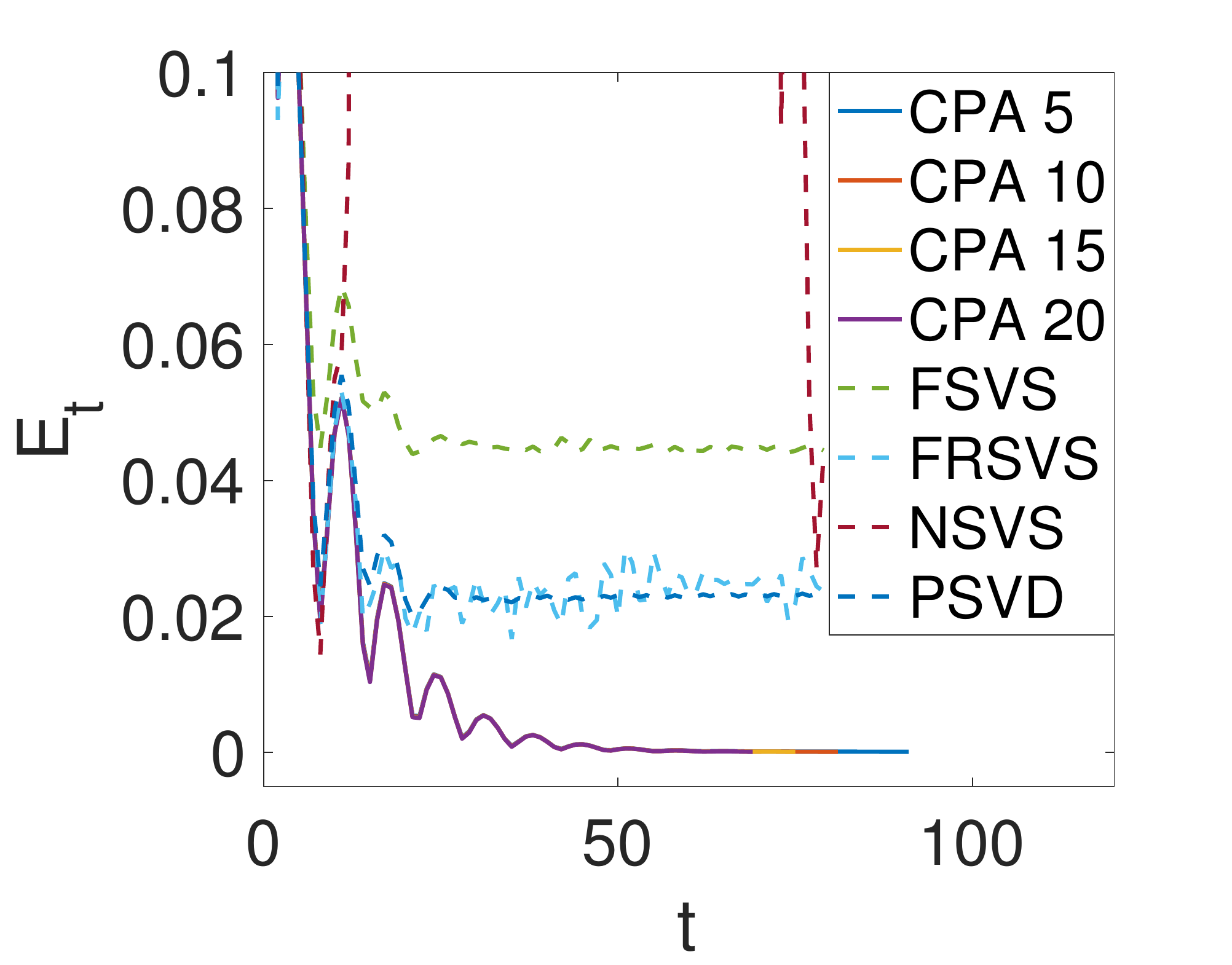}&
\includegraphics[trim=0 0 50 0, width=0.3\linewidth, clip]{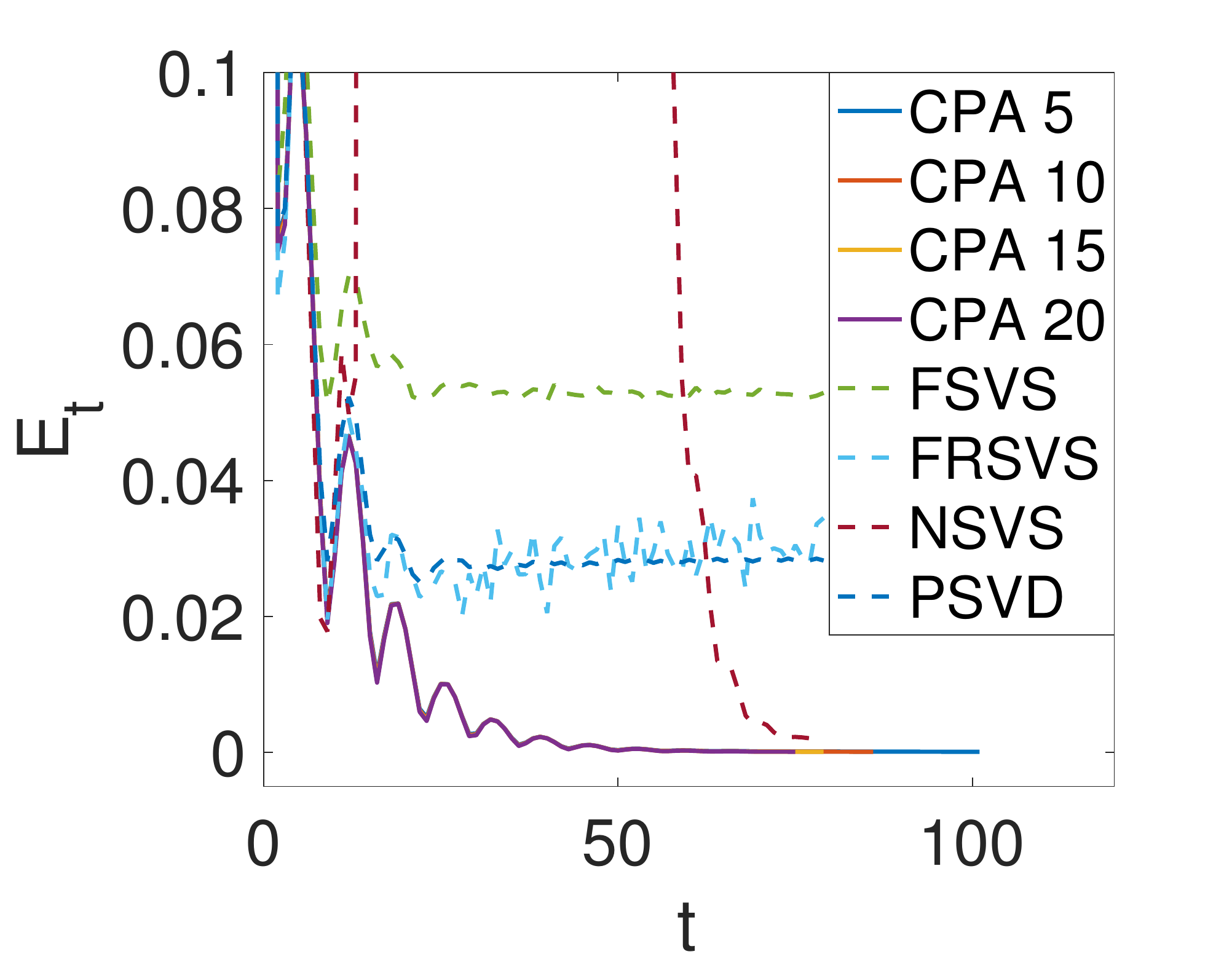}
\\
(a) & (b)
\\
\multicolumn{2}{c}{\includegraphics[trim=0 0 50 0, width=0.3\linewidth, clip]{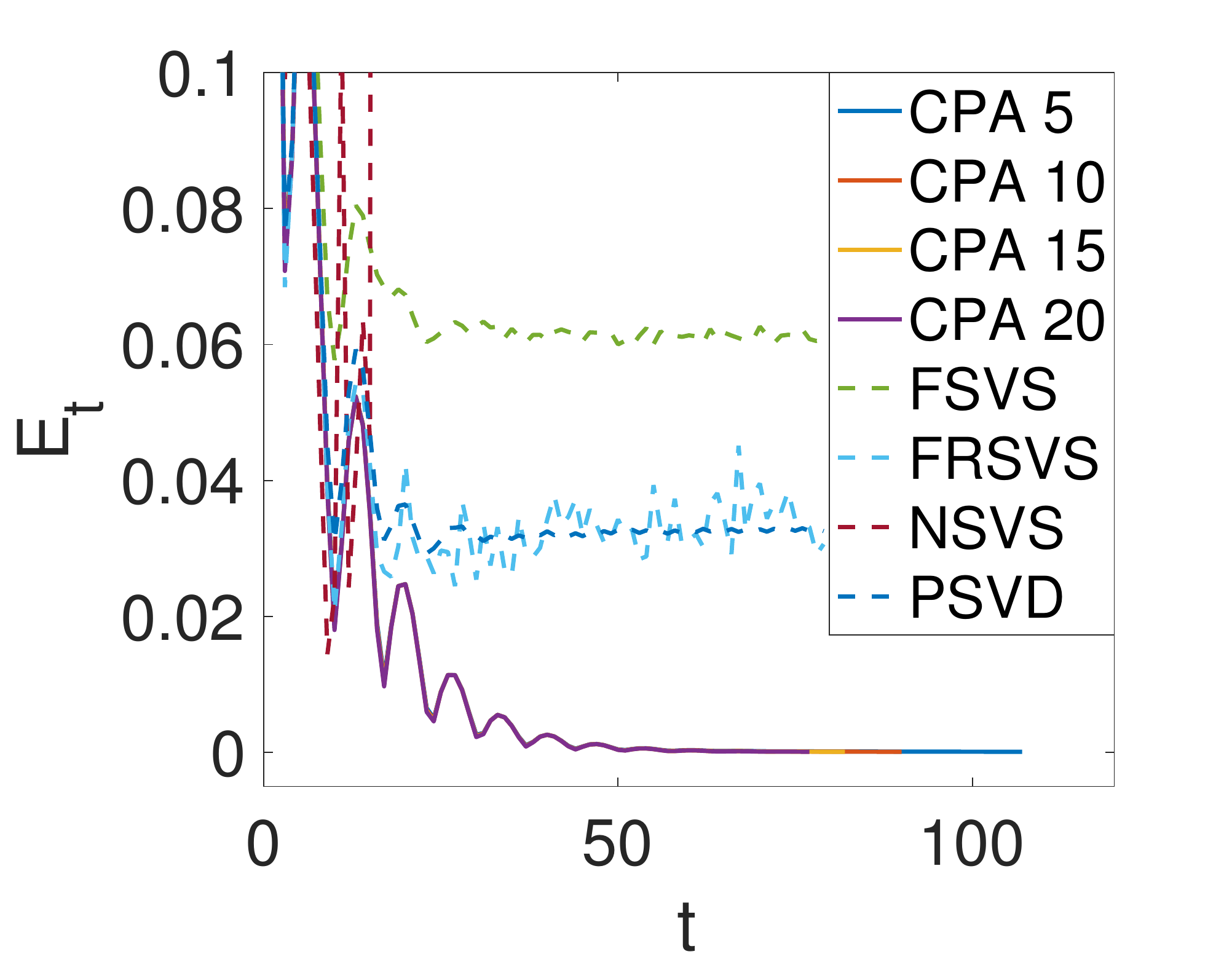}}
\\
\multicolumn{2}{c}{(c)}
\end{tabular}
\caption{Figures represent $E_t$ for \textit{Bricks}, whose (a), (b) and (c) are for the component of red, green and blue, respectively.
Maximum number of iterations for CPA-based method is indicated as follows.
5th-order approximation: 108th iteration;
10th-order approximation: 91st iteration;
15th-order approximation: 83rd iteration;
and 20th-order approximation: 78th iteration.}
\label{fig:Comparison_results_Error_RMSE}
\end{figure}
\begin{figure*}[htb]
\small
\begin{tabular}{c||c|c|c|c}
\tiny{Reduction rate\textbackslash Rank}&10&100&200&500\\\hline \hline
1$\%$&
\begin{minipage}[h]{0.195\linewidth}
  \centering
  \centerline{\includegraphics[trim=3mm 3mm 3mm 3mm, clip,width=3.5cm]{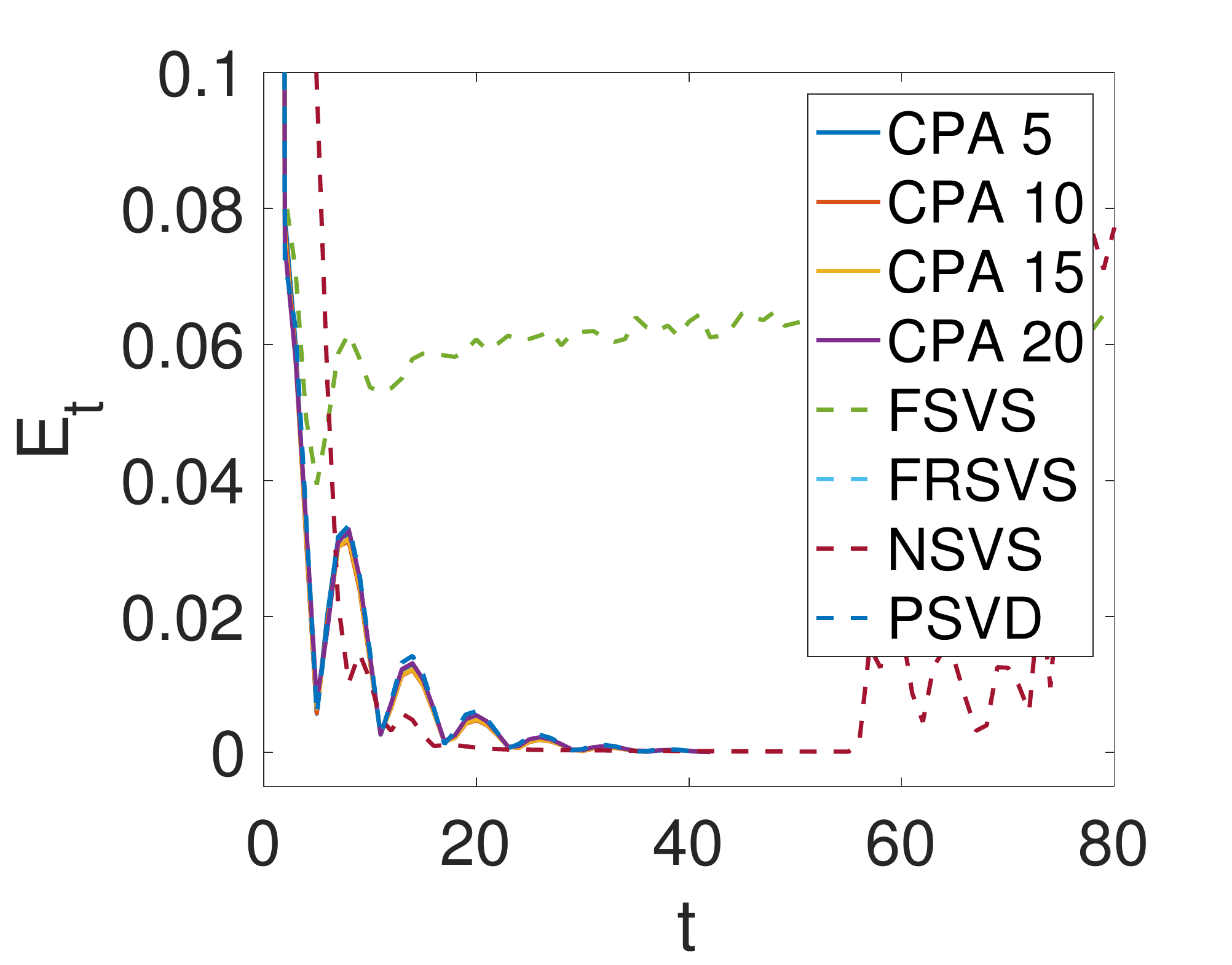}}
\end{minipage}
&
\begin{minipage}[h]{0.195\linewidth}
  \centering
  \centerline{\includegraphics[trim=3mm 3mm 3mm 3mm, clip,width=3.5cm]{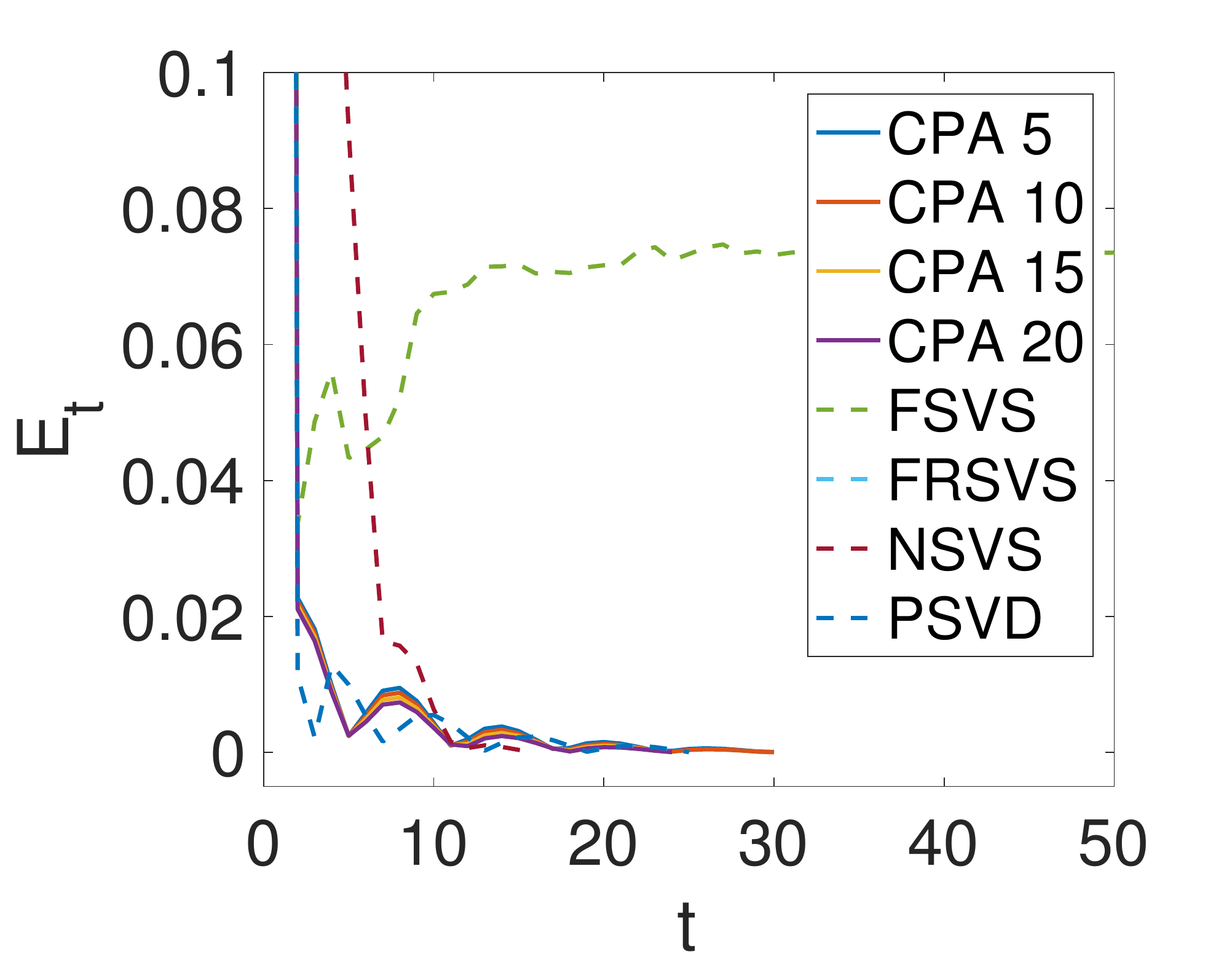}}
\end{minipage}
&
\begin{minipage}[h]{0.195 \linewidth}
  \centering
  \centerline{\includegraphics[trim=3mm 3mm 3mm 3mm, clip,width=3.5cm]{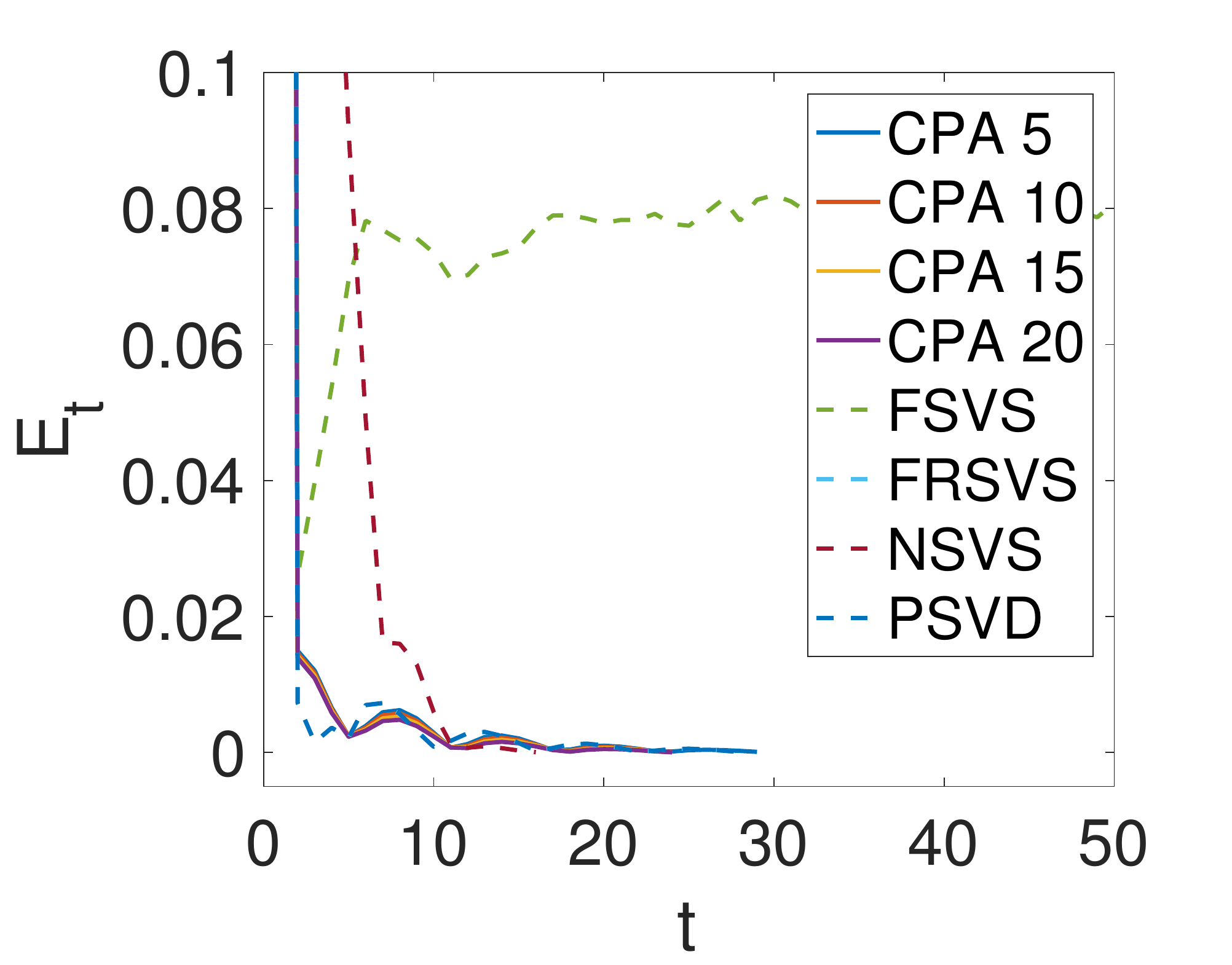}}
\end{minipage}
&
\begin{minipage}[h]{0.195 \linewidth}
  \centering
  \centerline{\includegraphics[trim=3mm 3mm 3mm 3mm, clip,width=3.5cm]{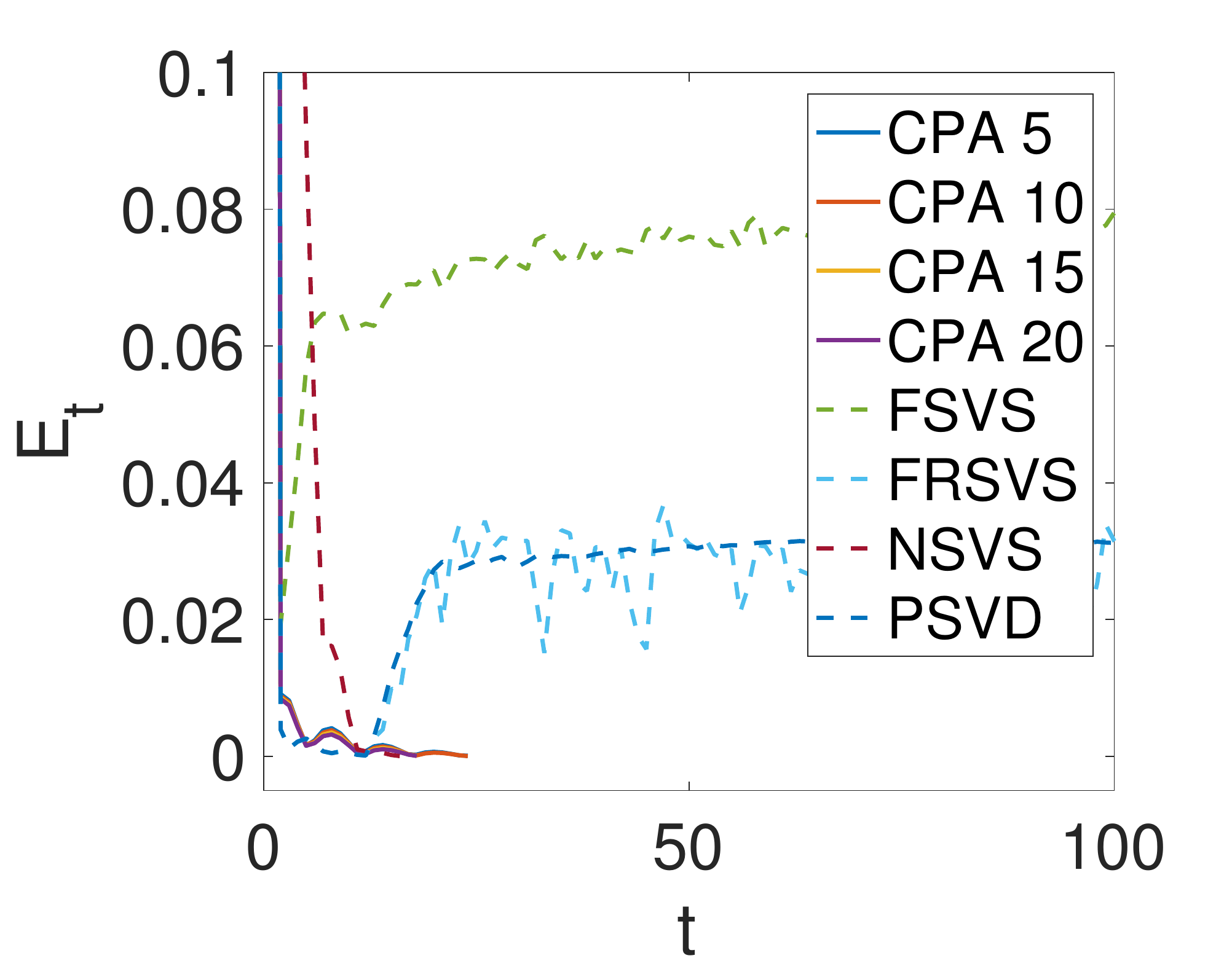}}
\end{minipage}
\\\hline
10$\%$&
\begin{minipage}[h]{0.195 \linewidth}
  \centering
  \centerline{\includegraphics[trim=3mm 3mm 3mm 3mm, clip,width=3.5cm]{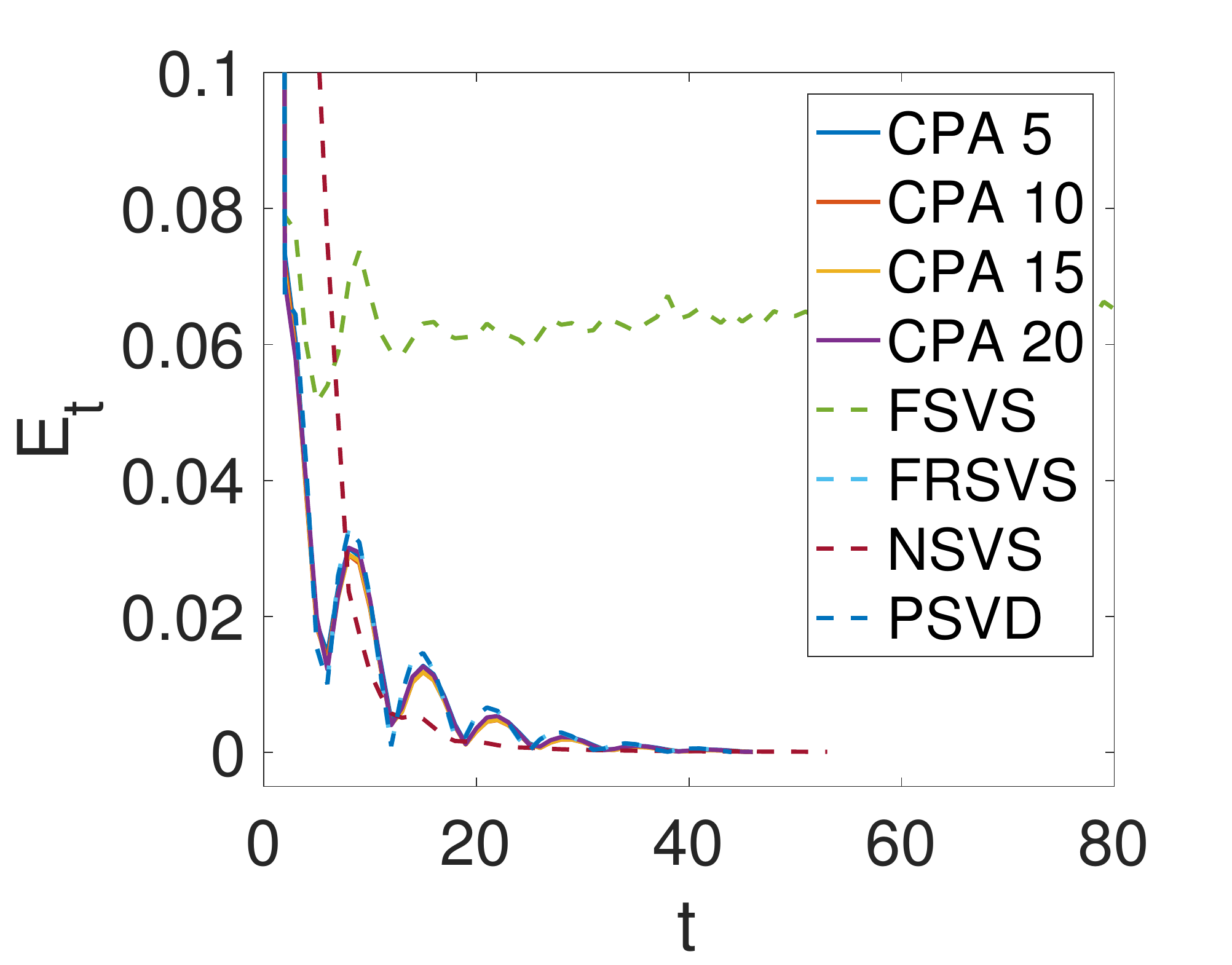}}
\end{minipage}
&
\begin{minipage}[h]{0.195 \linewidth}
  \centering
  \centerline{\includegraphics[trim=3mm 3mm 3mm 3mm, clip,width=3.5cm]{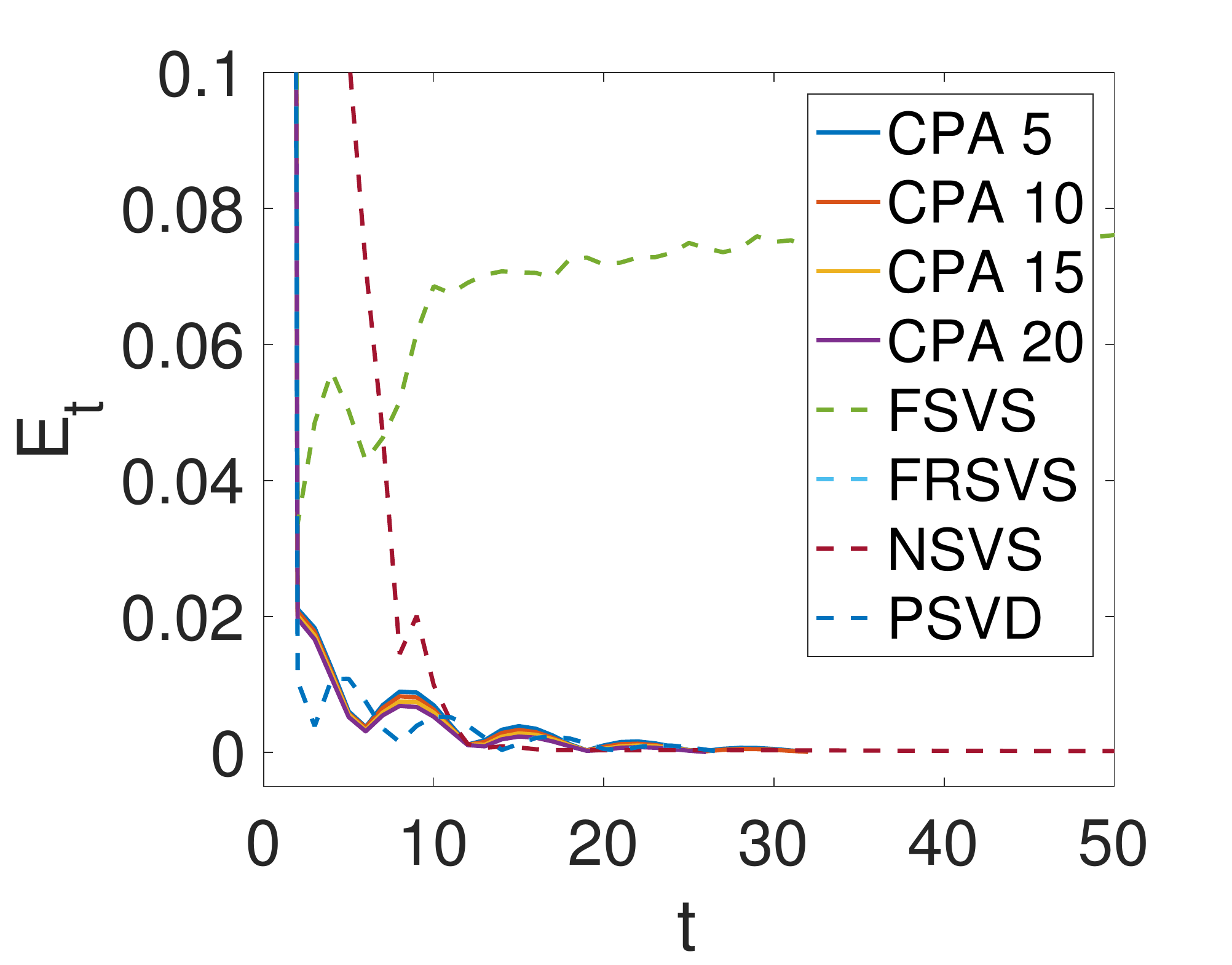}}
\end{minipage}
&
\begin{minipage}[h]{0.195 \linewidth}
  \centering
  \centerline{\includegraphics[trim=3mm 3mm 3mm 3mm, clip,width=3.5cm]{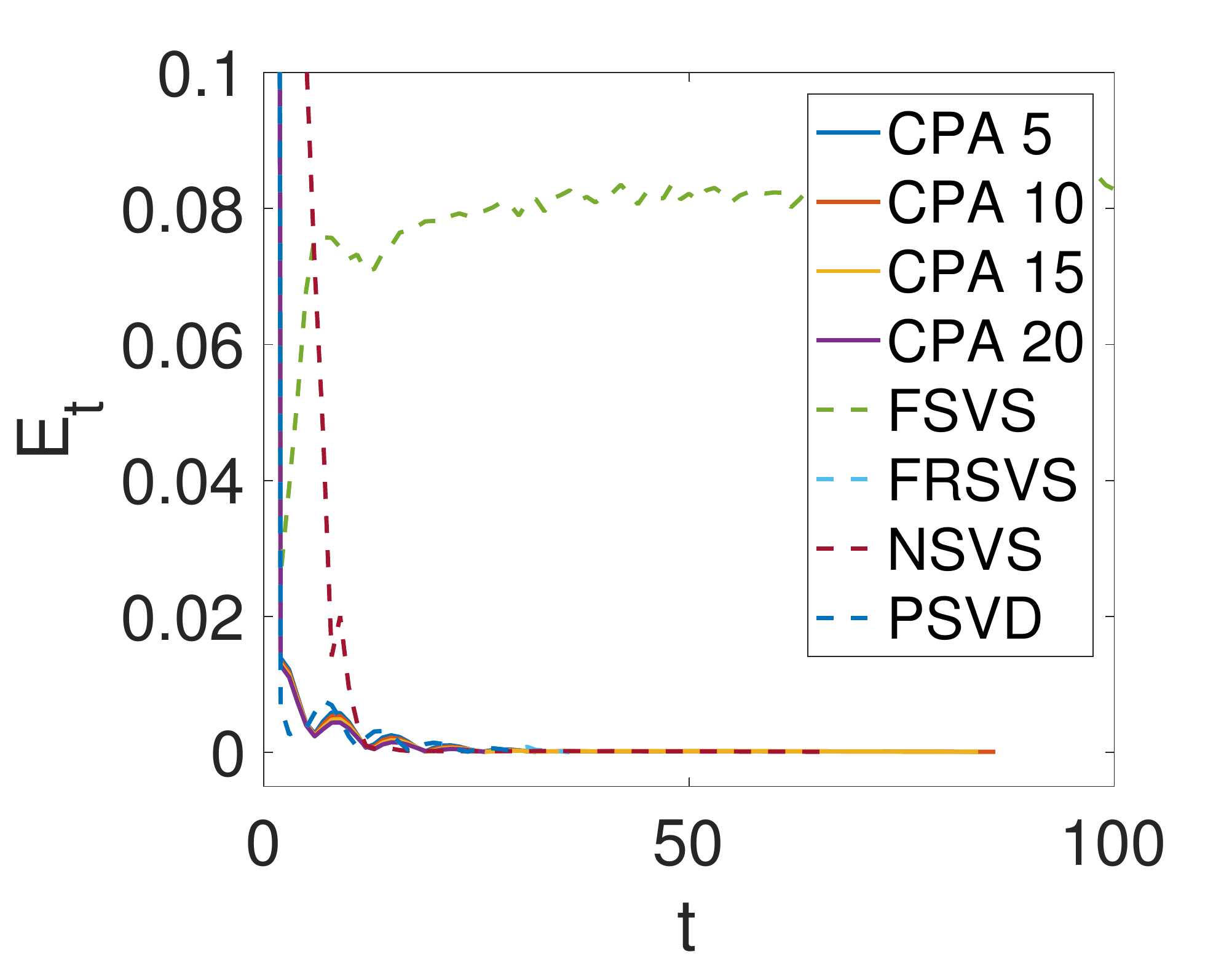}}
\end{minipage}
&
\begin{minipage}[h]{0.195 \linewidth}
  \centering
  \centerline{\includegraphics[trim=3mm 3mm 3mm 3mm, clip,width=3.5cm]{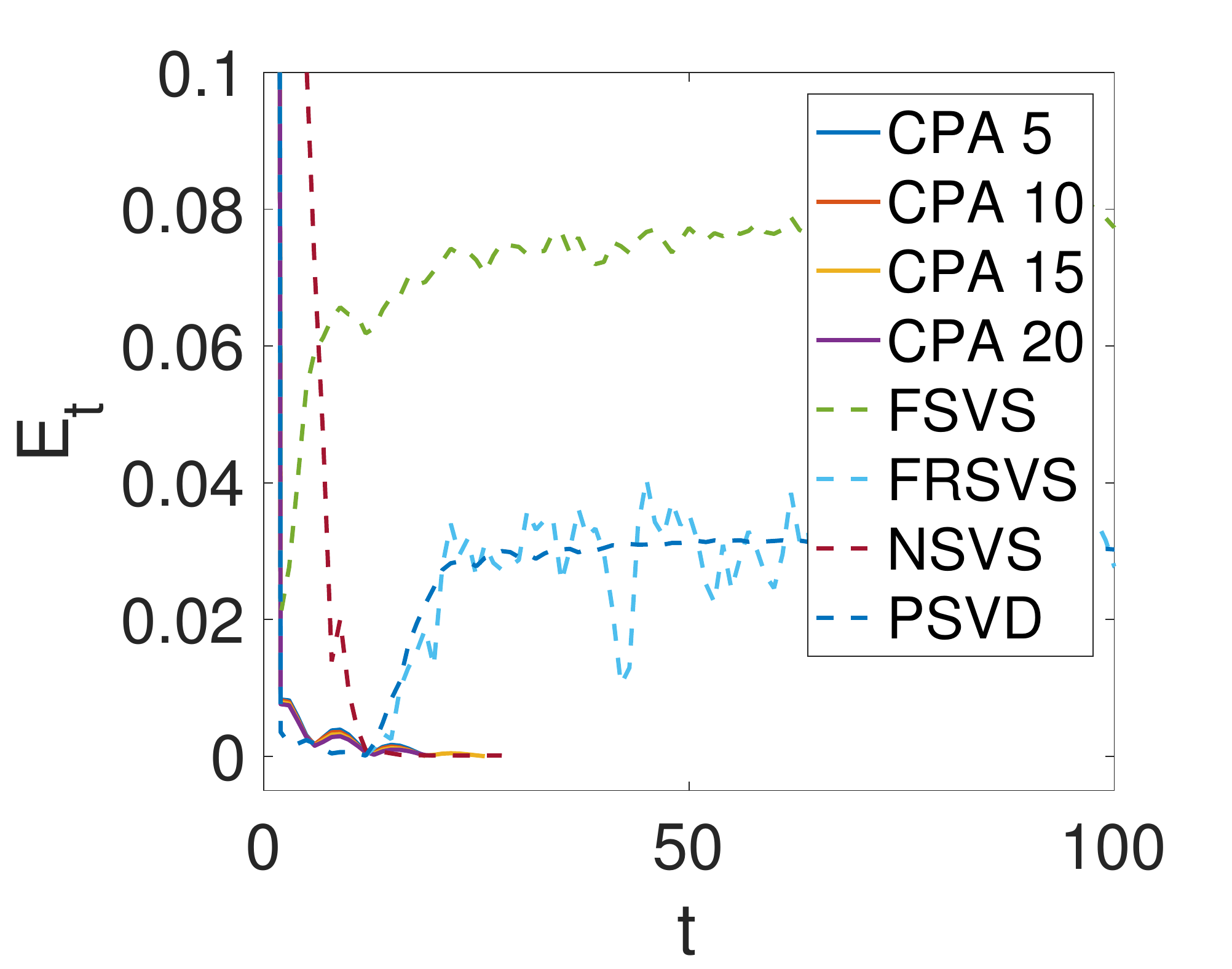}}
\end{minipage}
\\\hline
20$\%$&
\begin{minipage}[h]{0.195 \linewidth}
  \centering
  \centerline{\includegraphics[trim=3mm 3mm 3mm 3mm, clip,width=3.5cm]{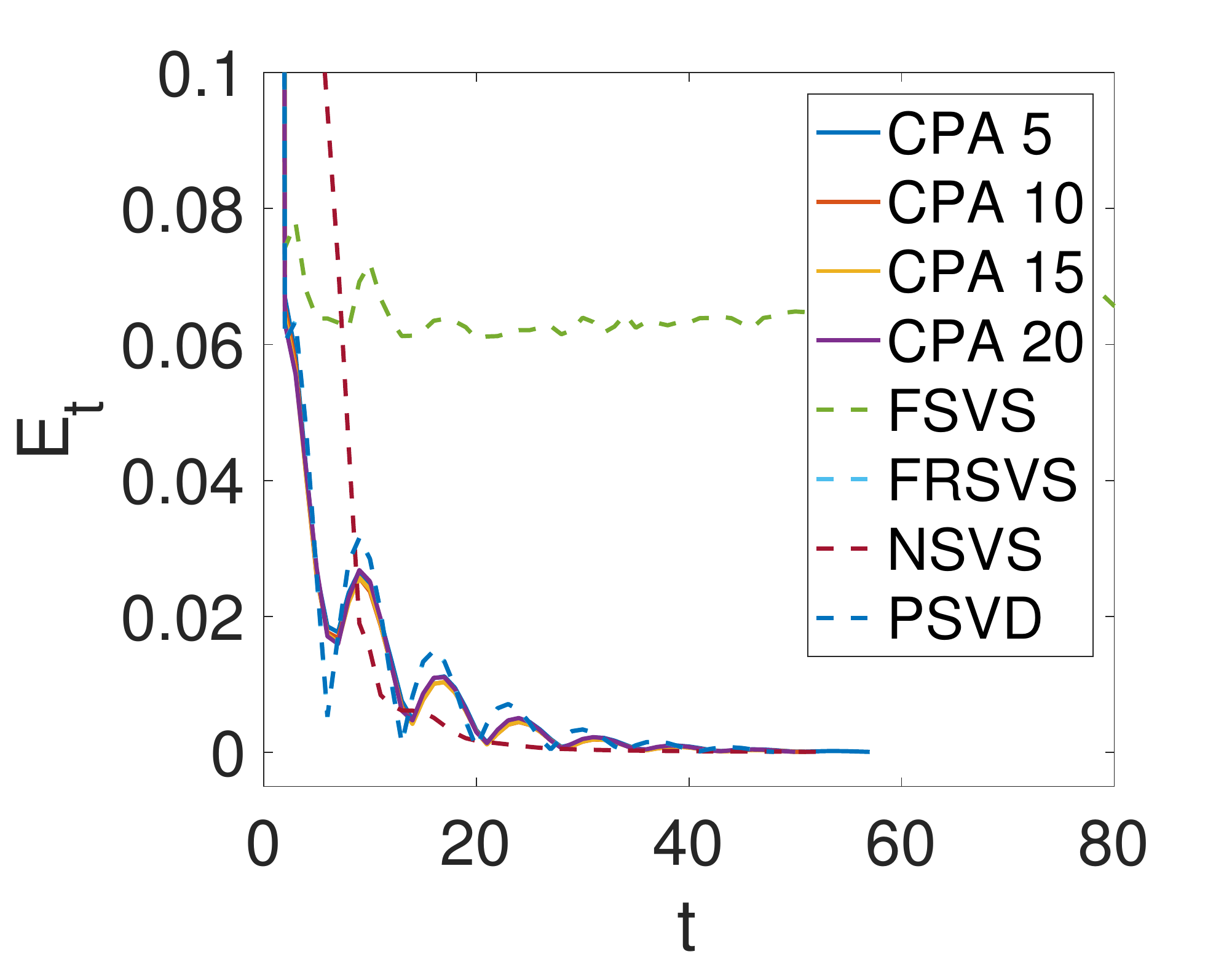}}
\end{minipage}
&
\begin{minipage}[h]{0.195 \linewidth}
  \centering
  \centerline{\includegraphics[trim=3mm 3mm 3mm 3mm, clip,width=3.5cm]{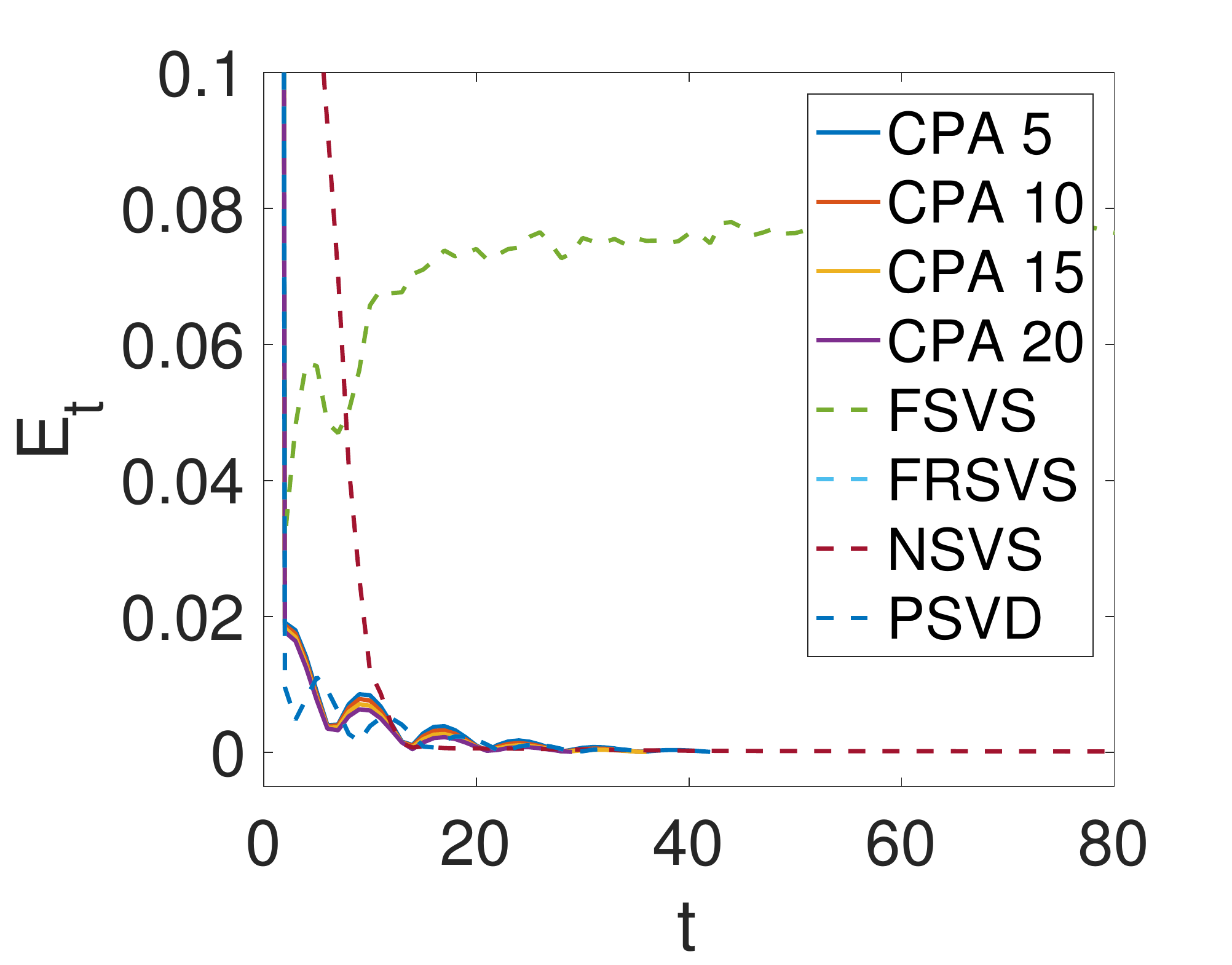}}
\end{minipage}
&
\begin{minipage}[h]{0.195 \linewidth}
  \centering
  \centerline{\includegraphics[trim=3mm 3mm 3mm 3mm, clip,width=3.5cm]{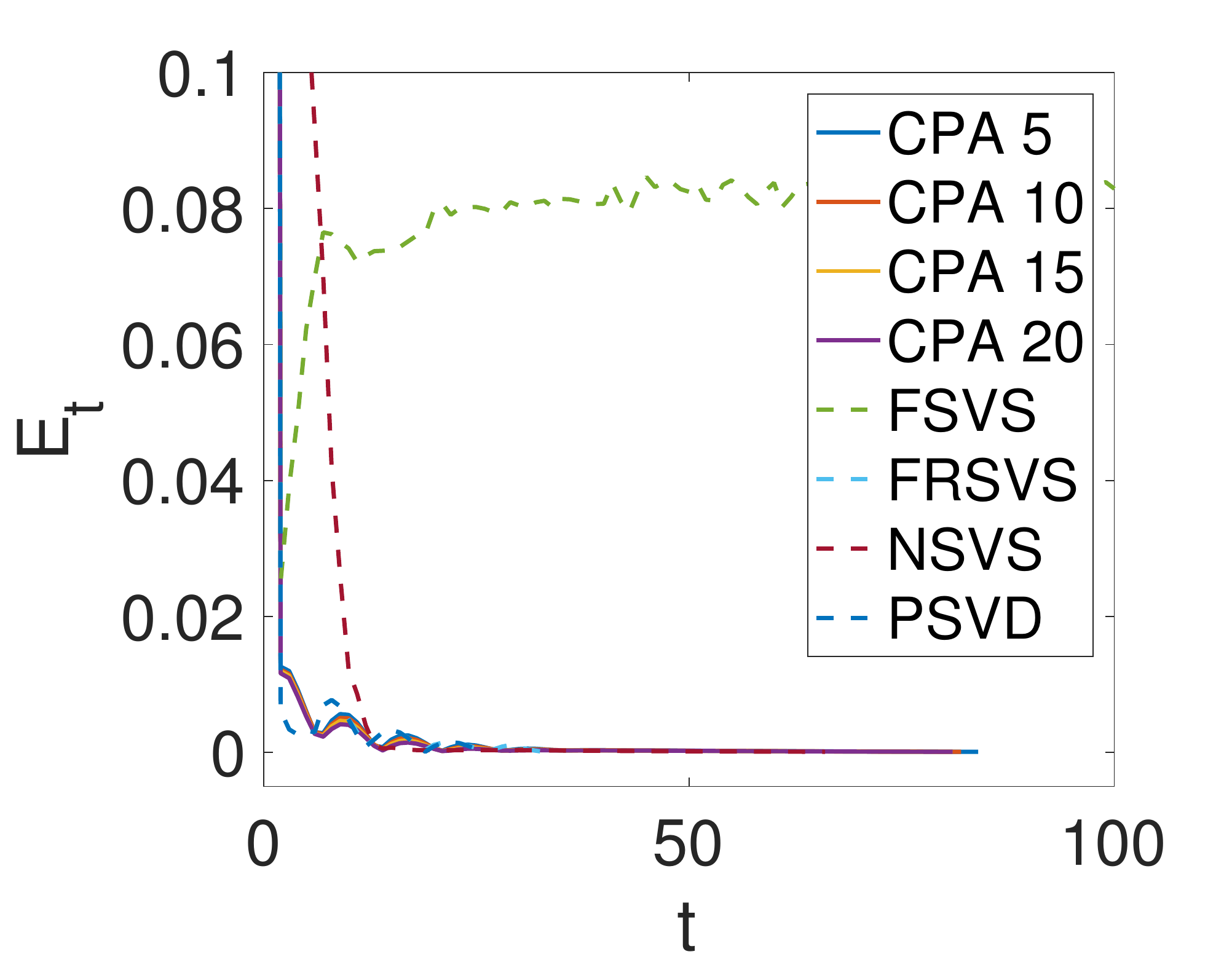}}
\end{minipage}
&
\begin{minipage}[h]{0.195 \linewidth}
  \centering
  \centerline{\includegraphics[trim=3mm 3mm 3mm 3mm, clip,width=3.5cm]{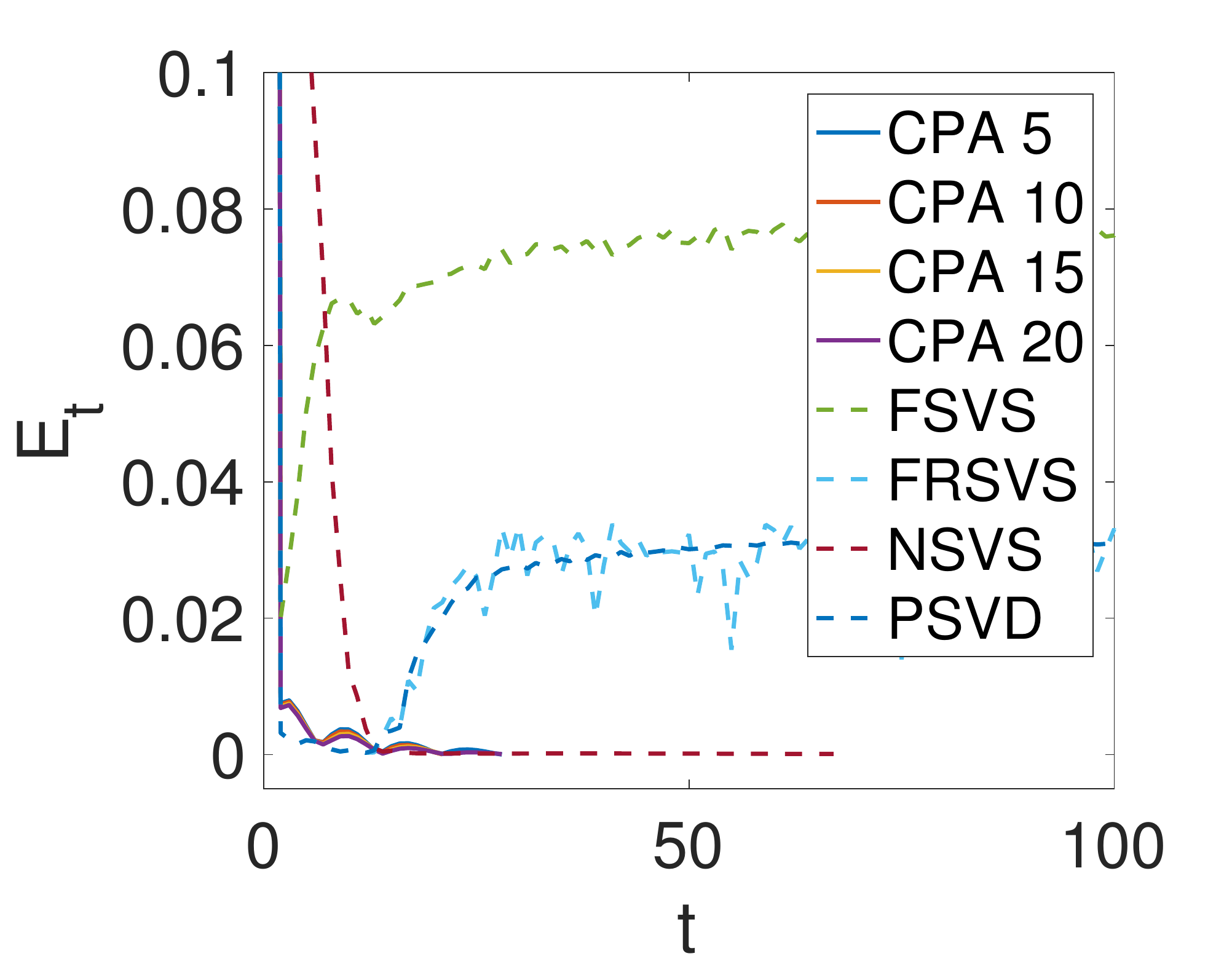}}
\end{minipage}
\\\hline
\end{tabular}
\caption{Figures represent $E_t$ for the systhetic data in the various conditions with respect to the reduction rates of elements and the matrix ranks.}
\label{fig:Comparison_results_Error_Syn}
\end{figure*}where $E_\ell \!:=\! 0.3\!\times\!10^{-1}$ and $E_\mathrm{m} \!:=\! 6\!\times\!10^{-4}$.
This was determined experimentally.
Recall that the size of $\underline{\mathbf{\Phi}}$ is $n\times n$.
To verify \eqref{recommendation_epsilon}, it was compared with fixed threshold.
Four values $k_1=(0.5\!\times\!10^{-4})n^2$, $k_2=(1.5\!\times\!10^{-4})n^2$, $k_3=(0.5\!\times\!10^{-3})n^2$, and $k_4 = n^2$ were used for this comparison, where $k_4$ means that all components of a matrix are retained.
Additionally, the DCT was exploited for the sparsifying matrix.

The results are listed in Table~\ref{computation_time_adaptive_DCT}.
The adaptive method was faster than the fixed method.
Our method with $k_1$ was the fastest but it did not converge well for the large size of approximation errors, where $k_1$ is considered as the limitation of the thresholding percentage in our method.
In addition, the singular value shrinkage of our CPA-based method by thresholding the matrix components was about five times faster than that by maintaining those.

%
%
\subsection{Comparison with Existing Methods}
\label{sec:comparison_existing_meth}

As previously mentioned, there are several fast singular value shrinkage methods \cite{SVD_echon2,SVD_echon3,Fast_singular_thresh}.
To illustrate the advantage of our method, we compared it with the singular value shrinkage by using the exact PSVD, fast randomized singular value shrinkage (FRSVS) \cite{SVD_echon2}, singular value shrinkage by using the Nystr\"om method (NSVS) \cite{SVD_echon3}, and the fast singular value shrinkage without the exact SVD (FSVS) \cite{Fast_singular_thresh}.
The experiments were conducted using the image inpainting method for \textit{Bricks} and a synthetic data.
The synthetic data is constructed as a block diagonal matrix whose number of the main diagonal blocks is equal to its rank.
Let $\mathbf{D} \in \mathbb{R}^{1000\times 1000}$ be the matrix form of the synthetic data with rank $n$, and this is defined as $\mathbf{D}:=\mathbf{J} - \mathrm{blkdiag}(\mathbf{D}_{\mathrm{s}},n)$, where $\mathbf{J}:=\mathbbm{1}_{1000}\mathbbm{1}^\top_{1000}$, $\mathbf{D}_{\mathrm{s}}:= 0.5\mathbbm{1}_{1000/n}\mathbbm{1}^\top_{1000/n}$ and $\mathrm{blkdiag}(\mathbf{D}_{\mathrm{s}},n):=\mathrm{diag}(\underbrace{\mathbf{D}_{\mathrm{s}},\mathbf{D}_{\mathrm{s}},\ldots,\mathbf{D}_{\mathrm{s}}}_{n})$.
In the experiment, $n=$10, 100, 200, and 500 were used, and $x\%$ of elements in $\mathbf{D}$ were randomly replaced with zero, where $x\in \{1,10,20\}$.
To restore the corrupted data, the optimization problem in \eqref{texture_repaire4} was solved without using the average term, i.e., $\mathrm{ave}\bigl( \mathrm{vec}(\mathbf{M}) \bigr) = \mathrm{ave}\bigl( \mathrm{vec}(\mathbf{M}_\partial) \bigr)$, because the matrices $\mathbf{M}$ and $\mathbf{M}_\partial$ can hardly be defined for the random corruption.
The function \textit{svdsechon}\footnote{Available at \url{https://www.mathworks.com/matlabcentral/fileexchange/47132-fast-svd-and-pca}} was used for carrying out the exact PSVD, and the 500 largest singular values were calculated for \textit{Bricks}.
All preferences of the FRSVS \cite{SVD_echon2} was determined in the original code\footnote{Available at \url{http://thohkaistackr.wixsite.com/page/projectfrsvt}} provided by the authors.
In the FRSVS, the 550 largest singular values were approximately derived for \textit{Bricks}, whose number of singular values was experimentally determined for carrying out precise and fast singular value shrinkage.
In addition, the 200 largest singular values were derived for the synthetic data in the exact PSVD and the FRSVS.
The partial singular values derived in the exact PSVD and the FRSVS were soft-thresholded: The $i$th partial singular value $\sigma^\mathrm{p}_i$ is shrunk to $\max(\sigma^\mathrm{p}_i-1/\rho,0)$.
Since the Nystr\"om method requires a square matrix to derive partial eigenvalues, it was applied to $\mathbf{X}^{\!\top}\!\mathbf{X}$, as in the EVD-based method.
The 500 and 200 largest eigenvalues were calculated for the \textit{Bricks} and the synthetic data, respectively.
Those eigenvalues were then shrunk in the same way as with the EVD-based method, whose number of calculated eigenvalues was experimentally determined from the same reason as the FRSVS.
All preferences used in the FSVS were directly used as suggested in \cite{Fast_singular_thresh}.

The results for \textit{Bricks} are indicated in Table~\ref{computation_time_dct_wav_blk} and Fig.~\ref{fig:Comparison_results_Error_RMSE}.
Note that the experiments of the existing methods were stopped at the 80th iterations since these methods did not converge.
The concept of the FSVS is similar to our method, but it requires longer computation time as shown in Table~\ref{computation_time_dct_wav_blk} and does not converge.
For the FRSVS and NSVS, although their average computation times are slightly less than ours, they result in much larger errors.
This would be because many singular values or eigenvalues above the threshold $1/\rho$ were reduced to zero, so that the exact PSVD, FRSVS, and NSVS produced the large errors in each iteration leading to unstable convergence.
In contrast, our method is stable and does not affect the convergence of the optimization method because the CPA-based method can shrink the entire singular values.

Figure~\ref{fig:Comparison_results_Error_Syn} shows the comparison of errors in the case of the synthetic data in several conditions according to the reduction rate of data elements and the matrix rank.
The optimization methods using the existing methods do not converge well when the matrix rank is 500.
This is because many singular values above $1/\rho$ are discarded, i.e., enormous computation errors are produced in each iteration.
From the results, the matrix rank of target data should be estimated beforehand, and then the numbers of partial singular values and vectors should be estimated to be larger than the matrix rank, in order to make the optimization method converged.
In contrast, our method can lead the optimization method to stable convergence.
It certainly generates some approximation errors, but it can process all singular values, which means that most singular values above $1/\rho$ are remained.

%
%
\section{Conclusion}
\label{sec:conclusion}

We proposed a fast thresholding method of singular values without computing singular values and vectors.
The key tool of the proposed method is CPA.
From CPA characteristics, singular value shrinkage could be computed by a multiplication of matrices.
The proposed method was further accelerated using the sparsity of a signal, where the frequency transform was used for obtaining sparse coefficients.
Moreover, we studied the approximation order for reducing the size of approximation errors.
The experimental results revealed that our method was much faster than the exact methods with high approximation precision in the case of a large data size.
In addition, our method can lead the optimization method to be stable convergence in comparison of the existing fast singular value shrinkage methods because of its approximation precision.

\appendices

%
%
\section{ADMM Applicable Forms}
\label{sec:admm_applicable_forms}

%
%
\subsection{Texture Image Inpainting}
\label{sec:admm_inpaint}

Let $\mathbf{i}\!:=\!\mathrm{vec}(\mathbf{I})$, $\mathbf{l}\!:=\!\mathrm{vec}(\mathbf{L})$, $\mathbf{m}\!:=\!\mathrm{vec}(\mathbf{M})$, and $\mathbf{m}_\partial\!:=\!\mathrm{vec}(\mathbf{M}_\partial)$.
The 2-D DCT matrix is represented as $\mathbf{\Psi}$, and the matrix form of $P_{\Omega}$ and $P_{\overline{\Omega}}$ are defined as $\mathbf{\Omega}$ and $\overline{\mathbf{\Omega}}$.
In addition, the indicator functions of the sets $\mathcal{I} \!:=\! \{\mathbf{x}\!\in\! \mathbb{R}^{mn}|~\mathbf{x}\!=\!\mathbf{\Omega}\mathbf{i} \}$, $\mathcal{M} \!:=\! \{\mathbf{x}\!\in\! \mathbb{R}^N\ |\ \mathrm{ave}(\mathbf{x})\!=\!\mathrm{ave}(\mathbf{m}_\partial)\}$ and $\mathcal{D}$ are denoted as $\iota_{\mathcal{I}}$, $\iota_\mathcal{M}$, and $\iota_\mathcal{D}$, respectively, where $N$ is the size of $\mathrm{vec}(\mathbf{M})$.
By using the above definitions, \eqref{texture_repaire4} is redefined as
\begin{equation}
\min_{\mathbf{l}} \| \mathbf{l}\|_{*}+\eta  \| \mathbf{\Psi}\mathbf{l}\|_1+\iota_{\mathcal{I}}(\mathbf{\Omega}\mathbf{l}) + \iota_\mathcal{D}(\mathbf{l}) + \iota_\mathcal{M}(\overline{\mathbf{\Omega}}\mathbf{l}).
\label{texture_repaire5}
\end{equation}
Let the vector $\mathbf{z} \!\in\! \mathbb{R}^{5mn}$ be
\begin{equation}
\mathbf{z}:=
\begin{bmatrix}
\mathbf{z}^{(1)}\\
\mathbf{z}^{(2)}\\
\mathbf{z}^{(3)}\\
\mathbf{z}^{(4)}\\
\mathbf{z}^{(5)}
\end{bmatrix}
=
\begin{bmatrix}
\mathbf{Id}\\
\mathbf{\Psi}\\
\mathbf{\Omega}\\
\mathbf{Id}\\
\overline{\mathbf{\Omega}}
\end{bmatrix}
\mathbf{l}
=
\mathbf{K}
\mathbf{l}
.
\label{z_insert}
\end{equation}
Finally, \eqref{texture_repaire5} is represented as
\begin{equation}
\scalebox{0.95}[0.95]{$
\begin{aligned}
\min_{\mathbf{l},\mathbf{z}} & \ \| \mathbf{z}^{(1)}\|_{*} \!+\! \eta  \| \mathbf{z}^{(2)}\|_1
 \!+\! \iota_{\mathcal{I}}(\mathbf{z}^{(3)}) \!+\! \iota_\mathcal{D}(\mathbf{z}^{(4)}) \!+\! \iota_\mathcal{M}(\mathbf{z}^{(5)}) \\
\text{s.t.} & \quad \mathbf{z}=\mathbf{K}\mathbf{l},
\end{aligned}
$}
\label{admm_eq}
\end{equation}
Equation \eqref{admm_eq} can be applied to the ADMM algorithm in \eqref{ADMM_algorithm} which is indicated in Appendix~\ref{sec:optimization_tool}.
Let  $\mathbf{u}_0 \!:=\! [(\mathbf{u}^{(1)}_0)^\top,(\mathbf{u}^{(2)}_0)^\top,(\mathbf{u}^{(3)}_0)^\top,(\mathbf{u}^{(4)}_0)^\top,(\mathbf{u}^{(5)}_0)^\top]^\top$ be an arbitrary auxiliary vector, where $\mathbf{u}^{(i)}_0 \!\in\! \mathbb{R}^{mn}$ in $i\!=\!1,\ldots,5$.
Applying the ADMM to \eqref{admm_eq} yields the following algorithm:
\begin{equation}
\left\lfloor{\begin{split}\mathbf{l}_{t+1} &:= (\mathbf{K}^\top \mathbf{K})^{-1}\mathbf{K}^\top(\mathbf{z}_t - \mathbf{u}_t)
\\
\mathbf{z}^{(1)}_{t+1} &:= \mathrm{prox}_{1/\rho \| \cdot \|_*} (\mathbf{l}_{t+1}+\mathbf{u}_t^{(1)})\text{~~------~($\ast$)}
\\
\mathbf{z}^{(2)}_{t+1} &:= \mathrm{prox}_{\eta / \rho \| \cdot \|_1} (\mathbf{\Psi}\mathbf{l}_{t+1}+\mathbf{u}_t^{(2)})
\\
\mathbf{z}^{(3)}_{t+1} &:= \Pi_{\mathcal{I}}(\mathbf{\Omega}\mathbf{l}_{t+1}+\mathbf{u}_t^{(3)})
\\
\mathbf{z}^{(4)}_{t+1}&:= \Pi_\mathcal{D}(\mathbf{l}_{t+1}+\mathbf{u}_t^{(4)})
\\
\mathbf{z}^{(5)}_{t+1} &:= \Pi_\mathcal{M}(\overline{\mathbf{\Omega}}\mathbf{l}_{t+1}+\mathbf{u}_t^{(5)})
\\
\mathbf{u}_{t+1}&:=\mathbf{u}_t+\mathbf{K}\mathbf{l}_{t+1} - \mathbf{z}_{t+1},\end{split}}\right.
\label{admm_algorithm_inp}
\end{equation}
where the update of $\mathrm{prox}_{1/\rho\| \cdot \|_*}$ in ($\ast$) of \eqref{admm_algorithm_inp} can be computed by singular value shrinkage, which is performed using our CPA-based method.
Additionally, the update of $\mathrm{prox}_{\eta/\rho\| \cdot \|_1}$ is approximated by soft-thresholding, i.e., $\mathrm{sgn}(X_{ij})\max(|X_{ij}|\!-\!\eta/\rho,0)$, where $X_{ij}$ is the entry of an arbitrary matrix $\mathbf{X}$ and $\mathrm{sgn}(\cdot)$ is the sign function.
In the following applications, the same calculation is used for the updates of $\mathrm{prox}_{1/\rho\| \cdot \|_*}$ and $\mathrm{prox}_{\eta/\rho\| \cdot \|_1}$.
In \eqref{admm_algorithm_inp}, $\Pi_{\mathcal{I}}(\cdot)$, $\Pi_\mathcal{D}(\cdot)$, and $\Pi_\mathcal{M}(\cdot)$ are the metric projections onto $\mathcal{I}$, $\mathcal{D}$, and $\mathcal{M}$, respectively.
Practically, $\Pi_{\mathcal{I}}(\cdot)$ is given by maintaining the assigned pixels, and $\Pi_\mathcal{D}(\cdot)$ is calculated by pushing the entries outside $[0,1]$ into $0$ or $1$ (the nearest is chosen).
Additionally, the auxiliary value is calculated as the difference between the average value of $\mathbf{m}_\partial$ on $\mathbf{i}$ and the average value of $\overline{\mathbf{\Omega}}\mathbf{l}_{t+1}+\mathbf{u}_t^{(5)}$ on the recovered region.
The $\Pi_\mathcal{M}(\cdot)$ is derived by adding the auxiliary value to $\overline{\mathbf{\Omega}}\mathbf{l}_{t+1}+\mathbf{u}_t^{(5)}$ on the recovered region.

%
%
\subsection{Background Modeling}
\label{sec:admm_background}

Let $\mathbf{i}\!:=\!\mathrm{vec}(\mathbf{I})$, $\mathbf{l}\!:=\!\mathrm{vec}(\mathbf{L})$, and $\mathbf{s}\!:=\!\mathrm{vec}(\mathbf{S})$.
The indicator function of the set $\mathcal{I} \!:=\! \{\mathbf{x}\!\in\! \mathbb{R}^{mnK}|~\mathbf{x} \!=\! \mathbf{i} \}$ is defined as $\iota_{\mathcal{I}}$.
By using the above definitions, \eqref{convex_opt_back_mod} is rewritten as
\begin{equation}
\min_{\mathbf{l}, \mathbf{s}} \, \| \mathbf{l} \|_* + \eta \| \mathbf{s} \|_1 + \iota_{\mathcal{I}} (\mathbf{l}+\mathbf{s}).
\label{convex_opt_back_mod2}
\end{equation}
When an auxiliary vector $\mathbf{z}$ is represented as
\begin{equation}
\mathbf{z}:=
\begin{bmatrix}
\mathbf{z}^{(1)}\\
\mathbf{z}^{(2)}\\
\mathbf{z}^{(3)}\\
\end{bmatrix}
=
\begin{bmatrix}
\mathbf{Id} & \mathbf{O}\\
\mathbf{O} & \mathbf{Id}\\
\mathbf{Id} & \mathbf{Id}
\end{bmatrix}
\begin{bmatrix}
\mathbf{l}\\
\mathbf{s}
\end{bmatrix}
=\mathbf{K}\mathbf{l}'.
\label{aux_back_mod}
\end{equation}
Problem \eqref{convex_opt_back_mod2} is further rewritten as
\begin{equation}
\min_{\mathbf{z}, \mathbf{l}'}
\ \| \mathbf{z}^{(1)} \|_* + \eta \| \mathbf{z}^{(2)} \|_1 + \iota_{\mathcal{I}}(\mathbf{z}^{(3)})
\quad
\text{s.t.} \ \mathbf{z}=\mathbf{K}\mathbf{l}'.
\label{convex_opt_back_mod3}
\end{equation}
Let $\mathbf{u}_0 \!:=\! [(\mathbf{u}^{(1)}_0)^\top,(\mathbf{u}^{(2)}_0)^\top,(\mathbf{u}^{(3)}_0)^\top]^\top$ be an auxiliary vector for the ADMM, where $\mathbf{u}^{(i)}_0 \!\in\! \mathbb{R}^{mnK}$ in $i\!=\!1,2,3$.
Applying the ADMM to \eqref{convex_opt_back_mod3} yields the following algorithm:
\begin{equation}
\left\lfloor{
\begin{aligned}
\mathbf{l}'_{t+1} &:= (\mathbf{K}^\top \mathbf{K})^{-1}\mathbf{K}^\top(\mathbf{z}_t - \mathbf{u}_t)
\\
\mathbf{z}^{(1)}_{t+1} &:= \mathrm{prox}_{1/\rho \| \cdot \|_*} (\mathbf{l}_{t+1}+\mathbf{u}_t^{(1)})\text{~~------~($\ast$)}
\\
\mathbf{z}^{(2)}_{t+1} &:= \mathrm{prox}_{\eta/\rho \| \cdot \|_1} (\mathbf{s}_{t+1}+\mathbf{u}_t^{(2)})
\\
\mathbf{z}^{(3)}_{t+1} &:= \Pi_{\mathcal{I}}(\mathbf{l}_{t+1}+\mathbf{s}_{t+1}+\mathbf{u}_t^{(3)})
\\
\mathbf{u}_{t+1}&:=\mathbf{u}_t+\mathbf{K}\mathbf{l}'_{t+1} - \mathbf{z}_{t+1},
\end{aligned}
}\right.
\label{admm_algorithm_mod}
\end{equation}
where the update of $\mathbf{z}^{(1)}$ in ($\ast$) of \eqref{admm_algorithm_mod} is calculated using our CPA-based method.
In \eqref{admm_algorithm_mod}, $\Pi_{\mathcal{I}}(\cdot)$ is the metric projection onto $\mathcal{I}$, which is given by maintaining the observed pixel values of the original sequences.

%
%
\bibliographystyle{IEEEbib}
\bibliography{myrefs}

\end{document}